\documentclass[tbtags,reqno]{amsart}

\usepackage{mathrsfs}
\usepackage{latexsym}
\usepackage{amsthm,amsopn,amsfonts,amssymb,epsfig,tabmac}
\numberwithin{equation}{section}

\makeatletter
\renewcommand{\subsubsection}{\@startsection
{subsubsection}
{3}
{0mm}
{\baselineskip}
{-0.5\baselineskip}
{\normalfont\normalsize\bfseries}}
\makeatother

\newtheorem{theorem}{Theorem}
\newtheorem{lemma}[theorem]{Lemma}
\newtheorem{proposition}[theorem]{Proposition}

\newtheorem{conjecture}[theorem]{Conjecture}

\newtheorem{corollary}[theorem]{Corollary}

\newtheorem{definition}[theorem]{Definition}

\theoremstyle{remark}
\newtheorem{remark}[theorem]{Remark}
\newtheorem*{acknow}{Acknowledgments}

\def\la{{\lambda}}

\def\cal L{{\mathcal L}}

\newcommand{\cercle}[1]{\ensuremath{\setlength{\unitlength}{1ex}\begin{picture}(2.8,2.8)\put(1.4,0.7){\circle{2.8}\makebox(-5.6,0){#1}}\end{picture}}}
\newcommand{\tcercle}[1]{\ensuremath{\setlength{\unitlength}{1ex}\begin{picture}(2.8,2.8)\put(1.4,1.4){\circle{2.8}\makebox(-5.6,0){#1}}\end{picture}}}
\newcommand{\gcercle}{\ensuremath{\setlength{\unitlength}{1ex}\begin{picture}(5,5)\put(2.5,2.5){\circle{5}}\end{picture}}}

\newcommand{\ds}{\displaystyle}
\newcommand{\mc}{\mathcal}

\newcommand{\spar}{\ensuremath{\mbox{\textnormal{SPar}}}}
\newcommand{\versg}[1]{\ensuremath{\overleftarrow{\,\phantom{|}#1\,\phantom{|}}}}
\newcommand{\versd}[1]{\ensuremath{\overrightarrow{\,\phantom{|}#1\,\phantom{|}}}}
\newcommand{\aand}{\ensuremath{ \quad\textrm{and}\quad}}
\newcommand{\wwhere}{\ensuremath{ \quad\textrm{where}\quad}}
\let\d\partial
\let\n\noindent

\def\T{{\mathcal T}}
\def\M{{\mathcal M}}

\def\ee{{\mathbf e}}
\let\la\lambda
\let\La\Lambda
\let\Om\Omega
\let\om\omega
\let\ta\theta
\let\ka\kappa
\let\rw\rightarrow

\newcommand{\LL}{\ensuremath{\langle\!\langle}}
\newcommand{\RR}{\ensuremath{\rangle\!\rangle}}

\begin{document}

\title[Symmetric functions in superspace]
{Symmetric functions in superspace}

\author{P. Desrosiers} \thanks{pdesrosi@phy.ulaval.ca}
\address{D\'epartement de physique, de g\'enie physique et
d'optique, Universit\'e Laval, Qu\'ebec, Canada, G1K 7P4.}
\author{L. Lapointe}
\thanks{lapointe@inst-mat.utalca.cl }
\address{Instituto de Matem\'atica y F\'{\i}sica, Universidad de
Talca, Casilla 747, Talca, Chile.}
\author{P. Mathieu} \thanks{pmathieu@phy.ulaval.ca}
\address{D\'epartement de physique, de g\'enie physique et
d'optique, Universit\'e Laval,  Qu\'ebec, Canada, G1K 7P4.}

\maketitle

\begin{quote}
\small{\sc Abstract.}
We construct a generalization of the theory of symmetric
functions involving functions  of commuting and anticommuting
(Grassmannian) variables.
These new functions,
called symmetric functions in superspace, are invariant under the
diagonal action of
the symmetric group acting on the sets of commuting and anticommuting
variables.  We first obtain superspace analogues of
a  number of standard objects and concepts in
the theory of symmetric functions:
partitions, monomials, elementary
symmetric functions, completely symmetric functions, power sums, involutions,
generating functions, Cauchy formulas, and scalar products.  We then
consider a one-parameter  extension of the combinatorial scalar product.
It provides
the natural setting for the definition of a family of
``combinatorial'' orthogonal Jack
polynomials in superspace.  We show that  this family coincides with
that of  ``physical'' Jack polynomials in
superspace that were previously introduced by the authors
as orthogonal eigenfunctions of a supersymmetric quantum mechanical many-body problem.
The equivalence of the two families is established by showing that the
``physical'' Jack polynomials are also orthogonal with respect to the
combinatorial scalar product.
This equivalence is also directly
demonstrated for particular values
of the free parameter.
\end{quote}

\tableofcontents

\section{Introduction}

Grassmannian variables refer to anticommuting variables $\{\ta_i\}$,
that is, to variables obeying the relation
\begin{equation}\label{grassmannvar}
\ta_i\ta_j=-\ta_j\ta_i\, ,
\end{equation}
and in particular,
\begin{equation}
\ta_i^2=0\, .
\end{equation}
In the first two subsections of this introduction, we have tried to
review the origin of Grassmannian variables in mathematics and
physics. Aside from its intrinsic interest, our aim in doing so is
to show that our approach to superpolynomials is well
established within the conceptual framework of supersymmetry.
However, these first two subsections can be safely skipped. For the
interested readers, we  stress that no background in quantum field
theory is assumed;  we have simply tried to  give a flavor of the
underlying physics through simple illustrations.  The
introduction pertaining
to the present work starts in subsection 1.3.

\subsection{Grassmannian variables in physics and mathematics}

As suggested by their name, the introduction of Grassmannian
variables goes back to Grassmann in the framework of his theory of
extension (1844), an ancestor of vector analysis (see e.g.
\cite{crowe}, chapter 3).\footnote{This subsection on the emerging
role of Grassmannian variables in mathematics and physics is based
in part on \cite{rosa}.} More precisely, Grassmann introduced  in
this context basis elements $\{e_x,e_y,e_z\}$ and a new outer
product satisfying
\begin{equation}
e_x e_y = -e_y e_x\, , \quad e_x e_z = -e_z e_x \, ,\quad e_y e_z =
-e_z e_y \, , \qquad e_x e_x= e_y e_y = e_z e_z=0 \,  .
\end{equation}
In the modern reinterpretation of this construction, $e_x, e_y,e_z$
are replaced by $dx,dy,dz$, and  the product is antisymmetrized (it
corresponds to the wedge product); the anticommuting character of
the variables is thus recaptured by a new product structure.

Anticommuting quantities {\it per se} have been rediscovered in
physics. In 1928, Jordan and Wigner realized that in order to
implement the Pauli exclusion principle (i.e., two electrons cannot
share the same set of quantum numbers), fermionic fields had to be
quantized with modes subject  to anticommutation relations instead
of the usual commutation relations \cite{jor}. This amounted to
introduce operators $b_n$ (where $n$ is an integer labeling a normal
mode) subject to the following rules
\begin{equation}
\{b_n,b_m\}:=b_nb_m+b_mb_n= \delta_{n+m,0} \, ,
\end{equation}
the exclusion principle being then a consequence of $b_n^2=0$ for
$n\not=0$. (By contrast, the usual commutation relations would take
the form $[a_n,a_m]= \delta_{n+m,0} $ where $a_n$ is an (ordinary)
operator, i.e., a harmonic oscillator mode). Within this scheme, the
new anticommuting quantities that were introduced were operators and
not variables.

The early development of both quantum theory and quantum field
theory relied entirely on the canonical quantization method, i.e.,
the lifting of a classical structure at the operator level and the
decomposition of
the operators in modes subject to a simple commutation or
anticommutation relation according to their statistics (that of
bosons or fermions respectively). The advent of Feynman path
integration as an alternative method of quantization (a method, we
stress, that involves variables instead of operators) posed the
problem of integrating Grassmannian variables. It was first pointed
out in \cite{kalat} that the analog of the (bosonic)
multidimensional Gaussian integral
\begin{equation}
I_b=\int \left(\prod_i dx_i\right) \, e^{\sum_{i,j} x_i b_{ij} x_j}
\propto \, ({\rm det} \, b)^{-1/2}
\end{equation}
(with integration over all space), would have to take the following
form for Grassmannian variables
\begin{equation} \label{bere}
I_f=\int \left(\prod_i d\ta_i\right)\, e^{\sum_{i,j} \ta_i f_{ij}
\ta_j} \propto \, ({\rm det} \, f)^{1/2}
\end{equation}
This prompted the contribution of  Berezin to the theory of
Grassmannian algebra (whose results are summarized in \cite{ber}),
and in particular, his proposal for the basic integral relations:
\begin{equation}
\int d\ta = 0\, , \qquad \int d \ta \, \ta = 1\, ,
\end{equation}
which ensured (\ref{bere}). With (left) differentiation defined in
the natural way,
\begin{equation}
\frac{\d}{\d\ta_i} \ta_j\ta_k= \delta_{ij}\ta_k- \delta_{ik}\ta_j\,
,
\end{equation}
it is readily seen that the Berezin integration of a Grassmannian
variable is essentially equivalent to its differentiation:
\begin{equation}
\int d\ta \, [ \alpha(x)+\ta \, \beta(x)]= \beta(x)= \frac{\d}{ \d
\ta} [\alpha(x)+\ta \beta(x)]
\end{equation}
This work on Grassmannian algebras is at the root of the Lie algebra
extensions with fermionic generators (see e.g., the review
\cite{CNS}).

Quite interestingly, at about the same time, another
     graded algebra arose in physics in a totally different setting.
The superconformal (or Ramond-Neveu-Schwarz) algebra \cite{RNS} is a
graded version of the Virasoro algebra that appears  when fermionic
degrees of freedom are inserted in the dual model (at the time, an
alternative to quantum field theory and subsequently understood as a
string theory). Gervais and  Sakita showed that this graded
superconformal algebra comes from a sort of supergauge
transformation, a transformation involving  anticommuting parameters
and which, in retrospect, precisely reflects the supersymmetric
invariance of the dual model on the two-dimensional world-sheet
\cite{GS}.

\subsection{Supersymmetry, superfields and superspace}

Supersymmetry is certainly one of the most spectacular and profound
ideas that has emerged from theoretical physics over the last thirty
years \cite{SUSY}. This is a symmetry that relates bosons and
fermions.  And in the context of a quantum field theory, it
corresponds to a fermionic symmetry that changes the statistics of
the fields. Schematically, and for a one-dimensional space, this
transformation takes the form\footnote{This is a sample
transformation. One could have chosen the analogous transformation
in which  the derivative acts on $F(x)$  instead of $B(x)$.  The
proper choice is dictated by the relative dimension of the bosonic
and fermionic fields.}
\begin{equation} \label{susy}
\delta B(x)= \eta F(x)\, ,\qquad  \delta F(x)= \eta \d_x B(x)\, ,
\end{equation}
     where $\eta$ is an anticommuting constant
$(\eta^2=0$) and where $B$ and $F$ are respectively bosonic and
fermionic fields.

The discovery of supersymmetry within the context of
four-dimensional quantum field theory, at a time when confidence in
quantum field theory had been resurrected after a period of doubts,
created a highly favorable situation for the rapid expansion of this
area. To illustrate the depth and significance of these early
developments, it suffices to mention the astonishing observation
that whenever the supersymmetric transformation is  localized (i.e.,
the parameter $\eta$ is no longer regarded as a constant),
supergravity emerges automatically.  Gravity can thus be viewed as a
consequence of supersymmetry \cite{Free}.

In the study of supersymmetric quantum field theories, an important
technical tool was introduced by Salam and Strathdee \cite{SS}: the
concept of superspace.\footnote{It would be fair to indicate that
superspace was actually  first discovered in the context of dual
models in \cite{Mon}. Note that in one of the pioneering papers on
supersymmetry (the second reference in \cite{SUSY}), the authors
constructed a fermionic field theory whose  translation invariance
was enlarged to accommodate translations of the form
$x_\nu\rightarrow x_\nu +\eta \psi_\nu $
     (where $\nu$ is a space-time index and $\psi_\nu$ is a fermionic
     field).  In a sense, they investigated a field theory defined in a
larger space: to the
     usual variables $x_\nu$ describing the four dimensional Minkowski space,
     they added new fermionic coordinates $\psi_\nu(x)$. The difference was
that these $\psi_\nu(x)$
     were the very dynamical fields in terms of which the theory was
     defined.} In its simplest setting (still keeping the
illustrative formulas  within the context of a one-dimensional
space), a bosonic field $B(x)$ and its fermionic partner $F(x)$ are
collected together within a superfield $\Phi(x,\ta)$, regarded as a
function of $x$ and a new anticommuting space variable $\ta$. Since
$\ta^2=0$, the Taylor expansion of the superfield in the $\ta$
variable contains only two terms, the two ``component fields''
\begin{equation}
\Phi(x,\ta)= F(x)+\ta B(x)\, .
\end{equation}
The supersymmetric transformation can now be interpreted
geometrically within superspace, the space described by the doublet
($x, \, \ta$), as a simple translation of the form
\begin{equation} x\rw x-\eta \ta \; ,\qquad \ta \rw\ta+\eta\, .
\end{equation}
With \begin{equation} \delta\Phi= \Phi(x-\eta \ta,
\ta+\eta)-\Phi(x,\ta)= \ta \eta \d_x F(x) +\eta B(x)\,
,\end{equation} and by setting
\begin{equation}
\delta\Phi = \delta F(x) + \ta \delta B(x)
\end{equation}
     we recover (\ref{susy}).\footnote{Observe that $\delta$ commutes
with $\ta$ since it is bosonic: $\Phi$ and $\delta \Phi$ have the
same statistics. For our illustrative example, we chose to construct
a fermionic superfield. Another option would have been to set
${\tilde \Phi}= B+\ta F$, which is bosonic (the product of
two fermionic variables is bosonic, i.e., $\ta F$ is bosonic like
$B$). This would have led to modified transformation rules: $\delta
B(x)= \eta \d_x F(x)$ and $ \delta F(x)= \eta B(x)$.}  Note that if
we denote the supersymmetry transformation generator by ${\mathcal Q}$,
i.e., if we set $\delta\Phi= \eta \{{\mathcal Q},\Phi\}$, we easily
obtain the superspace differential realization ${\mathcal Q}= \d_\ta-\ta
\d_x$.

Superfields (or, for the present matter, superfunctions) and
superspace are the concepts we wanted to introduce before
formulating the objectives of the present work.

\subsection{Symmetric superpolynomials}

Our program here is to lay down the foundation of the theory of
symmetric polynomials in superspace. The superspace we are
interested in is the superextension of the Euclidean space in $N$
variables, which we shall denote $\mathscr{E}^{N|N}$  and whose
coordinates will be $(x_1, \cdots x_N; \ta_1, \cdots \ta_N)
$, with $x_ix_j=x_jx_i, \, x_i\ta_j=\ta_jx_i$ and $\ta_i\ta_j=-\ta_j\ta_i$.~\footnote{Superspaces with much more complicated structures can
be constructed. For instance, the theory of supermanifolds is
treated in \cite{bere,DeWitt}. To be precise, it should be said that
supermanifolds have been discovered in \cite{Leites} before the
introduction of superspace and actually even before the seminal
Wess-Zumino paper \cite{SUSY}.} Superfunctions, or functions in
superspace, are thus functions of two types of variables. For
instance, all superfunctions in $\mathscr{E}^{2|2}$ are combinations
of the following expressions
\begin{equation}
f_0(x_1,x_2)\, , \qquad \ta_1 f_1(x_1,x_2) + \ta_2 f_2(x_1,x_2)\, ,
\qquad \ta_1\ta_2 f_{3}(x_1,x_2)
\end{equation}
where the $f_i$'s stand for arbitrary functions of $x_1$ and $x_2$.
Superfunctions of the second type are fermionic (alternatively said
to be odd) while those of the first and third types are bosonic
(even).

The immediate question we have to address is the following: what is
the meaning of a symmetric superfunction? Observe that there are two
copies of the symmetric group at our disposal: the usual one, acting
on the commuting variables $x_i$, spanned by the exchange operators
$K_{ij}$ defined by
\begin{equation}
K_{ij}x_j = x_i K_{ij}
\end{equation}
     and another one acting on the anticommuting variables $\ta_i$, generated
by the new exchange operators $\ka_{ij}$:
\begin{equation}
\ka_{ij}\ta_j = \ta_i \ka_{ij}
\end{equation}
We stress that a symmetric superfunction is not a function invariant
under each type of symmetry transformation.\footnote{There is only
one symmetric function of $\theta$, which is
$\theta_1+\ldots+\theta_N$.} It is a function invariant under the
action of the diagonal subgroup of the tensor product of these two
copies of the symmetric group, i.e., superfunctions are invariant
under the simultaneous interchange of $(x_i,\ta_i)$ and
$(x_j,\ta_j)$.  In other words, a symmetric superfunction $f$
satisfies the condition
\begin{equation}
K_{ij}\ka_{ij} f = f
\end{equation}
Examples of symmetric superpolynomials in $\mathscr{E}^{2|2}$ are
\begin{equation}
x_1^2x_2^2\, , \quad\ta_1 x_1^4+\ta_2x_2^4\, , \quad \ta_1x_2^2+
\ta_2x_1^2\,, \quad  \ta_1\ta_2(x_1^3x_2-x_1x_2^3)\,.
\end{equation}

The enforced interconnection between the transformation properties
of the bosonic and the fermionic variables is a direct consequence
of the definition of a supersymmetric transformation as a
translation in superspace.  This is also what makes the resulting
object most interesting and novel. In particular, it ensures that
the resulting symmetric superpolynomials are completely different
from the ``supersymmetric polynomials'' previously considered in the
literature.

Recall that what is  called a supersymmetric polynomial (see e.g.,
\cite{Stem}) is first of all a doubly symmetric polynomial in two
distinct sets of ordinary (commuting) variables $x_1,\cdots x_m$ and
$y_1,\cdots, y_n$, i.e.,  invariant under
independent permutations of the $x_i$'s and the $y_i$'s. It is said
to be supersymmetric if, in addition, it satisfies the following
cancelation condition: by substituting $x_1=t$ and $y_1=t$, the
polynomial becomes independent of $t$. An example of a generating
function for such polynomials is
\begin{equation}
\prod_{i=1}^m (1-qx_i)\prod_{j=1}^n(1-qy_j)^{-1}= \sum_{r\geq 0}
p_{(r)}(x,y) q^r
\end{equation}
This generating function is known to appear in the context of
classical Lie superalgebras (as a superdeterminant) \cite{Kac}.
Actually most of the work on supersymmetric polynomials is motivated
by its connection with superalgebras. For an example of such an
early work, see \cite{jarvis}. More precisions and references are
also available in \cite{brenti, moens}.

The key differences between  these supersymmetric polynomials and
our symmetric superpolynomials  should
     be clear. In our case, we symmetrize two sets of variables with
respect to the diagonal action of the symmetric group.  And
moreover, one of the two sets is made out of Grassmannian variables.

\subsection{Toward a theory of symmetric polynomials in superspace}

The first step in the elaboration of a theory of symmetric
polynomials in superspace is the introduction of a proper labeling
for bases of the ring of symmetric superpolynomials, that is, a
superversion of partitions. With this concept in hand, the
construction of the superextension of the symmetric monomial basis
(supermonomial basis for short) is rather immediate.
  From there on, there are two natural routes that can be followed.

\subsubsection{Orthogonal symmetric superpolynomials from Cauchy superformulas}

     The first approach amounts
to extend to superspace the other classical symmetric functions.
This could be done via the extension of their generating functions
which, for the elementary $e_n$, homogeneous $h_n$ and power sum
$p_n$ symmetric functions are respectively given by \cite{Mac}:
\begin{equation}
\sum_{n\geq 0} e_n\, t^n=\prod_{i\geq 1}(1+x_it)\; ,\quad
\sum_{n\geq 0} h_n\, t^n=\prod_{i\geq 1}\frac{1}{(1-x_it)}\;,\quad
\sum_{n\geq 1} p_n\, t^n=\prod_{i\geq 1}\frac{x_i t}{(1-x_it) }\, .
\end{equation}
The basis elements are generated from the product of these
functions:
\begin{equation}
f_\la = f_{\la_1}\cdots f_{\la_n}
\end{equation}
where $\la$ denotes a partition and $f$ is any of $e,\,h$ or $p$.
Recall also the Cauchy formula
\begin{equation}
\prod_{i,j}(1-x_iy_j)^{-1}= \sum_\la h_\la(x)\, m_\la(y)= \sum_\la
z_\la^{-1} p_\la(x)\, t\la(y)= \sum_\la s_\la(x)s_\la(y)
\end{equation}
where $s_\la$ stands for the Schur functions and, for $\la=
(1^{m_1}2^{m_2}\cdots)$, $z_\la= \prod_i i^{m_i} m_i!$. From this, a
``combinatorial'' scalar product can be defined:
\begin{equation}\label{SPco}
\LL s_\la\,| \,  s_\mu \RR=  \LL m_\la\,| \,  h_\mu \RR=
\delta_{\la,\mu} \, \qquad \LL p_\la\,| \,  p_\mu \RR= z_\la\,
\delta_{\la,\mu}\, .
\end{equation}
Note that given the bases $e_\la$ and $m_\la$, we can recover the
definition of the partition conjugation $\la'$ by enforcing that
$e_{\la'}$ has the triangular decomposition
$e_{\lambda'}=m_{\lambda}+{\rm smaller~terms}$ when expanded in the
monomial basis.

The idea is to lift all this structure to superspace.

\subsubsection{Physical construction of the  Jack
superpolynomials}

Another line of attack is to start with a
superspace extension of a general class of superpolynomials out of
which all other simple bases can be extracted.  A sufficiently rich
basis for this purpose is the one given by the superanalogs of the
Jack polynomials \cite{Stan}.  Recall that the ordinary Jack
polynomials depend upon a free parameter, denoted $\beta$, and that
various interesting bases are recovered in the appropriate limits:
the symmetric monomial basis when $\beta\rw 0$, the elementary basis
(up to conjugation) when $\beta \rw \infty$  and the Schur
polynomials when $\beta=1$.

The first difficulty with this approach is to have a well-defined
way of generating the proper superextension of the Jack polynomials.
But here again, the pathway is dictated by supersymmetry, or more
precisely, by the consideration of a problem in supersymmetric
quantum mechanics \cite{Witten}.\footnote{In that vein, it is
interesting to point out that it is precisely within the framework
of supersymmetric quantum mechanics that supersymmetry first played
a significant role in mathematics with the influential contribution
of Witten in Morse theory \cite{Morse} (see \cite{Gieres} for more
on the relation between mathematics and supersymmetry).} Recall that
the Jack polynomials $J_\la$ are uniquely characterized by their
triangular decomposition in the monomial basis together with the
following eigenfunction property:
\begin{equation}
H_0J_\la =   \left\{\sum_i (x_i \partial_{x_i})^2+\beta
\sum_{i<j}\frac{x_i+x_j}{x_{i}-x_j}(x_i
\partial_{x_i}-x_j\partial_{x_j})\right\}J_\la = \epsilon_\la
J_\la\, .
\end{equation}

It turns out that $H_0$ is exactly the Hamiltonian for the
trigonometric Calogero-Moser-Sutherland (tCMS) model without the
contribution of its ground-state wave function \cite{For, LV}. The
tCMS model is a completely integrable quantum $N$-body problem
\cite{CMS} that has a unique supersymmetric extension \cite{sCMS}.
Jack superpolynomials are thus naturally defined from the
supersymmetric tCMS eigenvalue problem \cite{DLM1, DLM2,
DLM3}.\footnote{The terms ``superanalogs of Jack polynomials'',
``super-Jack polynomials'' and ``Jack superpolynomials'' have also been
used in the literature for somewhat different polynomials. In
\cite{Serg}, superanalogs of Jack polynomials designated the
eigenfunctions of the CMS Hamiltonian constructed from the root
system of the Lie superalgebra $su(m,N-m)$ (recall that to any root
system corresponds a CMS model \cite{OP}). (The same objects are
called super-Jack polynomials in \cite{VS}.)
   But we
stress that such a Hamiltonian does not contain anticommuting
variables, so that the resulting eigenfunctions are quite different
from our Jack superpolynomials.  Notice also that in
\cite{DLM1,DLM2}, we used the term ``Jack superpolynomials'' for
eigenfunctions of the stCMS model that decompose triangularly in the
supermonomial basis. However, these are not necessarily orthogonal.
The construction of orthogonal Jack superpolynomials was presented in
\cite{DLM3} and from now on, when we refer to ``Jack
superpolynomials'', we refer to the orthogonal ones.}

Quantum mechanics is a special quantum field theory in $0+1$ (no
space and one time) dimension. For a non-supersymmetric $N$-body
problem, we introduce $N$ position operators ${\hat x}_i$
($i=1,\cdots N$) and their canonical conjugate ${\hat p}_i$, the
momenta operators, subject to the commutation relations
\begin{equation}
[{\hat x}_j, {\hat x}_k]=[{\hat p}_j, {\hat p}_k]= 0\, , \qquad
[{\hat x}_j, {\hat p}_k]= \mathrm{i}\delta_{jk}
\end{equation}
(setting $\hbar=1$). In  the Schr\"odinger picture, we work with a
differential realization of these variables: ${\hat x}_j$ is
replaced by the ordinary variable $x_j$ while $p_j$ is replaced by $
-\mathrm{i}\partial_{x_j}$. Now given a Hamiltonian, how do we
supersymmetrize it? The signature of a supersymmetric system is the
presence of a fermionic conserved operator, conserved in the sense
that it commutes with the Hamiltonian. But it is clearly impossible
to construct a fermionic operator without having fermionic
coordinates.  The first step  thus amounts to introduce Grassmannian
partners to our phase space variables (i.e., positions and momenta).
Call them ${\hat\ta}_i$ and ${\hat\ta}_i^\dagger$, with $i=1,\cdots
N$. We enforce the canonical anticommutation relations;
\begin{equation}
     \{{\hat \ta}_j, {\hat \ta}_k\}=  \{{\hat \ta}_j^\dagger, {\hat
\ta}_k^\dagger\}= 0 \, ,\qquad  \{{\hat \ta}_j, {\hat \ta}_k^\dagger
\}=\delta_{jk}\, .
\end{equation}
Again, it is convenient to work with a differential realization:
${\hat \ta_i}\rw {\ta_i}$ and ${\hat \ta}_i^\dagger\rw
\partial_{\ta_i}$.  Now, we still have to specify a procedure for
constructing the supersymmetric Hamiltonian. This turns out to be
rather simple. Introduce two fermionic quantities:
\begin{equation}
Q= \sum_{i=1}^N \ta_i A_i(x)\, , \qquad Q^\dagger= \sum_{i=1}^N
\ta_i^\dagger A^\dagger_i(x)\, ,
\end{equation}
with $A_i$ and $A_i^\dagger$ yet to be determined and set
\begin{equation}
H= \{ Q, Q^\dagger\}= H_0+H_1\, ,
\end{equation}
where $H_0=H(\ta_i=0)$. In other words, we impose that $H_0$, the
part of the Hamiltonian independent of the fermionic variables, be
equal to the non-supersymmetric Hamiltonian that we are trying to
supersymmetrize. This fixes  $A_i$ and $A_i^\dagger$, which in turn
specifies $H_1$, and thus $H$. The construction ensures that both
$Q$ and $Q^\dagger$ are conserved, e.g.,
\begin{equation}
[Q, \{ Q, Q^\dagger\}]= Q(Q  Q^\dagger+  Q^\dagger Q)-(Q Q^\dagger+
Q^\dagger Q) Q=0
\end{equation}
since $Q^2= ( Q^\dagger)^2=0$.

We have thus a clear procedure for supersymmetrizing the Jack
polynomial eigenvalue problem, that is, for obtaining the
supersymmetric tCMS (called stCMS for short) model.  That this
indeed leads to orthogonal superpolynomials that decompose
triangularly in the supermonomial basis has been established in
\cite{DLM3}.

Now the orthogonality just alluded to is with respect to the
so-called physical scalar product:
\begin{equation}
\label{physca}\langle
A(x,\theta)|B(x,\theta)\rangle_{\beta,N}=\prod_{1\leq j\leq
N}\frac{1}{2\pi \mathrm{i}} \oint \frac{ dx_j}{x_j}\int
d\theta_j\,\theta_j\prod_{\substack{1\leq k, l\leq N\\k\neq
l}}\left(1-\frac{x_k}{x_l}\right)^\beta A({\bar x},{\bar \theta})\,
B(x,\theta)\, ,
\end{equation} where the ``bar conjugation'' is defined
as
\begin{equation}\label{defcomplex} {\bar x}_j= 1/x_j\quad\mbox{and}\quad
\overline {(\theta_{i_1}\cdots\theta_{i_m})}\theta_{i_1}\cdots\theta_{i_m}=1
\;.
\end{equation}

\subsubsection{Physical vs combinatorial Jack superpolynomials}

We construct in this article a superspace extension of the
classical bases $m_{\la},e_\la,\, h_\la$ and
$p_\la$ by standard combinatorial methods. But we had also
previously constructed a
one-parameter family of orthogonal polynomials in superspace
reducing to the Jack polynomials
when the fermionic variables are equal to zero.
The question is thus whether
these constructions are arbitrary or somehow belong to the ``proper'' superspace extension
of symmetric function theory.  We will give two reasons why we believe the latter holds.

The first reason has already been mentioned at the beginning of the
previous subsection: the aforementioned combinatorial bases are
recovered as special cases of the Jack polynomials in superspace
just as they are in the non-supersymmetric case.

A second and stronger reason comes from observing that
the $\beta$-deformation of
the combinatorial supersymmetric scalar product can be used to provide
another definition of the Jack polynomials in superspace.
Let us recall that there exists a
one-parameter deformation of the scalar product (\ref{SPco}) between the
$p_\la$'s. The usual Jack polynomials $J_\la$ can be defined purely
combinatorially by enforcing orthogonality with respect to this
deformed scalar product, in addition with a triangularity requirement
(i.e., the $J_\la$'s decompose triangularly in the monomial basis
$\{m_\la\}_\la$ with respect to the
  dominance ordering).
Quite remarkably, the Jack polynomials are also orthogonal
with respect to the physical scalar product induced by the CMS
model (which is simply (\ref{physca}) without the
$\theta$ dependence).  In this sense, one could say that the
physical and combinatorial scalar products are compatible, as
the physical and combinatorial definitions give rise to the very same
objects.

It occurs that it is rather immediate to $\beta$-deform
the superspace extension of (\ref{SPco}). The question is thus whether
the physical  Jack superpolynomials, eigenfunctions of the stCMS quantum
many-body problem, are also orthogonal with respect to this combinatorial
product. We will show in this article that this is indeed the case.
We thus end up with the remarkable conclusion that
our two lines of approach for building a theory of symmetric functions in
superspace, the combinatorial and physical ones,
  yield the very same objects.

\subsection{Organization of the article}


     The article is organized as follows. Section 2 first introduces the concept of superpartition.
     Then relevant results concerning the Grassmann algebra and
   symmetric superpolynomials are reviewed.  A simple interpretation
     of the later, in terms of differentials forms, is also given.
     This section also includes  the definition of
     supermonomials and a formula for their products.

     Section 3 gives the superspace analog of the well known
     elementary symmetric functions, completely symmetric functions  and
     power-sum bases. The generating function for each of them
is displayed.  Determinantal
     formulas that generalize classical formulas describing basic
relations between the
basis elements are presented. Furthermore, orthogonality and duality
relations are
established. In the final subsection, we present a one-parameter deformation of
 the scalar product, the duality transformation and the homogeneous basis.

    Section 4 starts with a review of basic facts concerning our previous
(physical) construction  of  Jack polynomials in superspace.  These
functions are then linked
to the combinatorial theory of symmetric superpolynomials elaborated
in Section 3 in two
different and independent ways. First, it is shown that the
physical Jack superpolynomials are also orthogonal with respect to
the combinatorial product introduced in Section~3.5. And later,
it is shown
that in  non-trivial limiting cases (i.e., special
values of the free parameter or particular superpartitions),
the physical Jack superpolynomials reduce to the symmetric
superfunctions constructed previously in Section~3.

    We finally present, in the conclusion, some natural extensions of
this work. In
particular, we give a precise conjecture concerning
 the existence of (combinatorial)
Macdonald superpolynomials.

As already indicated, this work concerns, to a large extent, a
generalization of symmetric function theory. In laying down its
foundation, we generalize a vast number of basic results from
this theory which can be found for instance in
\cite{Mac} and \cite{StanBook} (Chap. 7).  Clearly, the core of
most of our derivations is bound to be a variation around the proofs
of these older results. We have chosen not to refer everywhere to
the relevant ``zero-fermionic degree'' version of the stated results.
But we acknowledge our debt in that regard  to these two classic
references.  For the  results pertaining specifically to the Jack
polynomials, we have relied heavily on the seminal paper
\cite{Stan} without complete credit in the bulk of the paper, again
to avoid overquoting.


\section{Foundations}

\subsection{Superpartitions}

We recall that a partition
$\lambda=(\lambda_1,\lambda_2,\ldots,\lambda_\ell)$ of $n$, also
written as $\la\vdash n$, is an ordered set of integers such that:
$\lambda_1\geq\lambda_{2}\geq\ldots\geq\la_{\ell}\geq 0$ and
$\sum_{i=1}^{\ell}\la_i=n$. A particular
juxtaposition of two partitions gives a superpartition.

\begin{definition}A
{\it superpartition} $\Lambda$ in the $ m $-fermion sector is a
sequence of non-negative integers separated by a semicolon such that
the sequence before the semicolon is
a partition with
$m$ distinct parts, and such that the remaining sequence is a usual
partition.
That is,  \begin{equation}\label{defsuperpart}
\Lambda:=(\Lambda_1,\ldots,\Lambda_m;\Lambda_{m+1},\ldots,\Lambda_{N})\,
,\end{equation} where $ \Lambda_i>\Lambda_{i+1}\geq 0$ for $ i=1,
\ldots m-1$ and $\Lambda_j \ge \Lambda_{j+1}\geq 0 $ for
$j=m+1,\dots,N-1$.
\end{definition}

Given $\Lambda=(\Lambda^a ;\Lambda^s)$,  the partitions $\Lambda^a$
and $\Lambda^s$ are respectively called the {\it antisymmetric and
the symmetric components} of the superpartition
$\Lambda$.\footnote{From now on,  superscripts $a$ and $s$ refer
respectively to strictly decreasing and decreasing sequences of non-negative
integers.} The {\it
bosonic and fermionic degrees} of $\La$ are
$|\Lambda|=\sum_{i=1}^{N}\Lambda_i$ and
$\overline{\underline{\Lambda}} = m$, respectively. Note that, in
the zero-fermion sector, the semicolon is usually omitted and $
\Lambda $  reduces then to $\La^s$.

We say that the ordered set $\Lambda$ in (\ref{defsuperpart}) is
a superpartition of $(n|m)$ if $|\Lambda|=n$ and
$\overline{\underline{\Lambda}} = m$; in symbols, this is written as
$\La\vdash(n|m)$. The set composed of all superpartitions of $(n|m)$
is denoted $\mbox{SPar}(n|m)$. When the fermionic degree is zero, we
recover standard partitions:
$\mbox{SPar}(n|0)=\mbox{Par}(n)$.



We also define
\begin{equation} \mbox{SPar}(n):=\bigcup_{m\geq
0}\mbox{SPar}(n|m)\quad\mbox{and}\quad \mbox{SPar}:=\bigcup_{m,n\geq
0}\mbox{SPar}(n|m)\, ,\end{equation} with $
\mbox{SPar}(0|0)=\emptyset$ and $\mbox{SPar}(0|1)=\{(0;0)\}$. For
example, we have
\begin{equation}\mbox{SPar}(3|2)=\{\,(3,0;0),\,
(2,1;0),\,(2,0;1)\,(1,0;2),\,(1,0;1,1)\, \}\, .\end{equation} Notice
that $\mathrm{SPar}(n|m)$ is empty for all $n<m(m-1)/2$.

    We
introduce the operator ${\ell}:\,\mbox{SPar}\rightarrow\mathbb{N}$
that gives the length of a superpartition as
\begin{equation}{\ell}(\La):=\overline{\underline{\Lambda}}+{\ell}(\Lambda^s)\quad\mbox{for}\quad
{\ell}(\Lambda^s):=\mbox{Card}\{\La_i\in\La^s \,:\,\La_i>0\}\,
.\end{equation}
     With this
definition, ${\ell}\bigl((1,0;1,1)\bigr)=2+2=4$ (i.e., a zero-entry in $\La^a$ contributes to the length of $\La$). To every superpartition
$\La$, we can also associate a unique partition $\La^*$ obtained by
deleting the semicolon and  reordering the parts in non-increasing
order.  For instance,
\begin{equation}
(5,2,1,0;6,5,5,2,2,1)^*=(6,5,5,5,2,2,2,1,1,0)=(6,5,5,5,2,2,2,1,1)\,.
\end{equation}
From this, we can introduce another notation for
superpartitions. A superpartition $\La=(\La^a;\La^s)$ can be viewed
as the partition $\Lambda^{*}$ in which every part of $\La^a$ is
circled.
 If a part ${\La^a}_j=b$ is equal to at least one part of
$\La^s$, then we circle the leftmost $b$ appearing in $\La^*$. We
shall use $C[\Lambda]$ to denote this special notation.
   For instance, \begin{equation}\Lambda=(3,1,0;4,3,2,1) \iff
C[\La]=(4,\cercle{3},3,2,\cercle{1},1,\cercle{0}).\end{equation}

This allows us to introduce a diagrammatic representation of
superpartitions. To each $\La$, we associate a unique diagram,
denoted by $D[\La]$.  It is obtained by first drawing the Ferrer's diagram
associated to $C[\La]$, that is, by drawing  a diagram with ${C[\La]}_1$ boxes
in the first row, ${C[\Lambda]}_2$ boxes in the second row and so
forth, all rows being left justified.  If, in addition, the integer
${C[\La]}_j=b$ is circled, then we add a circle at the end of the
$b$  boxes in the $j$-th row. For example,
    \begin{equation}
    D[3,1,0;4,3,2,1]=
{\tableau[scY]{&&&\\&&&\bl\tcercle{}\\&&\\&\\&\bl\tcercle{}\\&\bl\\\tcercle{}\bl\\
}}\end{equation} The {\it conjugate} of a superpartition $\La$,
denoted by $\La'$, is obtained by interchanging the rows and the
columns in the diagram;  in matrix notation, we have
$D[\La']=(D[\La])^{\mathrm{t}}$ if $\mathrm{t}$ stands for the
transpose operation. Hence, $(3,1,0;4,3,2,1)'=(6,4,1;3)$ since
\begin{equation}\left(\,{\tableau[scY]{&&&\\&&&\bl\tcercle{}\\&&\\&\\&\bl\tcercle{}\\&\bl\\\tcercle{}\bl\\
}}\,\right)^\mathrm{t}={\tableau[scY]{&&&&&&\bl\tcercle{}\\&&&&\bl\tcercle{}\\&&\\&\bl\tcercle{}\\
}}\end{equation}Recall that, for any partition $\la$, the
conjugation can be defined by ${\la'}_j=\mbox{Card}\{k\,:\,\la_k\geq
j\}$. Thus, in symbols, the conjugation of $\La\in\mbox{SPar}(n|m)$
reads \begin{equation}\La'=({\La'}^a ; {\La'}^s)\, ,\end{equation}
where\begin{equation}{\La'}^s=\left({\La^*}'\setminus{\La'}^a\right)^+\aand
{{\La'}^a}_j=\mbox{Card}\{k\, :
\,\La_k>\La_{m+1-j}\}\,,\end{equation} with $\la^+$ standing  for the
partition obtained by reordering the parts of $\la$ non-increasingly.
Obviously, the conjugation of any superpartition $\La$ satisfies
\begin{equation}(\La')'=\La\quad\mbox{and}\quad(\La^*)'=(\La')^*\,
.\end{equation}

\vskip0.2cm
\begin{remark}
The description of a superpartition in terms of partition with some parts circled makes clear that overpartitions \cite{CoLo}  are special cases of superpartitions. Indeed, overpartitions are circled superpartitions (with the circle replaced by an overbar) that do not contain a possible circled zero. If we denote by $s_N(n|m)$ the number of superpartitions
$\La\in\spar(n|m)$ such that ${\ell}(\La)\leq N$, then this connection makes clear that their generating function is
\begin{equation}\label{gefct}
 \sum_{n,m, p\geq 0} s_{m+p} (n|m)
\,z^m y^{p} q^n=\frac{(-z;q)_\infty}{ (yq;q)_\infty} \quad {\rm with }\quad
(a;q)_\infty:=\prod_{n\geq 0}(1-aq^n)
\end{equation}
\end{remark}

To complete this subsection, we consider the natural ordering on
superpartitions, which is defined in terms of the
Bruhat order on
compositions. Recall that a composition of $n$ is simply a sequence
of non-negative integers whose sum is equal to $n$; in symbols
$\mu=(\mu_1,\mu_2,\ldots)\in \mathrm{Comp}(n)$ iff $\sum_i\mu_i =n$
and $\mu_i\geq0$ for all $i$.  The Bruhat ordering on compositions is
defined as follows. Given a composition $\la$, we let $\la^+$ denote
the partition obtained by reordering its parts in non-increasing
order. Now, $\la$ can be obtained from $\la^+$ by a sequence of
permutations. Among all permutations $w$ such that $\lambda= w
\lambda^+$, there exists a unique one, denoted $w_{\lambda}$, of
minimal length. For two compositions $\la$ and $\mu$, we say that
$\lambda \geq \mu $ if either $\lambda^+ > \mu^+$ in the usual dominance ordering or
$\lambda^+=\mu^+$ and $w_{\lambda} \leq w_{\mu}$ in the sense that
the word $w_{\lambda}$ is a subword of $w_{\mu}$ (this is the Bruhat
ordering on permutations of the symmetric group). Recall that for two partitions $\lambda$ and $\mu$ of the
same degree, the dominance ordering is: $\lambda \geq \mu$ iff
${\lambda}_1+\ldots+{\lambda}_k\geq{\mu}_1+\ldots+{\mu}_k$ for all
$k$.

Let $\La$ be a superpartition of $(n|m)$.  Then, to $\La$ is
associated  a unique composition of $n$, denoted by $\La^c$,
obtained by replacing the semicolon in $\Lambda$ by a comma.  We
thus have $\mathrm{Spar}(n)\subset\mathrm{Comp}(n)$, which leads to
a natural Bruhat ordering on superpartitions.

\begin{definition}\label{bruhatorder}Let $\La,\,\Om\,\in\,\spar
(n|m)$. The {\it Bruhat order}, denoted by
$\leq$, is such that $\Omega\leq\Lambda$ if $\Om^c\leq\La^c$.
\end{definition}

We need two refinements of the previous order, namely the $S$ and
$T$ orders (the origin of these names will become clearer in the
following lines).
\begin{definition}\label{storder}Let $\La,\,\Om\,\in\,\spar
(n|m)$. The {\it $S$ and $T$ orders} are respectively defined as
follows:\footnote{The $S$ order is the precisely the ordering
introduced in \cite{DLM1} but it differs from the more precise
ordering of \cite{DLM2}, called there $\leq^s$. In \cite{DLM3}, it
is called the $h$ ordering.  See also appendix B of \cite{DLM4}.}
\begin{equation}\begin{array}{lllllll}
\Omega\leq_S\Lambda&\mbox{if either}&\Om=\La&\mbox{or}&
\Om^*<\La^*\, ,&\\
\Omega\leq_T\Lambda&\mbox{if either}&\Om=\La&\mbox{or}& \Om^*=\La^*
&\mbox{and}&\Om^c<\La^c\, .\end{array}\end{equation}
\end{definition}
\n In order to describe other characterizations of these orders, we
need the following operators on compositions (or superpartitions):
\begin{equation}\label{defopST}\begin{array}{ll}
S_{ij}(\ldots,\la_i,\ldots,\la_j,\ldots)=&\ds\left\{\begin{array}{ll}
(\ldots,\la_i-1,\ldots,\la_j+1\ldots)&\mbox{if}\quad\la_i-\la_j>1\,
,\\
(\ldots,\la_i,\ldots,\la_j,\dots) &\mbox{otherwise}\,
,\end{array}\right.\\&\\
T_{ij}(\ldots,\la_i,\ldots,\la_j,\ldots)=&\ds\left\{\begin{array}{ll}
(\ldots,\la_j,\ldots,\la_i,\ldots)&\mbox{if}\quad\la_i-\la_j>0\,
,\\
(\ldots,\la_i,\ldots,\la_j,\dots) &\mbox{otherwise}\,.
\end{array}\right.
\end{array}\end{equation}

\begin{lemma}\label{lemmaST}\cite{Mac, manivel}  Let $\la$ and $\om$
be two compositions of $n$.
Then, $\la^+>\om^+$ iff there exists a sequence $\{S_{i_1,j_1},
\ldots, S_{i_k,j_k}\}$ such that
\begin{equation} \om^+=S_{i_1,j_1}\cdots S_{i_k,j_k}\la^+\, .\end{equation}
Similarly, $\la^+=\om^+$ and $\la>\om$ iff there exists a sequence
$\{T_{i_1,j_1}, \ldots, T_{i_k,j_k}\}$ such that
\begin{equation} \om=T_{i_1,j_1}\cdots T_{i_k,j_k}\la\, .\end{equation}
This last property can be translated for superpartitions as:
$\La>_T\Om$ iff $D[\Om]$ can be obtained by moving step by step in
the down-left direction the circles of $D[\Lambda]$.
\end{lemma}

At this stage, we are in a position to establish the fundamental
property relating conjugation and Bruhat order which is that
  the Bruhat order is anti-conjugate (in the sense of the following proposition).

\begin{proposition}\label{uty}
Let $\La,\,\Om\, \in\, \spar (n|m)$. Then \begin{equation}\La\geq\Om
\quad\Longleftrightarrow\quad \Om'\geq\La'\,
.\end{equation}\end{proposition}
\begin{proof}It suffices to prove the result for the $S$ and $T$
orderings. The case $\La>_S\Om$, that is, $\La>\Om$
and $\La^*\neq\Om^*$, is a well-known result
on partitions (see for instance (1.11) of \cite{Mac}).

We now consider $\La>_T\Om$, that is, $\La>\Om$ with $\La^*=\Om^*$.
Lemma \ref{lemmaST} tells us that in this case $D[\Om]$ can be
obtained by moving step by step in the down-left direction the
circles of $D[\Lambda]$. Under transposition, this amounts to saying
that $D[\Om']$ can be obtained by moving step by step in the
up-right direction the circles of $D[\Lambda']$, or equivalently
that  $D[\La']$ can be obtained by moving step by step in the
down-left direction the circles of $D[\Omega']$. Hence,
\begin{equation}\La>_T\Om\quad\Longrightarrow\quad \La'<_T\Om'\,
.\end{equation}
We have thus proved $\La>\Om\Rightarrow \La'<\Om'$. Since
conjugation is an involution, the claim holds.\end{proof}

   Consider for instance
$\La= (3,0;4,1)$ and $\Omega= (2,0;4,2)= S_{14}\La <_S \La$. We find that
$\La'= (3,1;2,2)$ and $\Om'=(3,1;3,1)$ so that $\La'=S_{34}\Om'$,
i.e. $\Om'>_S\La'$. Consider also $\Gamma= (1,0;4,3)= T_{14}\La <_T \La$,
which gives $\Gamma'=(3,2;2,1)>_T\La'=T_{24}\Gamma'$.

\vskip0.2cm
\begin{remark}Notice that we could have introduced
 as an alternative ordering
the dominance ordering on superpartitions,  denoted by $\leq_D $ and defined as follows: $\Omega\leq_D\Lambda$ if either $\Omega^*<\Lambda^*$ or $\Omega^*=\Lambda^*$ and
${\Omega}_1+\ldots+{\Omega}_k\leq {\Lambda}_1+\ldots{\Lambda}_k,$
$\forall\,\, k\, .$
 The
  usefulness of this ordering in special contexts lies in its simple
description in terms of inequalities.   However, it is not the
proper generalization of the dominance order on partitions
because it is not anti-conjugate
as will be illustrated later
by an example.
In fact, it is not as strict
as  the  Bruhat ordering (i.e., more superpartitions are
comparable in this order than in the Bruhat ordering).
This follows from the second property of Lemma \ref{lemmaST} which
obviously implies that for superpartitions, the
Bruhat ordering  is a weak subposet of the dominance
ordering, that is,
$ \La
\geq\Om\,\Rightarrow\,\La \geq_D\Om$.
However the converse is not true.   For instance, if
$\La=(5,2,1;4,3,3)$ and $\Om=(4,3,0;5,3,2,1)$ we easily verify that
$\La>_D \Om$.
But $\La'=(5,4,0;6,2,1)$ and $\Om'=(6,2,1;5,4)$, so that
$\Om'\not >_D\La'$.  Proposition \ref{uty}
 implies that the two orderings need to be distinct. This corrects a loose
implicit statement in \cite{DLM3} concerning the expected
equivalence of these two orderings.
\end{remark}

\subsection{Ring of symmetric superpolynomials}

Let $\mathscr{B}=\{B_j\}$ and $\mathscr{F}=\{F_j\}$ be the formal
and  infinite sets composed of all bosonic (commutative) and
fermionic (anticommutative) quantities respectively:
\begin{equation}[B_j,B_k]=[B_j,F_k]=\{F_j,F_k\}=0\, .\end{equation} Thus,
$\mathscr{S}=\mathscr{B}\oplus\mathscr{F}$ is $\mathbb{Z}_2$-graded
over any ring $\mathbb{A}$ when we identify $\,^0\!\mathscr{S}$ with
$\mathscr{B}$ and $\,^1\!\mathscr{S}$ with $\mathscr{F}$. The degree
of any element $s$ of $\mathscr{S}$, written $\hat{\pi}(s)$, is
defined via \begin{equation}\hat{\pi}(s)=\left\{\begin{array}{ll}
0,& s\in \mathscr{B},\cr 1, & s\in \mathscr{F}.
\end{array}\right.\end{equation}
Consequently, $\mathscr{S}$ possesses a parity operator (involution)
$\hat{\Pi}$ defined by
\begin{equation}\hat{\Pi}(s)=(-1)^{\hat{\pi}(s)} s\end{equation} and
satisfying \begin{equation}\hat{\Pi}^2=1,\quad
\hat{\Pi}(st)=\hat{\Pi}(s)\hat{\Pi}(t),\quad
\hat{\Pi}(as+bt)=a\hat{\Pi}(s)+b\hat{\Pi}(t)\, ,\end{equation} for
all $a,b\in \mathbb{A}$ and $s,t\in\mathscr{S}$.
The second relation implies that the product of two elements of $\mathscr{F}$ belongs to
$\mathscr{B}$, i.e.,  ``the product of two fermions is a boson''.

As an example of such a structure, we consider the Grassmann algebra
over a ring $\mathbb{A}$, denoted $\mathscr{G}_M(\mathbb{A})$.  It
is a non-commutative algebra with identity $1\in\mathbb{A}$,
generated by the $M$ elements $\theta_1,\ldots,\theta_M $ that
belong to the ring of Grassmannian variables $\mathbb{G}_M$, with
product defined in (\ref{grassmannvar}). Every element $g$ in
$\mathscr{G}_N(\mathbb{A})$ can be decomposed as $ g= \,^0\!g+ \,^1\!g$, with \footnote{It should
be stressed that the following sum is formal since physically it
makes no sense to add fermionic and bosonic quantities. More
precisely, an equality between two elements of $g$ is an equality
between the $\,^0\!g$ and $\,^1\!g$ parts of these elements, much like the
real and imaginary parts of complex equations.}
\begin{equation}\label{decompograd}
\,^\epsilon\!g= \ds
 a(1-\epsilon)+\sum_{\substack{k=\epsilon ~{\rm mod}\; 2,\, k>0,\\ 1\leq i_1<\ldots<i_k\leq
 M}}a^{i_1,\ldots,i_k}\theta_1\cdots\theta_k\;,
\end{equation}
where $\epsilon=0,1$ and
with the constants $a$ and $a^{i_1,\ldots,i_k}$ belonging to
$\mathbb{A}$. The dimension of this Grassmann algebra is thus
$\sum_{j=0}^N {\binom{N}{j}}=2^N$. When $\mathbb{A}$ is a field, the
subalgebra $\,^0\mathscr{G}_M(\mathbb{A})$ is also a field in which
the inverse is defined by (writing $\,^0\!g= a+f$)
\begin{equation}\label{invGrass}( \,^0\!g)^{-1} =
\frac{1}{a+f}=\frac{1}{a}\sum_{n\geq
0}(-1)^{n}\left(\frac{f}{a}\right)^{n}\, . \end{equation} Note that,
due to the nilpotency of the $\theta_j$'s, the term $(a+f)^{-1}$ is
finite for all $M<\infty$. For instance
\begin{equation}
({1-\theta_1\theta_2-\theta_3\theta_4})^{-1}=1+\theta_1\theta_2+\theta_3\theta_4+2\theta_1\theta_2\theta_3\theta_4\,
.\end{equation}

Before going further, we define another involution on the Grassmann
algebra:\begin{equation}\overleftarrow{\,\theta_{j_1}\cdots\theta_{j_m}\,}:=\overrightarrow{\,\theta_{j_m}\cdots\theta_{j_1}\,}\,
, \end{equation} where
\begin{equation}\overrightarrow{\,\theta_{j_1}\cdots\theta_{j_m}\,}:=\theta_{j_1}\cdots\theta_{j_m}\,
.\end{equation}
    In words, the operator $\overleftarrow{\phantom{abc}}$ reverses the order
of the anticommutative variables while
$\overrightarrow{\phantom{abc}}$  is simply the
identity map. \footnote{The explicit use of
$\overrightarrow{\phantom{abc}}$ is not essential.  Nevertheless, it
will make many formulas more symmetric and transparent.}  Using
 induction, we  get
\begin{equation}\overleftarrow{\,\theta_{j_1}\cdots\theta_{j_m}\,}=(-1)^{m(m-1)/2}\,
\overrightarrow{\,\theta_{j_1}\cdots\theta_{j_m}\,}\,.\end{equation}
    This result immediately  implies the following simple properties.

\begin{lemma}\label{ordredestheta} Let $\{\theta_1,\ldots,\theta_M\}$
and $\{\phi_1,\ldots,\phi_M\}$ be two sets of Grassmannian
variables.  Then \begin{equation} (\theta_{j_1}\phi_{j_1})\cdots
(\theta_{j_m}\phi_{j_m})=\overleftarrow{(\theta_{j_1}\cdots\theta_{j_m})}\overrightarrow{(\phi_{j_1}\cdots\phi_{j_m})}
=\overrightarrow{\,(\theta_{j_1}\cdots\theta_{j_m})\,}\overleftarrow{\,(\phi_{j_1}\cdots\phi_{j_m})\,}\,
\end{equation} and
    \begin{equation}
\overleftarrow{\overleftarrow{(\theta_{j_1}\cdots\theta_{j_m})}\overrightarrow{(\phi_{j_1}\cdots\phi_{j_m})}}
=\overrightarrow{(\phi_{j_1}\cdots\phi_{j_m})}\overleftarrow{(\theta_{j_1}\cdots\theta_{j_m})}\,
.\end{equation}\end{lemma}

Now, let $x=\{x_1,\ldots,x_N\}\subset\mathscr{B}$ and
$\theta=\{\theta_1,\ldots,\theta_M\}\subset\mathscr{F}$. The
superpolynomial algebra over a unital ring $\mathbb{A}$, denoted by
$\mathscr{P}^{N|M}(\mathbb{A})$ or by
$\mathbb{A}[x_1,\dots,x_N,\theta_1,\ldots,\theta_M]$, is the
Grassmann algebra $\mathscr{G}_M$ over the ring $\mathscr{P}^{N}$ of
polynomials in $x$. Recall that $\mathscr{P}^{N}$ is also a graded,
unital and commutative algebra over any unital ring $\mathbb{A}$.

Every $f(x)\in\mathscr{P}^{N}=\bigoplus_{n\geq 0}\mathscr{P}_{(n)}$
can be written as \begin{equation}\label{polx} f(x)=\sum_{n\geq 0}
f^{(n)}(x)\, ,\quad f^{(n)}(x)=\sum_{\alpha\in\mathbb{N}^N,\,
    |\alpha|=n}a_\alpha \, x^\alpha\,\in\mathscr{P}_{(n)}\,
,\end{equation} where
$
    x^\alpha=x_1^{\alpha_1}\cdots x_N^{\alpha_N}$ and $
    a_\alpha\in\mathbb{A}$.
Note that $|\alpha|=\sum_i\alpha_i$ is the degree of the weak
composition $\alpha$; it also corresponds to the degree of
$f^{(n)}(x)$.  Every superpolynomial $f(x,\theta)$ in
$\mathscr{P}^{N|M}(\mathbb{A})$ possesses  a bosonic and a fermionic
part, i.e., \begin{equation}\label{polxtheta}f(x,\theta)=
\,^0\!f(x,\theta)+ \,^1\!f(x,\theta)\, ,\end{equation} where the
components of $f$ are similar to those given in (\ref{decompograd}),
apart from the fact that the $a^{i_1,\ldots,i_k}$'s now belong to
$\mathscr{P}^{N}$. From decompositions (\ref{polx}) and
(\ref{polxtheta}), it is obvious that $\mathscr{P}^{N|M}$ is
bi-graded with respect to the bosonic and fermionic degrees, that
is,
\begin{equation}\mathscr{P}^{N|M}=\bigoplus_{n,m\geq0}\mathscr{P}^{N|M}_{(n|m)}
\, .\end{equation} Each submodule $\mathscr{P}^{N|M}_{(n|m)}$ is
finite dimensional.  It is composed of all homogeneous
superpolynomials $f(x,\theta)$ with degrees $n$ and $m$ in $x$ and
$\theta$, respectively.  We shall write $\mbox{deg} f= (n|m)$ for
any such a superpolynomial $f$.

Pure fermionic polynomials (i.e., elements of
$\mathscr{P}^{N|M}_{(n|m)}$ with $m$ odd) have nice properties. As
an example, consider the following proposition that shall be useful
in the subsequent sections.

\begin{proposition}\label{expfermions}Let
$\tilde{f}=\{\tilde{f}_0,\tilde{f}_1,\ldots\}$ and
$\tilde{g}=\{\tilde{g}_0,\tilde{g}_1,\ldots\}$ be two sequences of
fermionic polynomials parametrized by non-negative integers.  Let
also
\begin{equation}\tilde{f}_{\mu}:=\tilde{f}_{\mu_1}\tilde{f}_{\mu_2}\cdots\aand
\tilde{g}_{\mu}:=\tilde{g}_{\mu_1}\tilde{g}_{\mu_2}\cdots
\end{equation}where $ \mu$ belongs to $
\mathrm{Par}_a(n)$, the set of partitions of $n$ with strictly
decreasing parts.    Then
\begin{equation}\exp\left[\sum_{n=0
}^{M-1}\tilde{f}_n\,\tilde{g}_n\right]=\sum_{n=0}^{M(M-1)/2}\sum_{\mu\in\mathrm{Par}_a(n)}
\overleftarrow{\phantom{\Big|}\tilde{f}_\mu\phantom{\Big|}}\,\overrightarrow{\phantom{\Big|}\tilde{g}_\mu\phantom{\Big|}}\,
.\end{equation}
\end{proposition}
\begin{proof} Due to the fermionic character of $\tilde{f}$ and
$\tilde{g}$, we have
\begin{equation}\begin{array}{lcl}\exp\left[\sum_{n=0
}^{M-1}\tilde{f}_n\,\tilde{g}_n\right]&=&\ds \prod_{0\leq  n\leq
M-1} (1+\tilde{f}_n\,\tilde{g}_n)\\&=&\ds 1+\sum_{0\leq n\leq
M-1}\tilde{f}_n\,\tilde{g}_n+\sum_{0\leq m<n\leq
M-1}\tilde{f}_m\,\tilde{g}_m
\tilde{f}_n\,\tilde{g}_n+\ldots\end{array}\end{equation} Since every
term in the last equality can be reordered by Lemma
\ref{ordredestheta}, the proof follows.\end{proof}

    We finally consider the
symmetric superpolynomials.  We specialize to the case in which the
number of bosonic and fermionic variables is the same, i.e., $N=M$
and $\mathscr{P}:=\mathscr{P}^{N|N}$. The algebra of symmetric
superpolynomials over the ring $\mathbb{A}$, denoted by
$\mathscr{P}^{S_N}(\mathbb{A})$ or by
    $\mathbb{A}[x_1,\dots,x_N,\theta_1,\ldots,\theta_N]^{S_N}$, is a
subalgebra of $\mathscr{P}$.  As mentioned in the introduction,
$\mathscr{P}^{S_N}$ is made out of all $f(x,\theta)\in\mathscr{P}$
invariant under the diagonal action of the symmetric group $S_N$.

To be more explicit, we introduce $K_{ij}$ and $\kappa_{ij}$, two
distinct and Abelian superpolynomial realizations of the
transposition $(i,j)\in
S_N$:\begin{equation}K_{ij}f(x_i,x_j,\theta_i,\theta_j)=f(x_j,x_i,\theta_i,\theta_j)
\, ,
\kappa_{ij}f(x_i,x_j,\theta_i,\theta_j)=f(x_i,x_j,\theta_j,\theta_i)\,
,\end{equation} for all $f\in \mathscr{P}$. Note that $\kappa_{ij}$
has the following realization
\begin{equation}\kappa_{ij}:=1-(\theta_{i}-\theta_{j})
(\d_{\theta_{i}}-\d_{\theta_{j}}).\end{equation}
Since every permutation is generated by products of elementary
transpositions $(i,i+1)\in S_N$, we can define symmetric
superpolynomials as follows.
\begin{definition} \label{defsupersym}
A superpolynomial $f(x,\theta)\in\mathscr{P}$
is {\it symmetric} if and only if \begin{equation}\mc{K}_{i,i+1}
f(x,\theta)=f(x,\theta)\quad\mbox{where}\quad\mc{K}_{i,i+1}:=\kappa_{i,i+1}K_{i,i+1}\end{equation}
for all $i\in\{1,2,\ldots,N-1\}$.\end{definition}

    \n But, every
monomial $\theta_J=\theta_{j_1}\cdots\theta_{j_m}$ is completely
antisymmetric, that is, \begin{equation}
\kappa_{ik}\theta_J=\left\{\begin{array}{ll}
-\,\theta_J,&\mbox{if}\quad i\,, k\,\in J,\\
\phantom{-}\,\theta_J,&\mbox{if}\quad i\,, k\,\not\in
J.\end{array}\right.\end{equation}
    This observation immediately
implies the following result.

\begin{lemma}\label{antisymsym}Let $f(x,\theta)\in \mathscr{P}$ be
expressed as:
\begin{equation} f(x,\theta)=\sum_{m\geq 0}\,\sum_{{1\leq j_1<\ldots<j_m\leq
N}}\,f^{j_1,\ldots,j_m}(x)\,\theta_{j_1}\cdots\theta_{j_m}\, .
\end{equation} If  $f(x,\theta)$ is symmetric, then each polynomial
$f^{j_1,\ldots,j_m}(x)$ is completely antisymmetric in the set of
variables $y:=\{x_{j_1},\ldots,x_{j_m}\}$ and completely symmetric
in the set of variables $x\setminus y$.\end{lemma}

\subsection{Geometric interpretation of superpolynomials}

Symmetric functions can be interpreted as symmetric 0-forms $f$
acting on a manifold: $K_{ij}f(x)=f(x)$ where $x$ is a local
coordinate system.  Similarly, symmetric superfunctions in the
$p$-fermion sector can be reinterpreted as symmetric $p$-forms $f^p$
acting on the same manifold: ${\mathcal K}_{ij}f^p(x)=f^p(x)$.  Thus,
the set of all symmetric superfunctions is in correspondance with
the completely symmetric de Rham complex.  This geometric point
of view is   briefly explained in this subsection. (Note that none of our results relies on this observation.)

We consider a Riemannian manifold $\M$ of dimension $N$ with metric
$g_{ij}$ and  let $x=\{x^1,\ldots,x^N\}$ denote a coordinate system
on $U\subset\M$.
  Let $\T\M$ be the
tangent bundle on which we choose an orthonormal coordinate frame
$\ee=\{\partial_{1},\ldots,\d_{N}\}$. As usual,
$\ee^*=\{dx^1,\ldots,dx^N\}$ denotes the dual basis to $\ee$
belonging to the cotangent bundle $\T^*\M$:
$dx^i(\d_{x^{j}})=\delta^{i}_j$. The set of all $p$-form fields on
$M$ is a vector space denoted by ${\bigwedge}^p$.  Each $p$-form can
be written as
\begin{equation}\label{developform}\alpha^p(x)=\sum_{1\leq j_1<\ldots<j_p\leq N}
\alpha_{j_1,\ldots,j_p}(x) dx^{j_1}\wedge\cdots\wedge dx^{j_p}\,  ,
\end{equation} where the exterior  (wedge)  product is antisymmetric
: $dx^i\wedge dx^j=-dx^j\wedge dx^i$. Let $d$ be the exterior
differentiation on forms, whose action is
\begin{equation}d\, \alpha^p(x)=\sum_{1\leq k,j_1,\ldots,j_p\leq
N}[\d_{x^k}\alpha_{j_1,\ldots,j_p}(x)]dx^k\wedge
dx^{j_1}\wedge\cdots\wedge dx^{j_p}\,  . \end{equation}
    This operation is used to
define the de Rham complex of $\M$:
\begin{equation}0\longrightarrow\mathbb{R}
\longrightarrow{\bigwedge}^0\stackrel{d}{\longrightarrow}
{\bigwedge}^1\stackrel{d}{\longrightarrow}\cdots\stackrel{d}{\longrightarrow}{\bigwedge}^N
\stackrel{d}{\longrightarrow}0\, .\end{equation}

In order to represent our Grassmannian variables $\theta^j$ and
$\theta^j \,^\dagger$  in terms of forms, we introduce the two operators
$\hat{e}_{dx^j}$ and $\hat{\i}_{\d_{x_k}}$,
   where $\hat{e}_\alpha$ and $\hat{\i}_{v}$ respectively stand for the
  left exterior product by the form $\alpha$ and
  the interior product (contraction) with respect to the vector field $v$.
  These operators satisfy a
  Clifford (fermionic) algebra
\begin{equation}\{\hat{e}_{dx^i}\,,\,\hat{\i}_{\d_{x_k}}\,\}=\delta^i_k\aand
\{\hat{e}_{dx^i}\,,\,\hat{e}_{dx^j}\}=0=\{\hat{\i}_{\d_{x_j}}\,,\,
\hat{\i}_{\d_{x_k}}\,\}\, .\end{equation} This implies that the
$\theta^j$'s and  $\theta^j\,^\dagger$, as operators,
    can be realized as follows:
  \begin{equation} \theta^j\sim\hat{e}_{dx^j}\aand \theta^j
  \,^\dagger \sim g^{jk}\,\hat{\i}_{\d_{x_k}}
   \,,\end{equation}
that is,
\begin{equation} dx^j\sim \theta^j(1)\aand\theta^j\,^\dagger\sim
g^{jk}\,\partial_{\theta^k}\wwhere
\partial_{\theta^k}:=\frac{\partial\,}{\partial \theta^k}\, .\end{equation}
Note that introducing the Grassmannian variables as operators is
needed to enforce the wedge product of the forms $dx_j$. Moreover,
if $\alpha^p$ is a generic $p$-form field and $\hat{\pi}\,
:\,{\bigwedge}^p\rightarrow {\bigwedge}^p$ is the operator defined
by
\begin{equation} \hat{\pi}_p:=
\theta^j\partial_{\theta^j}=\theta^jg_{jk}\theta^k\,^\dagger\quad\Longrightarrow\quad
\hat{\pi}_p\alpha^p=p\,\alpha^p\,,  \end{equation}then
$\hat{\Pi}_p:=(-1)^{\hat{\pi}_p}$ is involutive. (Manifestly,
$\hat{\Pi}_p=\hat{\Pi}$, the parity operator introduced previously.)
This involution is also an isometry in the Hilbert space scalar
product.  The operator $\hat{\Pi}_p$ induces a natural $\mathbb{Z}_2$
grading in the de Rham complex.

The construction of  the symmetric de Rham complex, denoted
$\mathrm{SRham}$, is obtained as follows. We make a change of
coordinates: $x\rightarrow f(x)$, where
$f=\{f^n\}:=\{f^1,\ldots,f^N\}$ is an $N$-tuple of symmetric and
independent functions of $x$.  For instance, $f^n$ could be an
elementary symmetric function $e_n$, a complete symmetric function
$h_n$, or a power sum $p_n$ (see Section 3). This implies a change
of basis in the cotangent bundle: $dx\rightarrow df(x)$. Explicitly,
\begin{equation} df^n=\sum_i (\partial_i f^n)(x) \,dx_i\sim
\tilde{f}^n=\sum_i (\partial_i f^n)(x) \,\theta_i\, .\end{equation}
In other words, $d{f}$ is a new set of  ``fermionic'' variables
invariant under any permutation of the $x_j$'s.\footnote{Of course,
this change of basis is well defined in $U\subset\M$ if and only if
the Jacobian determinant $J(f,x)$ is not zero.  For any standard
basis of the symmetric function space (i.e., the $e_n$'s, $h_n$'s or
$p_n$'s), we easily verify that the Jacobian is, up to a sign, equal
to the Vandermonde determinant.   Thus, the
change of coordinates is well defined for any locus $U$  such as the
following one: $L_U:=\{x(p)\,:\,x^1(p)>x^2(p)>\ldots\, ;\,\forall\,
p\in U\}\, .$}

  These remarks explicitly show that
   symmetric polynomials in superspace can be interpreted as
  symmetric differential forms.  We stress that the diagonal
   action of the symmetric group $S_N$ comes naturally in the
   geometric perspective.  Note finally
that for a Euclidian superspace (relevant to our context), the
position (upper or lower) of the indices does not matter.

\subsection{Supermonomial basis}

The symmetric supermonomial function, denoted by $m_{\Lambda}=
m_{\Lambda}(x,\theta)$, is  the superanalog of the monomial
symmetric function.  It is defined as follows \cite{DLM1}.
\begin{definition} \label{defsupermono}To each superpartition
$\La\in\mathrm{SPar}(n|m)$,
we associate a superpolynomial
$m_{\Lambda}\in\mathscr{P}^{S_N}_{(n|m)}$, called {\it symmetric
supermonomial}, defined by
\begin{equation}
m_{\Lambda}={\sum_{\sigma\in S_{N}}}' \theta^{\sigma(1, \ldots,
m)}x^{\sigma(\Lambda)} \, ,
\end{equation} where the prime indicates that the  summation is
restricted to distinct terms, and where
\begin{equation}
x^{\sigma(\Lambda)}=x_1^{\Lambda_{\sigma(1)}} \cdots
x_m^{\Lambda_{\sigma(m)}} x_{m+1}^{\Lambda_{\sigma(m+1)}} \cdots
x_{N}^{\Lambda_{\sigma(N)}}\quad\mbox{and}\quad\theta^{\sigma(1,
\ldots, m)} = \theta_{\sigma(1)} \cdots \theta_{\sigma(m)} \, .
\end{equation}\end{definition}

\noindent Obviously, the previous definition can be replaced by the
following:
\begin{equation} \label{monod}
         m_{\Lambda}= \frac{1}{n_{\Lambda}!} \sum_{\sigma \in S_N}
\mc{K}_{\sigma}
         \left( \theta_1 \cdots \theta_m x^{\Lambda}\right) \end{equation}
    for
\begin{equation} \label{augmented}
 n_{\Lambda}! =n_{\Lambda^s}!:=
n_{\Lambda^s}(0)!\, n_{\Lambda^s}(1)! \,
         n_{\Lambda^s}(2) ! \cdots \, ,
\end{equation}
    where $n_{\Lambda^s}(i)$ indicates  the number of $i$'s in
$\Lambda^s$,
         the symmetric part of $\Lambda=(\Lambda^a;\Lambda^s)$.
Moreover, $\mc{K}_{\sigma}$ stands for $\mc{K}_{i_1,i_1+1} \cdots
\mc{K}_{i_n,i_n+1}$ when the element $\sigma$ of the symmetric group
$S_N$ is written  in terms of elementary transpositions, i.e.,
$\sigma = \sigma_{i_1} \cdots \sigma_{i_n}$.  Notice that a
symmetric supermonomial $m_\La$, with $\La\vdash(n|m)$,  belongs to
$\mathscr{P}_{(n|m)}(\mathbb{Z)}$, the module of superpolynomials of
degree $(n|m)$ with integer coefficients.

\begin{theorem}
The set $\{m_\La\}_{\La\vdash(n|m)}:=\{m_\La \,:\,
\La\in\spar(n|m)\}$ is a basis of
$\mathscr{P}^{S_N}_{(n|m)}(\mathbb{Z})$.
\end{theorem}
\begin{proof}
Each superpolynomial $f(x,\theta)$ of degree
 $(n|m)$, with $N$ variables and with integer coefficients, can be
 expressed as a sum of monomials of the type
 $ \theta_{j_1} \dots \theta_{j_m} x^{\mu}$,
 with coefficient $a_{\mu}^{j_1,\dots,j_m} \in \mathbb Z$,
 and where $\mu$ is a composition of $N$.  Let $\Omega^a$ be the
 reordering of the entries
 $({\mu_{j_1},\dots,\mu_{j_m}})$, and let
 $\Om^s$ be the reordering of the remaining entries of $\mu$.
 Because the superpolynomial $f(x,\theta)$ is also symmetric,
 $f(x,\theta)$
 is by definition invariant under the action of $\mathcal K_{\sigma}$,
 for any $\sigma \in S_N$.   Therefore,
 $a_{\mu}^{j_1,\dots,j_m}$
 must be equal, up to a sign, to the coefficient
 $a_{\Omega}^{1,2,\dots,m}$ of
 $\theta_1 \cdots \theta_m x^{\Om}$ in $f(x,\theta)$,
 where $\Om=(\Om^a;\Om^s)$.
 Note that  from Lemma~\ref{antisymsym}, $\Omega^a$ needs to have
 distinct parts,
 which means that $\Omega$ is a superpartition.
 This gives that $f(x,\theta)-a_{\Omega}^{1,2,\dots,m} m_{\Omega}$
 does not contain any monomial that is also a monomial of $m_{\Om}$,
 since otherwise
 it would need by symmetry to contain the monomial $\theta_1 \cdots
 \theta_m x^{\Om}$.

 Now, consider any total order on superpartitions, and let $\La$ be
 the highest superpartition in this order such that there is a
 monomial of $m_{\La}$ appearing in $f(x,\theta)$.
   By the previous argument, $f(x,\theta)- a_{\Lambda}^{1,\dots,m}m_{\La}$ is
 a symmetric superpolynomial such that no monomial belonging to
 $m_{\La}$ appears
 in its expansion.  Since no monomial of $m_{\La}$ appears in any other
 monomial of $m_{\Om}$, for $\Om\neq \La$, the proof follows by induction.
\end{proof}

\begin{corollary} The set $\{m_\La\}_{\La}:=\{m_\La \,:\,
\La\in\spar\}$ is a basis of $\mathscr{P}^{S_N}(\mathbb{Z})$.
\end{corollary}

This corollary implies that $\mathscr{P}^{S_N}(\mathbb{Z})$ could
also be defined  as the free $\mathbb{Z}$-module spanned by the set
of symmetric supermonomials.

To end this section, we give a
formula for the expansion coefficients of the product of two supermonomials
 in terms of supermonomials.
This furnishes an illustrative example of  what
could be called supercombinatorics.  In this kind of calculation, the
standard counting of combinatorial objects (e.g., tableaux) is
affected by signs resulting from the reordering of fermionic
variables (represented by circles in the supertableaux).

\begin{definition} \label{deffilling} Let $\La \in \spar(n|m)$,
$\Omega \in \spar(n'|m')$ and $\Gamma \in \spar(n+n'|m+m')$. In each
box or circle of $D[\La]$ , we write a letter $a$. In its $i$-{th}
circle (the one corresponding to $\Lambda_i$), we add the label $i$ to
the letter $a$.
We do the same process for  $D[\Omega]$ replacing $a$ by $b$. We
then define $\mathcal T[ \La,\Om;\Gamma]$ to be the set of distinct
fillings of $D[\Gamma]$ with the letters of $D[\La]$ and $D[\Omega]$
obeying the following rules:
\begin{enumerate}
\item the circles of $D[\Gamma]$ can only be filled with labeled
letters (an $a_i$ or a $b_i$);
\item each row of the filling of $D[\La]$ is reproduced
in a single and distinct row of the filling of $D[\Gamma]$; in other words,
rows of $D[\La]$ cannot be split
and two rows of $D[\La]$ cannot be put within a single row of $D[\Gamma]$;
\item rule 2 also holds when $D[\La]$ is replaced by $D[\Omega]$;
\item in each row, the unlabeled $a$'s appear to the left of the
unlabeled $b$'s.
\end{enumerate}
\end{definition}

 For instance, there are three possible fillings of
$(2,1,0;1^3)$ with $(1,0;1)$ and
$(0;2,1^2)$:
\begin{equation}\label{exfilling}
{\tableau[scY]{b&b&\bl\tcercle{$a_2$}\\a&\bl\tcercle{$a_1$}\\a\\b\\b\\\bl\tcercle{$b_1$}}}
\qquad
{\tableau[scY]{b&b&\bl\tcercle{$a_2$}\\a&\bl\tcercle{$a_1$}\\b\\a\\b\\\bl\tcercle{$b_1$}}}
\qquad
{\tableau[scY]{b&b&\bl\tcercle{$a_2$}\\a&\bl\tcercle{$a_1$}\\b\\b\\a\\\bl\tcercle{$b_1$}}}\,
.
  \end{equation}\label{exone}
    There are also three possible fillings of
$ (3,1,0;1^2)$ with $(1,0;1)$ and $(0;2,1^2)$:
\begin{equation}\label{exfilling2}
{\tableau[scY]{a&b&b&\bl\tcercle{$a_1$}\\b&\bl\tcercle{$a_2$}\\a\\b\\\bl\tcercle{$b_1$}}}
\qquad
{\tableau[scY]{a&b&b&\bl\tcercle{$a_1$}\\b&\bl\tcercle{$a_2$}\\b\\a\\\bl\tcercle{$b_1$}}}
\qquad
{\tableau[scY]{a&b&b&\bl\tcercle{$a_1$}\\a&\bl\tcercle{$b_1$}\\b\\b\\\bl\tcercle{$a_2$}}}
\,.\end{equation}

\begin{definition}\label{defewightfilling} Let
$T \in  \mathcal T[ \La,\Om;\Gamma]$, with
$\overline{\underline{\Lambda}}=m$ and
$\overline{\underline{\Om}}=m'$. The weight of $T$, denoted by
$\hat{w}[T]$, corresponds to the sign of the permutation needed to
reorder the content of the circles in the filling of $D[\Gamma]$ so
that from top to bottom they read as $a_1 \dots a_m b_1 \dots
b_{m'}$.
\end{definition}

In the example (\ref{exfilling}), each term has weight $-1$ (odd parity). The oddness of these fillings comes from
the transposition
that is needed to reorder $a_1$ and $a_2$. In the second example (\ref{exone}),
the two first fillings are even while the last
filling is odd due to the needed transposition of $a_2$ and $b_1$.
As we shall see in the next proposition, the two previous sets lead
respectively to the coefficients of $m_{(2,1,0;1^3)}$ and
$m_{(3,1,0;1^2)}$ in the product of $m_{(1,0;1)}$ and
$m_{(0;2,1^2)}$, that is,
\begin{equation}
m_{(1,0;1)}\,m_{(0;2,1^2)}=\underbrace{(-3\,)}_{-1-1-1}\times\,m_{(2,1,0;1^3)}+\underbrace{(\,1\,)}_{+1+1-1}\times\,m_{(3,1,0;1^2)}+\mbox{other
terms}\, .\end{equation}
\begin{proposition}
Let $m_\La$ and $m_\Om$ be any two supermonomials.  Then
\begin{equation}\label{propmono}m_\La \,m_\Om
=\sum_{\Gamma\in\mathrm{SPar}}
N^\Gamma_{\La,\Om}\,m_\Gamma\, ,\end{equation}
where the integer $N^\Gamma_{\La,\Om}=(-1)^{
\overline{\underline{\Lambda}}\, \cdot\, \overline{\underline{\Om}}}
\, N^\Gamma_{\Om,\La}$ is given by
\begin{equation}N^\Gamma_{\La,\Om}:= \sum_{T \in  \mathcal T[
\La,\Om;\Gamma] } \hat w[T]\, .\end{equation}
\end{proposition}
  \begin{proof}
 From the symmetry property in Definition~\ref{defsupermono}, the
coefficient $N^\Gamma_{\La,\Om}$ is simply given by the coefficient
of $\theta_{\{1,\ldots,m+p\}}x^\Gamma$ in $m_\Lambda m_{\Om}$.  The terms contributing
to this coefficient correspond to all distinct permutations $\sigma$
and $w$ of the entries of $\Lambda$ and $\Omega$ respectively such
that
\begin{equation} \label{proofmono2}
\Gamma=
(\Lambda_{\sigma(1)}+\Omega_{w(1)},\dots,\Lambda_{\sigma(N)}+\Omega_{w(N)})
\, ,
\end{equation}
where the entries of $\Lambda^a$ and $\Omega^a$ are distributed
among the first $m+p$ entries (no two in the same position). But
this set is easily seen to be in correspondence with the fillings in
$\mathcal T[ \La,\Om;\Gamma]$ when realizing that labeled letters
simply give the positions of the fermions in $C[\Gamma]$ (the
circled version of $\Gamma$).
The only remaining problem is thus the ordering of the fermions.  In
(\ref{proofmono2}), from the definition of monomial symmetric
functions, the sign of the contribution is equal to the sign of the
permutation needed to reorder the fermionic entries of $\Lambda^a$
and $\Omega^a$ that are distributed among the first $m+p$ entries so
that they correspond to $({\La^a}, \Om^a)$.  But this is simply the
sign of the permutation that reorders the circled entries in the
corresponding filling of $D[\Gamma]$ such that they read as $a_1
\dots a_m b_1 \dots b_p$.
\end{proof}

\section{Generating functions and multiplicative bases}

In the theory of symmetric functions, the number of variables is
usually irrelevant, and can be set for convenience
to be equal to infinity.  In a similar way,
we shall consider from now on that, unless otherwise specified, the number of $x$
and $\theta$ variables is infinite, and denote the ring of symmetric superfunctions
as $\mathscr{P}^{S_\infty}$.

\subsection{Elementary symmetric superfunctions}

Let $J=\{j_1,\ldots,j_r\}$ with $1\leq j_1<j_2 <j_3 \cdots$ and let
$\# J:=\mbox{Card} \, J$.  The  $n$-th bosonic and fermionic elementary
symmetric superfunctions, for $n \geq 1$,  are defined respectively by
\begin{equation}e_n:=\sum_{J;\,\# J=n}x_{j_1}\cdots
x_{j_n}\qquad {\rm and}\qquad
\tilde{e}_n:=\sum_{i\geq 1}\sum_{\substack{J;\,\# J=n\\ i\not\in
J}}\theta_i\,x_{j_1}\cdots x_{j_n}\end{equation}
In addition, we impose \begin{equation}e_0=1\;\quad {\rm and }
\quad \tilde{e}_0=\sum_i\theta_i\, .
\end{equation}
     So, in terms of supermonomials, we have
\begin{equation}e_n=m_{(1^n)}\;, \qquad \tilde{e}_n=m_{(0;1^n)}\; .
\end{equation}

We introduce two parameters: $t\in\mathscr{B}$ and
$\tau\in\mathscr{F}$.  It is easy to verify that the generating
function for the elementary superfunctions is
\begin{equation}\label{defE} E(t,\tau):=\sum_{n=0}^\infty t^n(e_n
+\tau\tilde{e}_n) =\prod_{i=1}^\infty (1+t x_i + \tau \theta_i)\,
.\end{equation} Actually, to go from the usual generating function
$E(t):=E(t,0)$ to the new one, one simply replaces $x_i\rw
x_i+\tau'\ta_i$ and redefines $\tau'=\tau/t$, an operation that
makes manifest the invariance of $E(t,\tau)$ under the simultaneous
interchange of the $x_i$'s and the $\ta_i$'s.

   From an analytic point of view, the fermionic elementary
superfunctions are obtained by exterior differentiation:
\begin{equation}\tilde{e}_{n-1}(x,\theta)\sim\tilde{e}_{n-1}(x,dx)=d\,e_n(x)\,
,\end{equation}for
all $n \geq 1$.  How can we explain that the generating function (\ref{defE})
leads precisely to the  fermionic elementary superfunctions that are
obtained by the action of the exterior derivative of the elementary
symmetric function?  The rationale for this feature turns out to be
rather simple.  Indeed, let $\tau:= \, t \, dt$ and define $D$ to act on
a function $f(x,t)$ as a tensor-product  derivative:
\begin{equation} D f:= dt \wedge df \, .\end{equation}
In consequence, we formally have
\begin{equation}(1+t x_i +\tau \theta_i)\sim(1+D)(1+tx_i) \aand
E(t,\tau)\sim(1+D)\,
E(t)\,,\end{equation}
which is the desired link.

In order to obtain a new basis of the symmetric superpolynomial
algebra, we associate, to each superpartition
$\La=(\La_1,\ldots,\La_m;\La_{m+1},\ldots, \La_\ell)$ of $(n|m)$, a
superpolynomial $e_\La\in\mathscr{P}^{S_\infty}_{(n|m)}$ defined by
\begin{equation}e_\La:=\prod_{i=1}^m\tilde{e}_{\La_{i}}\prod_{j=m+1}^\ell
{e}_{\La_{j}}\, , \end{equation} Note that the product of
anticommutative quantities is always done from left to right:
$\prod_{i=1}^N F_i:=F_1F_2\cdots F_N$.
We stress that the ordering matters in the fermionic sector since for instance
\begin{equation}
{e}_{(3,0;4,1)}= {\tilde e}_{3}{\tilde e}_{0}{e}_{4}{e}_{1}= - {\tilde e}_{0}{\tilde e}_{3}{e}_{4}{e}_{1}
 \end{equation}

\begin{theorem} \label{theoebase}Let $\La$ be a superpartition of
$(n|m)$ and $\La'$ its conjugate. Then \begin{equation}\label{einm}
\versg{ e_\La}=\tilde{e}_{\La_m}\cdots
\tilde{e}_{\La_1}{e}_{\La_{m+1}}\cdots{e}_{\La_{N}}=m_{\La'}+\sum_{\Om
< \La'} N_{\La}^{\Om}\,m_\Om\, ,\end{equation} where $N_{\La}^{\Om}$
is an integer. Hence, $\{\,e_\La \,:\, \La\vdash(n|m)\}$ is a basis
of $\mathscr{P}^{S_\infty}_{(n|m)}(\mathbb{Z})$.\end{theorem}
\begin{proof}
We first observe that $\versg{ e_\La}=(-1)^{m(m-1)/2} e_{C[\La]}$,
where $C[\La]$ denotes as usual the partition $\La^*$ in which
fermionic parts of $\La$ are identified by a circle. Then, assuming
that we work in $N$ variables, the monomials $\theta_J x^\nu$ that
appear in the expansion of $e_{C[\La]}$ are in correspondence with
the fillings of $D[\Lambda']$ with the letters $1,\dots,N$ such
that:
\begin{enumerate}
\item the non-circled entries in the filling of $D[\Lambda']$
increase when going
       down in a column;
\item if a column contains a circle, then the entry that fills the
circle cannot
       appear anywhere else in the column.
\end{enumerate}
The correspondence follows because the reading of the $i$-th column
corresponds to one monomial of $e_{\Lambda_i}$ (or $\tilde
e_{\Lambda_i}$). To be more specific, if the reading of the column
is $j_1,\dots,j_{\Lambda_i}$ (with a possible  extra letter $a$), it
corresponds to the monomial $x_{j_1} \cdots x_{j_{\Lambda_i}}$ (or
$\theta_a x_{j_1} \cdots x_{j_{\Lambda_i}}$). The first condition
ensures that we do not count the permutations of $x_{j_1} \cdots
x_{j_{\Lambda_i}}$ as distinct monomials.  The second one ensures
that in the fermionic case, the index of $\theta_a$ is distinct from
the index of the variables $x_{j_1},\dots,x_{j_{\Lambda_i}}$.

Now, to obtain the coefficient $N^{\Omega}_\Lambda$, it suffices to
compute the coefficient of $\theta_{C[\Omega]} x^{C[\Omega]}$ in
$e_{C[\Lambda]}$, where $\theta_{C[\Omega]}$ represents the product
of the fermionic entries of $C[\Omega]$ read from left to right.
Note that this coefficient has the same sign as
$\theta_{\{1,\dots,m\}}x^{\Omega}$ in $m_{\Om}$ and there is thus no
need to compensate by a sign factor. The monomials contributing to
$N^{\Omega}_\Lambda$ are therefore fillings of $D[\Lambda']$
(obeying the two conditions given above) with the letter $i$
appearing $C[\Omega]_i$ times in non-circled cells with one
additional time in a circled cell if $C[\Omega]_i$ is fermionic.  We
will call the set of such fillings $\mathcal
T^{(e)}[\Omega;\Lambda']$.

Finally, given a filling $T \in \mathcal T^{(e)}[\Omega;\Lambda']$,
we read the content of the circles from top to bottom and obtain a
word $a_1\dots a_m$. The sign of the permutation needed to reorder
this word such that it be {\it increasing} gives the weight
associated to the filling $T$, denoted this time $\bar{w}[T]$.  The
weight of $T$ is the sign needed to reorder the monomial associated
to $T$ so that it corresponds
  to   $\theta_{C[\Omega]} x^{C[\Omega]}$ up to a factor $(-1)^{m(m-1)/2}$.
This is because to coincide with the product in $e_{C[\Lambda]}$
being done columnwise, we would have to read from bottom to top.
Reading from top to bottom provides the $(-1)^{m(m-1)/2}$ factor
needed to obtain the coefficient in $\versg{e_\La}$ instead of in
$e_{C[\La]}$. We thus have \begin{equation}N^{\Om}_{\La}=\sum_{T \in
\mathcal T^{(e)}[\Omega;\Lambda']} \bar{w}[T]\, .\end{equation} We
now use this equation to prove the theorem.

First, it is easy to convince ourselves that there is only one
element in $\mathcal T^{(e)}[\Lambda';\Lambda']$ and that it has a
positive weight. Because the rows of $D[\La']$ and $C[\La']$
coincide, for the filling to have increasing rows, we have no choice
but to put the $C[\La']_i$ letters $i$ in the $i$-th row of
$D[\La']$.  In the case when $C[\La']_i$ is fermionic, the extra $i$
has no choice but to go in the circle in row $i$ of $D[\La']$ for no
two $i$'s to be in a same column.  For example, given
$\Lambda'=(3,1;2,1)$ filling $D[\La']$ with the letters of $C[\La']$
leads to:
\begin{equation}{ \tableau[scY]{1&1&1&\bl\tcercle{1}\\2&2\\3&
\bl\tcercle{3}\\4}}
\quad\mbox{in}\quad{\tableau[scY]{&&&\bl\tcercle{}\\&\\&
\bl\tcercle{}\\&\bl}}\quad \longrightarrow\quad
{\tableau[scY]{1&1&1&\bl\tcercle{1}\\2&2\\3& \bl\tcercle{3}\\4}}
\end{equation}

This explains the first term in (\ref{einm}).

Second, let $\om=\Om^*$  and $\la=\La^*$. If $\Om \not <_S\La'$, a
filling of $\Om$ by $\La$ is obviously impossible because we would
need to be able to obtain in particular (forgetting about the
circles) a filling of the type $\mathcal T^{(e)}[\om;\la']$, which
would contradict the well known fact that the theorem holds in the
zero-fermion case.

Finally, we suppose that $\La^*=\Om^*$ and $\Om \not <_T\La'$. From
Lemma~\ref{lemmaST}, this implies that there is at least one circle
in $D[\La']$, let's say in row $i$, lower than its counterpart in
$D[\Om]$. Since to fill $D[\La']$ with $C[\Om]$, the non-circled
entries are filled with a row of 1's, then a row of 2's and so on,
we would need to be able to put an entry $j<i$ in the circle in row
$i$ of $D[\La']$. But this cannot happen since it
  would create a column with two $j$'s.
\end{proof}
\n Note that for the various examples that we have worked out, the
coefficients  $N^{\Om}_{\La}$ are non-negative.  So we may surmise
that a stronger version of the theorem, where $N^{\Om}_{\La}$ is
a non-negative integer, holds.

    The linear independence of the $e_\La$'s in
$\mathscr{P}^{S_N}$ implies that the first $N$ bosonic and
fermionic elementary superfunctions are algebraically independent
over $\mathbb{Z}$. Symbolically, \begin{equation}\label{fondtheo}
\mathbb{Z}[x_1,\dots,x_N,\theta_1,\ldots,\theta_N]^{S_N}\equiv
\mathbb{Z}[e_1,\dots,e_N,\tilde{e}_0,\dots\tilde{e}_{N-1}]\,
,\end{equation} which can be interpreted as the fundamental theorem
of symmetric superpolynomials.

\subsection{Complete symmetric superfunctions and involution}

The $n$-th bosonic and fermionic {\it complete symmetric
superfunctions} are given respectively by
\begin{equation}\label{bhn} h_n:=\sum_{\lambda \vdash n}\,m_\lambda
\qquad {\rm and}\qquad
\tilde{h}_{n}:=\sum_{\Lambda \vdash(n|1)}(\La_1+1)\,m_\La\,,\end{equation}
   From the explicit form of $h_n(x)$, namely,
$\sum_{1\leq i_1 \leq i_2 \leq \ldots \leq i_n}x_{i_1}\cdots
x_{i_n}$,
we see that its fermionic partner is again generated by the action
of $d$ in the form-representation:
\begin{equation}\tilde{h}_{n-1}(x,\theta)\sim
\tilde{h}_{n-1}(x,dx)=d\,h_n(x)\quad\mbox{for all}\quad n\geq 1\, .\end{equation}

The generating function for complete symmetric superpolynomials is
\begin{equation}\label{defH} H(t,\tau):=\sum_{n=0}^\infty
t^n(h_n +\tau\tilde{h}_n )=\prod_{i=1}^\infty \frac{1}{1-t x_i -
\tau \theta_i}\,.
\end{equation}
To prove (\ref{defH}), one simply uses
the inversion of even elements in the Grassmann algebra, given in
(\ref{invGrass}), which gives
  \begin{equation}\frac{1}{1-t x_i
- \tau \theta_i}= \sum_{n\geq 0}\left(t x_i+\tau \theta_i\right)^n=
\sum_{n\geq 0}\left[(t x_i)^n+n \tau \theta_i (tx_i)^{n-1}\right]\,
.\end{equation}

   From relations (\ref{defE}) and (\ref{defH}), we get
\begin{equation}H(t,\tau)\,E(-t,-\tau)=1\,.\end{equation} By
expanding the generating
functions in terms of  $e_n$, $\tilde{e}_n$, $h_n$ and $\tilde{h}_n$
in the last equation, we obtain recursion relations, of which the
non-fermionic one is a well known formula.

\begin{lemma}\label{recureh} Let $n\geq 1$, then
\begin{equation} \sum_{r=0}^n(-1)^r e_r h_{n-r}=0\, .\end{equation}
Let $n \geq 0$, then\begin{equation} \sum_{r=0}^n(-1)^r
(e_r\tilde{h}_{n-r}-\tilde{e}_r h_{n-r})=0\,
.\end{equation}\end{lemma}

\n Note that the second relation can be obtained from the first one
by the action of $d$ (representing, as usual, $\ta_i$ as $dx_i$).

We consider a homomorphism  $\hat{\omega}\,
:\,\mathscr{P}^{S_\infty}(\mathbb{Z})\rightarrow
\mathscr{P}^{S_\infty}(\mathbb{Z})$ defined by the following relations:
    \begin{equation}\label{definvolution}
    \hat{\omega} \,:\,
e_n\longmapsto h_n \quad\mbox{and}\quad
\tilde{e}_n\longmapsto\tilde{h}_n\, .\end{equation}
\begin{theorem}\label{involution}The homomorphism  $\hat{\omega}$ is
an involution, i.e., $\hat{\omega}^2=1$.  Equivalently, we have
\begin{equation}  \hat{\omega} \,:\, h_n\longmapsto e_n \quad\mbox{and}\quad
\tilde{h}_n\longmapsto\tilde{e}_n\,.\end{equation}\end{theorem}
\begin{proof} This comes from the application of transformation
(\ref{definvolution}) to the recursions appearing in
Lemma~\ref{recureh} followed by the comparison with the original
recursions.  Explicitly:
    \begin{equation}0=\sum_{r=0}^n(-1)^r
\hat{\omega}(e_r) \hat{\omega}(h_{n-r})=(-1)^n\sum_{r=0}^n(-1)^r
\hat{\omega}(h_{r})h_{n-r}\, ,\end{equation} which implies
$\hat{\omega}(h_{r})=e_r$.
Similarly, we have
\begin{equation}0=\sum_{r=0}^n(-1)^r
\left(\hat{\omega}(e_r)\hat{\omega}(\tilde{h}_{n-r})-\hat{\omega}(\tilde{e}_r)
\hat{\omega}(h_{n-r})\right)=(-1)^{n-1}\sum_{r=0}^n(-1)^r
\left( e_r\tilde{h}_{n-r}-\hat{\omega}(\tilde{h}_{r})h_{n-r}\right)\,
,\end{equation} leading to
$\hat{\omega}(\tilde{h}_{r})=\tilde{e}_{r}$.\end{proof}

Now, let
\begin{equation}h_\La:=\prod_{i=1}^{\underline{\overline{\La}}}\tilde{h}_{\La_{i}}\prod_{j=\underline{\overline{\La}}+1}^{\ell(\La)}
{h}_{\La_{j}}\, . \end{equation} Equation (\ref{definvolution}) and
Theorem \ref{involution} immediately give a bijection between two
sets of multiplicative superpolynomials:
\begin{equation}\hat{\omega}(e_\La)=h_\La\quad\mbox{and}\quad\hat{\omega}(h_\La)=e_\La\,
. \end{equation} We have thus obtained another $\mathbb{Z}$-basis
for the algebra of symmetric superpolynomials.
\begin{corollary} The set $\{\,h_\La \,:\, \La\vdash(n|m)\}$ is a basis of
$\mathscr{P}^{S_\infty}_{(n|m)}(\mathbb{Z})$.\end{corollary}

Finally, Lemma \ref{recureh} allows us to write determinantal
expressions for the elementary symmetric superfunctions in terms of
the complete symmetric superfunctions and {\it vice versa} using
the homomorphism $\hat \omega$.
\begin{proposition} \label{deteh} For $n\geq 1$, we have
\begin{equation}e_n=\left|%
\begin{array}{cccccc}
{h}_1  & {h}_2 & {h}_{3}&\cdots&h_{n-1}&{h}_n  \\
    1 & h_1 & h_2 & \cdots & h_{n-2}&h_{n-1} \\
    0 & 1 & h_1 & \cdots &h_{n-3}&h_{n-2}\\
    0&0&1&\ddots&h_{n-4}&h_{n-3}\\
     \vdots & \vdots & \ddots & \ddots& \ddots&\vdots \\
     0 & 0 & 0 & \cdots &1&h_1 \\
\end{array}%
\right|.\end{equation}For $n\geq 0$, we have
\begin{equation}\label{deuxiemedet}
\tilde{e}_n=\frac{1}{n!}\left|%
\begin{array}{cccccc}
     \tilde{h}_0 & \tilde{h}_1  & \tilde{h}_2 & \cdots
&\tilde{h}_{n-1}&\tilde{h}_n  \\
    n & (n+1)h_1 & (n+2)h_2 & \cdots & (2n-1) h_{n-1}&2n h_{n} \\
    0 & n-1 & n h_1 & \cdots &(2n-3) h_{n-2} &(2n-2)h_{n-1} \\
    0&0&n-2&\ddots&(2n-5) h_{n-3}&(2n-4)h_{n-2} \\
     \vdots & \vdots & \ddots & \ddots&\ddots&\vdots \\
     0 & 0 & 0 & \cdots & 1&2h_1 \\
\end{array}%
\right|.\end{equation}\end{proposition} \begin{proof} The first
relation is well known to be a simple application of Cramer's rule
to the linear system coming from Lemma~\ref{recureh}: $\mathbf{h}
=\mathbf{e}  \, \mathbf{H}$, where
\begin{equation} \mathbf{h}^\mathrm{t}=
\left(\begin{array}{c}
   h_1 \\ h_2 \\h_3\\ \vdots \\ h_N \\
\end{array}\right)\quad
\mathbf{e}^\mathrm{t}=\left(\begin{array}{c}
   e_1 \\ e_2 \\e_3\\ \vdots \\ e_N \\
\end{array}\right)
\quad \mathbf{H}=\left(\begin{array}{rrrrr}
   1&h_1&h_2&h_3&\ldots \\
   0&-1&-h_1&-h_2&\ldots  \\
   0&0&1&h_1&\ddots \\
   0&0&0&-1&\ddots\\
   \vdots & \ddots& \ddots& \ddots&\ddots\\
\end{array}\right).\end{equation}

To obtain the other determinant, we use the second formula of
Lemma~\ref{recureh} to obtain the linear system:
\begin{equation}
\tilde{\mathbf{h}} \, \overline{\mathbf{H}}=\tilde{\mathbf{e} }\, \mathbf{H}\,
,\end{equation} where $\mathbf{H}$ is as given above, and where
\begin{equation} \tilde{\mathbf{h}}^\mathrm{t}=
\left(\begin{array}{c}
   \tilde h_0 \\ \tilde h_1 \\ \tilde h_2\\ \vdots \\ \tilde h_{N-1} \\
\end{array}\right)\quad
\tilde{\mathbf{e}}^\mathrm{t} =\left(\begin{array}{c}
   \tilde e_0 \\ \tilde e_1 \\ \tilde e_2\\ \vdots \\ \tilde e_{N-1} \\
\end{array}\right)  \quad
  \overline{\mathbf{H}}=\left(\begin{array}{rrrrr}
   1&-e_1&e_2&-e_3&\ldots \\
   0&1&-e_1&e_2&\ldots  \\
   0&0&1&-e_1&\ddots \\
   0&0&0&1&\ddots\\
   \vdots & \ddots& \ddots& \ddots&\ddots\\
\end{array}\right). \end{equation}
Using the coadjoint formula for the inverse of a matrix, and the
determinantal expression for $h_n$ obtained by applying the
homomorphism $\hat \omega$ on the determinant of $e_n$ given above,
it is not hard to see that the $(i,j)$-th component of the inverse
of $\overline{\mathbf{H}}$ is simply $h_{j-i}$.  We are thus led to
the matrix relation:
\begin{equation}\tilde{\mathbf{h}}=\tilde{\mathbf{e}} \, \widetilde{\mathbf{H}} \,
,\end{equation}
where
\begin{equation}\widetilde{\mathbf{H}}_{i,j}=(-1)^{i+1}\sum_{k\geq0}^{i-j}h_{i-j-k}h_k=
-\widetilde{\mathbf{H}}_{i+1,j+1} \, .\end{equation} If we set
$H_i=\sum_{k=0}^i h_{i-k} h_k$, the matrix $\widetilde{\mathbf{H}}$
  can be expressed in a convenient form as
\begin{equation} \widetilde{\mathbf{H}}=\left(\begin{array}{rrrrr}
   1&H_1&H_2&H_3&\ldots\\
   0&-1&-H_1&-H_2&\ldots\\
   0&0&1&H_1&\ldots \\
   0&0&0&-1&\ldots\\
  \vdots & \ddots& \ddots& \ddots&\ddots\\
\end{array}\right) .\end{equation}
Using Cramer's rule and then multiplying rows $2,3,\dots n$ of the
resulting determinant by $n,n-1,\dots,1$ respectively, we obtain
\begin{equation}\label{finaldet}
\tilde{e}_n=\frac{1}{n!}\left|%
\begin{array}{cccccc}
     \tilde{h}_0 & \tilde{h}_1  & \tilde{h}_2 & \cdots
&\tilde{h}_{n-1}&\tilde{h}_n  \\
    n & nH_1 & n H_2 & \cdots & n H_{n-1}&n  H_{n} \\
    0 & n-1 & (n-1) H_1 & \cdots &(n-1) H_{n-2} &(n-1)H_{n-1} \\
    0&0&n-2&\ddots&(n-2) H_{n-3}&(n-2)H_{n-2} \\
     \vdots & \vdots & \ddots & \ddots&\ddots&\vdots \\
     0 & 0 & 0 & \cdots & 1&H_1 \\
\end{array}%
\right|.\end{equation} We will finally show that by manipulating
determinant (\ref{deuxiemedet}),
we obtain determinant (\ref{finaldet}).
Let $R_i$ be
row $i$ of
determinant  (\ref{deuxiemedet}).
If $R_2 \to R_2+R_3 h_1+\dots+R_{n} h_{n-2}$ in this
determinant, the second row becomes that of determinant  (\ref{finaldet})
due to the simple identity (for $i=1,\dots, n$)
$$
nH_i = (n+i) h_i +(n+i-2) h_{i-1}h_1+ \cdots + (n+i-2i+2)h_1
h_{i-1}+(n+i-2i)h_i \, .
$$
Doing similar operations on the lower rows, the two
determinants are seen to coincide.
\end{proof}

\subsection{Superpower sums}
We define the $n$-th bosonic and fermionic {\it superpower sums} as
follows: \begin{equation}p_n:=\sum_{i=1}^\infty x_{i}^n=m_{(n)}
\qquad {\rm and}\qquad \tilde{p}_n:=\sum_{i=1}^\infty\theta_i
x_i^n=m_{(n;0)}\, . \end{equation} Note that this time we will set $p_0=0$.
Obviously,
\begin{equation}n\,\tilde{p}_{n-1}(x,\theta)\sim
n\,\tilde{p}_{n-1}(x,dx)=d\,p_n(x)\end{equation}for
all $n \geq 1$.

Proceeding as in the complete symmetric superfunctions case, we
introduce products of power sums:
\begin{equation}
p_\La:=\prod_{i=1}^{\underline{\overline{\La}}}\tilde{p}_{\La_{i}}\prod_{j=\underline{\overline{\La}}+1}^{\ell(\La)}
{p}_{\La_{j}}\, .\end{equation} Also, we find that the generating
function for superpower sums is
\begin{equation}\label{defP}P(t,\tau):=\sum_{n\geq 0}t^n\,p_n
+\tau\sum_{n\geq 0}(n+1)t^n \,\tilde{p}_n =\sum_{i=1}^\infty \frac{t x_i
+ \tau \theta_i}{1-t x_i - \tau \theta_i}\, .\end{equation} One
directly verifies
\begin{equation}H(t,\tau)\,P(t,\tau)=(t\partial_t+\tau\partial_\tau)H(t,\tau)\end{equation}
and\begin{equation}E(t,\tau)\,P(-t,-\tau)=-(t\partial_t+\tau\partial_\tau)E(t,\tau)\,
.\end{equation} These expressions  lead after some manipulations to
the following recursion relations.
\begin{lemma}\label{recurep} Let $n\geq 1$.  Then
\begin{equation} h_n=\sum_{r=1}^n p_r h_{n-r}\,
, \qquad  ne_n=\sum_{r=1}^n (-1)^{r+1}p_r e_{n-r}\, .
\end{equation}
Let $n\geq 0$, and recall that $p_0=0$.  Then  
\begin{equation} (n+1)\,
\tilde{h}_n=\sum_{r=0}^n \bigl[ \, p_r\tilde{h}_{n-r}+ (r+1)\tilde{p}_r
h_{n-r} \, \bigr]\, ,
\end{equation}
\begin{equation}(n+1)\,
\tilde{e}_n=\sum_{r=0}^n (-1)^{r+1}  \bigl[ \, 
p_r\tilde{e}_{n-r}- (r+1)\tilde{p}_r e_{n-r} \, \bigr]
\,
.\end{equation}\end{lemma}

\begin{theorem}Let $\hat{\omega}$ be the involution defined in
(\ref{definvolution}). Then, for $n>0$,
    \begin{equation}\hat{\omega}\,:\,
p_n\longmapsto(-1)^{n-1}p_n\quad\mbox{and}\quad
     \tilde{p}_{n-1}\longmapsto(-1)^{n-1}\tilde{p}_{n-1}\end{equation}
or, equivalently,\begin{equation}\label{involup}
\hat{\omega}(p_\Lambda)=\omega_\La\,p_\Lambda\quad\mbox{with}\quad
\omega_\La:=(-1)^{|\La|+\overline{\underline{\Lambda}}-{\ell}(\La)}\,
.\end{equation}
\end{theorem}
\begin{proof}  We use  Lemma \ref{recurep} and  proceed as in the
proof of  Theorem \ref{involution}.\end{proof}

\begin{proposition} \label{detpe}For $n\geq 1$, we have
\begin{equation}p_n=\left|%
\begin{array}{ccccc}
{e}_1  & 2{e}_2 & 3{e}_{3}&\cdots&n{e}_n  \\
    1 & e_1 & e_2 & \cdots &e_{n-1} \\
    0 & 1 & e_1 & \cdots &e_{n-2}\\
     \vdots & \ddots & \ddots& \ddots&\vdots \\
     0 & 0 & 0&&e_1 \\
\end{array}%
\right|, \quad n!\,e_n=\left|%
\begin{array}{ccccc}
{p}_1  & p_2 &\cdots&p_{n-1}&p_n  \\
    1 & p_1  & \cdots &p_{n-2}&p_{n-1} \\
    0 & 2& \ddots &p_{n-3}&p_{n-2}\\
     \vdots & \ddots & \ddots& \ddots &\vdots \\
     0 & 0 &\cdots&n-1&p_1 \\
\end{array}%
\right|.\end{equation}For $n\geq 0$, we have
\begin{equation}\tilde{p}_n=\left|%
\begin{array}{ccccc}
     \tilde{e}_0 & \phantom{n}\tilde{e}_1  & \phantom{n}\tilde{e}_2 &
\cdots&\phantom{n}\tilde{e}_n  \\
    1 &\phantom{n} e_1 & \phantom{n}e_2 & \cdots &e_{n} \\
    0 & \phantom{n}1 & \phantom{n}e_1 & \cdots &e_{n-1} \\
     \vdots &  \ddots & \ddots&\ddots&\vdots \\
     0 & \phantom{n}0 & \phantom{n}0 & &e_1 \\
\end{array}%
\right|
,\quad {n!\,}\tilde{e}_n=\left|%
\begin{array}{ccccc}
     \tilde{p}_0 & \phantom{n}\tilde{p}_1   &
\cdots&\tilde{p}_{n-1}&\tilde{p}_n  \\
    n & \phantom{n}p_1 &  \cdots &p_{n-1}&p_{n} \\
    0 & n-1 &  \ddots &p_{n-2}&p_{n-1} \\
     \vdots &  \ddots & \ddots&\ddots &\vdots \\
     0 & \phantom{n}0 & \cdots & 1&p_1 \\
\end{array}%
\right|.\end{equation} Similar formulas for the complete symmetric
superfunctions are obtained by using the involution
$\hat\omega$.\end{proposition}
    \begin{proof} The proof is similar to that of Proposition
\ref{deteh}. \end{proof}

The  explicit formulas presented in Proposition \ref{detpe}
establish the correspondence between the sets
$\{p_n,\tilde{p}_{n-1}\}$ and $\{e_n,\tilde{e}_{n-1}\}$.
This implies, in particular, that
$e_\La=\sum_\Om c_{\La\Om}\,p_\Om$ for uniquely determined
coefficients  $c_{\La\Om}\in\mathbb{Q}$. Note that $c_{\La \Om}$ is
not necessarily an integer since, for instance, $e_2=p_1^2/2-p_2$.
Theorem \ref{theoebase} and Proposition \ref{detpe} thus imply
the following result.

\begin{corollary}\label{theobasep} The set $\{\,p_\La :
\La\vdash(n|m)\}$ is a basis of
$\mathscr{P}^{S_N}_{(n|m)}(\mathbb{Q})$ if  either $n\leq N$ and
$m=0$ or $n<N$ and $m>0$.  In particular,
the set $\{\,p_\La :
\La\vdash(n|m)\}$ is a basis of
$\mathscr{P}^{S_\infty}_{(n|m)}(\mathbb{Q})$.
\end{corollary}

 The
power sums will play a fundamental role in the remainder of the
article. For this reason, we will consider, from now on, only
symmetric superpolynomials defined over the rational numbers (or any
greater field):
\begin{equation} \mathscr{P}^{S_\infty}:=  \mathscr{P}^{S_\infty}(\mathbb{Q})\,
.\end{equation}

\subsection{Orthogonality}
Let $n_\la(i)$ denote the number of parts  equal to $i$ in the
partition $\lambda$.  We introduce a bilinear form, $\LL\,|\,\RR :
\mathscr{P}^{S_\infty}\times
\mathscr{P}^{S_\infty}\rightarrow\mathbb{Q}$, defined by
\begin{equation}\LL \versg{p_\La} | \versd{p_\Om }\RR:=z_\La
\delta_{\La,\Om}\, , \quad {\rm for} \quad
z_\La:=z_{\La^s}=\prod_{k\geq
1}\left[\,k^{n_{\La^s}(k)}\,n_{\La^s}(k)!\,\right]\, .\end{equation}

\begin{proposition} Let $f$ and $g$ be superpolynomials in
$\mathscr{P}^{S_\infty}$.
    Then $\LL \versg{f}|\versd{g}\RR$ is a scalar product,
that is, in addition to bilinearity, we have
\begin{equation}\begin{array}{ll}  \LL
\versg{f}| \versd{g}\RR=\LL \versg{g}|\versd{f}\RR=\LL
\versd{g}|\versg{f}\RR&\mbox{(symmetry)}\\
\LL \versg{f}|\versd{f}\RR\,> \,0 \quad \forall\quad f\neq
0&\mbox{(positivity)}\, .\\
    \end{array} \end{equation}\end{proposition}
\begin{proof} The symmetry property is a consequence of Lemma
\ref{ordredestheta}.  The positivity of the scalar product is proved
as follows.  By definition  $z_\La>0$ and by virtue of Corollary
\ref{theobasep}, there is a unique decomposition $f=\sum_\La f_\La
p_\La$.  Therefore $\LL \versg{f} | \versd{f}\RR=\sum_\La f_\La^2
z_\La>0$. \hfill \end{proof}

\begin{proposition}The involution $\hat{\omega}$ is an
isometry.\end{proposition}
    \begin{proof} Given that
$\{p_\La\}_\La$ is a basis of $\mathscr{P}^{S_\infty}$, for any symmetric
polynomials $f$ and $g$, we have $f=\sum_\La f_\La p_\La$ and
$g=\sum_\La g_\La p_\La$. Thus

\begin{eqnarray}
\LL \versg{\hat{\omega}f} |
\versd{\hat{\omega}g}\RR&=&\sum_{\La,\Om}f_\La g_\Om \LL
\versg{\hat{\omega}p_\La} | \versd{\hat{\omega}p_\Om}\RR\cr
&=&\ds \sum_{\La} f_\La g_\Om
(-1)^{|\La|+\overline{\underline{\Lambda}}-{\ell}(\La)}
(-1)^{|\Om|+\overline{\underline{\Om}}-{\ell}(\Om)}\LL\versg{p_\La}|\versd{p_\Om}\RR
\cr
&=&\sum_{\La} z_\La f_\La g_\La=\LL \versg{f} | \versd{g}\RR\,
,
\end{eqnarray}as claimed.\end{proof}

The following theorem is of particular importance since it gives
Cauchy-type formulas for  the superpower sums.

\begin{theorem} \label{cauchypp} Let
$K=K(x,\theta;y,\phi)$   be the
bi-symmetric formal superfunction given
by\begin{equation}\label{defK} K:=\prod_{
i,j}\frac{1}{1-x_iy_j-\theta_i\phi_j}  \; . \end{equation}
    Then
\begin{equation} K=\sum_{\La\in\mathrm{SPar}}z_\La^{-1}
\versg{p_\La(x,\theta)}\versd{p_\La(y,\phi)} \;.\end{equation} \end{theorem}
    \begin{proof} We have:
\begin{eqnarray}
\prod_{i,j}\frac{1}{1-x_iy_j-\theta_i\phi_j}&=&\exp\Big\{\sum_{i,j}
\ln\Big[(1-x_iy_j-\theta_i\phi_j)^{-1}\Big]\Big\}
    =\exp\Big\{\sum_{i,j}\sum_{n\geq
    1}\Big[\frac{1}{n}(x_iy_j+\theta_i\phi_j)^{n}\Big]\Big\}\cr
    &=&\exp\Big\{\sum_{n\geq
1}\Big[\frac{1}{n}p_n(x)\,p_n(y)\Big]+\sum_{n\geq0}\Big[\tilde{p}_n(x,\theta)\,\tilde{p}_n(y,\phi)\Big]\Big\}\cr
    &=&\prod_{n\geq 1}\sum_{k_n\geq
    0}\frac{1}{n^{k_n}k_n!}\Big[p_n(x)\,p_n(y)\Big]^{k_n}\exp
\Big[\sum_{n\geq0}\tilde{p}_n(x,\theta)\,\tilde{p}_n(y,\phi)\Big]\,
    .
\end{eqnarray}
Considering Proposition \ref{expfermions}, we find
\begin{equation}\prod_{i,j}\frac{1}{1-x_iy_j-\theta_i\phi_j}=\sum_{n,\,m\geq
0}\sum_{\substack{\lambda\in \mathrm{Part}_s(n)\\ \mu\in
\mathrm{Part}_a(m)}} \Big[z_{\lambda}^{-1}p_\lambda(x)
p_\lambda(y)\versg{ \tilde{p}_\mu(x,\theta) }\versd{ \tilde{p}_\mu
(y,\phi) }\Big]\, .\end{equation} This equation (together with Lemma
\ref{ordredestheta}) proves the theorem. \end{proof}

\begin{remark}  The inverse of the kernel satisfies:
\begin{equation}
K(-x,-\theta;y,\phi)^{-1} =\prod_{i,j}(1+x_iy_j+\theta_i\phi_j) \, =\sum_{\La\in\mathrm{SPar}}\omega_\La z_\La^{-1}
\versg{p_\La(x,\theta)}\versd{p_\La(y,\phi)}\, .
\end{equation}
The proof
of this  result is similar to that of  Theorem \ref{cauchypp}, apart from the presence of the
coefficient
$\omega_\La=(-1)^{|\La|+\overline{\underline{\Lambda}}-{\ell}(\La)}$,
which comes from the expansion of $\ln(1+x_iy_j+\theta_i\phi_j)$.  This shows that
\begin{equation}\label{deuxK}
K(x,\theta;y,\phi) = \omega^{(x,\ta)}\, K(-x,-\theta;y,\phi)^{-1} = \omega^{(y,\phi)}\, K(-x,-\theta;y,\phi)^{-1}
\end{equation}
where $\omega^{(x,\ta)}$ indicates that $\omega$ acts on the $(x,\ta)$ variables and similarly for $\omega^{(y,\phi)}$.
\end{remark}

We now give two direct consequences of Theorem \ref{cauchypp}.
\begin{corollary} $K$ is a reproducing kernel
in the space of symmetric superfunctions:
\begin{equation}\LL\,K(x,\theta;y,\phi)\,|\,f(x,\ta\,\RR=f(y,\phi)\,
,\quad\mbox{for all}\quad f\in\mathscr{P}^{S_{\infty}}\,
.\end{equation}\end{corollary}

\begin{proof} If $f\in\mathscr{P}^{S_{\infty}}$, there exist unique
coefficients $f_\La$ such that $f=\sum_\La f_\La p_\La$. Hence,
\begin{eqnarray}
\LL\,K(x,\theta;y,\phi)\,|\,f(x,\ta)\,\RR
&=&\sum_{\Om,\La}z_\Om^{-1}f_\La\LL
\versg{p_\Om(x,\ta)}|\versd{p_\La(x,\ta)}\RR \versd{p_\Om(y,\phi)}\cr &=&\sum_\La f_\La
\versd{p_\La(y,\phi)}=  f(y,\phi)\, ,\end{eqnarray} as
desired.\end{proof}

\begin{corollary} \label{einbasep}We have
\begin{eqnarray} h_n=\sum_{\Lambda\in \footnotesize{\spar}
(n|0)}z_\La^{-1}p_\Lambda\aand
    e_n=\sum_{\Lambda\in \footnotesize{\spar} (n|0)}z_\La^{-1}\omega_\Lambda
    p_\Lambda\, ,\end{eqnarray}
\begin{equation} \tilde{h}_n=\sum_{\Lambda\in\footnotesize{ \spar}
(n|1)}z_\La^{-1}p_\Lambda\aand
    \tilde{e}_n=\sum_{\Lambda\in \footnotesize{\spar}
(n|1)}z_\La^{-1}\omega_\Lambda
    p_\Lambda\, .\end{equation}\end{corollary}
\begin{proof} Using the definition of the generating function
$E(t,\tau)$, we  first make the following correspondence:
\begin{equation}E(t,0)=\sum_{n\geq 0}t^ne_n(x)=K(-x,0;y,0)^{-1}\Big|_{y=(t,0,0,\ldots)}\;.
\end{equation}
Thus, from Theorem \ref{cauchypp} and
$p_\la(y)|_{y=(t,0,0,\ldots)}=t^{|\lambda|}$, we have
\begin{equation}\sum_{n\geq0}t^ne_n(x)=\sum_{\lambda\in\mathrm{Par}}t^{|\lambda|}\omega_\lambda
z_\lambda^{-1} p_\lambda(x)\quad\Longrightarrow\quad
e_n=\sum_{\lambda\vdash n}z_\la^{-1}\omega_\lambda
    p_\lambda\, .\end{equation} Then, we observe that
    \begin{equation}\partial_\tau E(t,\tau)=\sum_{n=0}^{N-1}t^n
    \tilde{e}_n(x,\tau)=\partial_\tau
    K(x,\theta;y,\phi)\Big|_{\substack{y=(t,0,0,\ldots)\\
    \phi=(-\tau,0,\ldots)}}\; .\end{equation} Theorem
    \ref{cauchypp} and
\begin{equation}
p_\La(y,\phi)\Big|_{\substack{y=(t,0,0,\ldots)\\
    \phi=(-\tau,0,\ldots)}}= \left \{
\begin{array}{ll} \phantom{-\tau}t^{|\La|} & {\rm ~if~} \overline{\underline{\La}}=0 \\
-\tau t^{|\La|} & {\rm ~if~} \overline{\underline{\La}}=1  \\
   \phantom{-\tau}0 & {\rm ~otherwise~}
\end{array} \right.\phantom{\{}
\end{equation}
finally lead to
    \begin{equation}\sum_{n\geq0}t^n
\tilde{e}_n(x,\tau)=\sum_{\Lambda,\overline{\underline{\La}}=1}z_\La^{-1}\omega_\Lambda
    p_\Lambda(x,\theta)\quad
    \Longrightarrow\quad \tilde{e}_n=\sum_{\Lambda\in
\mathrm{SPar}(n|1)}z_\La^{-1}\omega_\Lambda
    p_\Lambda\, .\end{equation}
    Note that the minus sign disappears since $\partial_\tau$ and
$p_\La(x,\theta)$ anticommute when $\overline{\underline{\La}}=1$.
    Similar formulas relating the superpower sums to the
homogeneous symmetric superpolynomials
    are obtained using the involution $\hat{\omega}$.\end{proof}

\begin{lemma}\label{cauchybases}Let $\{u_\La\}$ and $\{v_\La\}$ be two
bases of $\mathscr{P}^{S_{\infty}}_{(n|m)}$.  Then
\begin{equation}K(x,\theta;y,\phi)=\sum_\La
\versg{u_\La(x,\theta)}\,\versd{v_\La(y,\phi)}\quad\Longleftrightarrow\quad\LL
\versg{u_\La}|\versd{v_\La}\RR=\delta_{\La,\Om}\,
.\end{equation}\end{lemma}
\begin{proof}
The proof is identical to the one in the case without
Grassmannian variables (see \cite{Mac} I.4.6).\end{proof}

\begin{proposition} \label{propkernelmh} Let $K$ be the
superfunction
defined in (\ref{defK}). Then,
\begin{equation}K=\sum_{\La\in\mathrm{SPar}}
\versg{m_\La(x,\theta)}\versd{h_\La(y,\phi)} \, .\end{equation}
\end{proposition}

\begin{proof} We start with the definition of the generating
function $E(t,\tau)$:
\begin{eqnarray}
K(-x,-\ta;y,\phi)^{-1}&=&\prod_{ i,j}(1+x_iy_j+\theta_i\phi_j) =  \prod_{i }E(x_i, \ta_i) \cr
&=&  \prod_{i }\Big[\sum_{n\geq0}x_i^n\,e_n(y)
+\theta_i\sum_{n\geq0} x_i^n\tilde{e}_n(y,\phi)\Big]\cr
&=&\sum_{\epsilon_1,\epsilon_2,\dots \in \{0,1\}}
\sum_{n_1,n_2,\dots \geq 0} \prod_i \left(\theta_i^{\epsilon_i}
x_i^{n_i} e_{n_i}^{(\epsilon_i)}(y,\phi)  \right) \cr
&=&
\sum_{\La\in\mathrm{SPar}}\versg{m_\La(x,\theta)}\versd{e_\La(y,\phi)}\;.
    \end{eqnarray}
In the third line we have set $e_n^{(0)}(y,\phi)=e_n(y,\phi)$ and $e_n^{(1)}(y,\phi)=\tilde
e_n(y,\phi)$. The fourth line follows by reordering the variables using Lemma
\ref{ordredestheta}. Using (\ref{deuxK}), we can recover $K(x,\ta;y,\phi)$ by acting
with  $\hat\omega^{\{y,\phi
\}}$ on $ K(-x,-\ta;y,\phi)^{-1}$. The identity then
follows from $\hat\omega(e_{\Lambda})= h_{\Lambda}$.
\end{proof}

The previous proposition and Lemma \ref{cauchybases} have the
following corollary.
\begin{corollary}
The supermonomials are dual to the complete symmetric
superfunctions:
\begin{equation} \LL \versg{\,h_\La} | \versd{m_\Om
}\RR=\delta_{\La,\Om}\, .\end{equation}\end{corollary}

\subsection{One-parameter deformation of the
scalar product and the homogeneous basis}

In this section, we introduce a natural one-parameter -- called $\beta$ --  deformation of the
scalar product, of the corresponding superspace kernel as well as the deformation of the
homogeneous symmetric basis in superspace.

   Let
$\mathscr{P}^{S_\infty}(\beta)$  stand for the algebra of symmetric
superpolynomials with coefficients in $\mathbb{Q}(\beta)$, i.e.,
rational functions of $\beta$.
We now introduce a mapping,
\begin{equation}\LL\,\cdot \, |\, \cdot\,\RR_{\beta} \,:
\,\mathscr{P}^{S_\infty}(\beta)\times
\mathscr{P}^{S_\infty}(\beta)\longrightarrow\mathbb{Q}(\beta)\end{equation}
defined by \begin{equation} \label{defscalprodcomb} \LL
\versg{p_\La} | \versd{p_\Om }\RR_\beta:=z_\La
(\beta)\delta_{\La,\Om}\,,\quad\mbox{where}\quad
z_\La(\beta):=\beta^{-{\ell}(\La)}z_{\La}\, .\end{equation}
This bilinear
form can again be shown to be a scalar product.

We also introduce a homomorphism  that generalizes the involution
$\hat{\omega}$. It is defined on the power sums as:
\begin{equation}\label{defendobeta}\hat{\omega}_\alpha(p_n)=(-1)^{n-1}\,\alpha\,
p_n\quad\mbox{and}\quad
\hat{\omega}_\alpha(\tilde{p}_{n})=(-1)^{n}\,\alpha\,\tilde{p}_{n}
\, , \end{equation} where $\alpha$ is some unspecified parameter.
This implies
\begin{equation}\hat{\omega}_\alpha(p_\Lambda)=\omega_\La(\alpha)\,p_\Lambda
\quad\mbox{with} \quad \omega_\La(\alpha):= \alpha^{\ell(\Lambda)}
(-1)^{|\La|-\overline{\underline{\La}}+{\ell}(\La)}\,.\end{equation}
Notice that $\hat{\omega}_1\equiv \hat{\omega}$.  This
homomorphism is still self-adjoint, but it is now
neither an involution
($\hat{\omega}_\alpha^{-1}=\hat{\omega}_{\alpha^{-1}}$) nor an
isometry ( $\|\hat{\omega}_\alpha p_\La \|^2=z_\La(\beta/\alpha^2)$
).
Note also that
\begin{equation}z_\La(\beta)\omega_\La(\beta)=z_\La\omega_\La\aand
z_\Lambda (\beta)^{-1} \omega_\Lambda(\beta^{-1}) = z_\La^{-1}
\omega_\Lambda\, .\end{equation}

\begin{theorem} \label{cauchyppbeta} With $K^\beta$ given by
 \begin{equation}
K^\beta=\prod_{i,j}\frac{1}{(1-x_iy_j-\theta_i\phi_j)^{\beta}}\,,
\end{equation}
we have
\begin{equation}K^\beta=\sum_{\La}z_\La(\beta)^{-1}
\versg{p_\La(x,\theta)}\versd{p_\La(y,\phi)}\,.
\end{equation}
\end{theorem}
\begin{proof}
Starting from
\begin{equation}\prod_{i,j}\frac{1}{(1-x_iy_j-\theta_i\phi_j)^{\beta}}=\exp\Big\{\beta\sum_{i,j}
\ln\Big[(1-x_iy_j-\theta_i\phi_j)^{-1}\Big]\Big\} \, ,\end{equation}
The above identity can be obtained straightforwardly proceeding as
in the proof of Proposition~\ref{cauchypp}. \end{proof}

\begin{remark} The inverse of $K^\beta$ satisfies:
\begin{equation}
 K(-x,-\ta;y,\phi)^{-\beta}=\prod_{i,j}(1+x_iy_j+\theta_i\phi_j)^{\beta}\,
=\sum_{\La} z_\La(\beta)^{-1}\omega_\La
\versg{p_\La(x,\theta)}\versd{p_\La(y,\phi)} \; ,\end{equation}
which is obtained by using
\begin{equation}p_\La(-x,-\theta)=(-1)^{|\La|+\overline{\underline{\La}}}\,
p_\La(x,\theta)\aand
  z_\La(-\beta)=(-1)^{\ell(\La)}z_\La(\beta)\,.\end{equation}
Notice also
the simple relation between the  kernel and its $\beta$-deformation
\begin{equation}
K^{\beta}(x,\ta;y,\phi)=\hat{\omega}_\beta K(-x,-\ta;y,\phi)^{-1}
\, ,\end{equation} where it is
understood that $\hat{\omega}_\beta$ acts either on $(x,\theta)$ or
on $(y,\phi)$.
\end{remark}

\begin{corollary}
\label{Kernelbeta}  $K^\beta(x,\theta;y,\phi)$ is a reproducing
kernel in the space of symmetric superfunctions with rational
coefficients in $\beta$:
\begin{equation}\LL\,K^\beta(x,\theta;y,\phi)\,|\,f(x,\ta)\,\RR_\beta=f(y,\phi)\,
,\quad\mbox{for all}\quad f\in\mathscr{P}^{S_{\infty}}(\beta)\,
.\end{equation}\end{corollary}

We now introduce a $\beta$-deformation of the bosonic and fermionic
complete homogeneous symmetric functions, respectively denoted as
$g_n(x)$ and $\tilde{g}_n(x,\theta)$ (the $\beta$-dependence being implicit). Their
generating function is
\begin{equation}\label{generatriceg} G(t,\tau;\beta):=\sum_{n\geq
0}t^n[g_n(x)+\tau \tilde{g}_n(x,\theta)]=\prod_{i\geq
1}\frac{1}{(1-tx_i-\tau\theta_i)^\beta}\, .\end{equation} Clearly,
$g_n=h_n$ and $\tilde{g}_n=\tilde{h}_n$ when $\beta=1$.  As usual,
we define
\begin{equation}\label{defgla}
g_\La:=\prod_{i=1}^{\underline{\overline{\La}}}\tilde{g}_{\La_i}
\prod_{j=\underline{\overline{\La}}+1}^{\ell(\La)}g_{\La_j}\,
  .\end{equation}

\begin{proposition}\label{propkernelmg}We have  $\ds
K^\beta(x,\theta;y,\phi)=\sum_{\La} \versg{m_\La(x,\theta)}
\versd{g_\La(y,\phi)}$.
\end{proposition}
\begin{proof}We proceed as in the proof of Proposition~\ref{propkernelmh}:
\begin{eqnarray}
K^\beta&=&\prod_{ i } G(x_i,\theta_i;\beta)=  \prod_{i }\Big[\sum_{n\geq0}x_i^n g_n(y)+\theta_i
\tilde{g}_n(y,\phi)\Big]\cr
&=&\sum_{\epsilon_1,\epsilon_2,\dots \in \{0,1\}}
\sum_{n_1,n_2,\dots \geq 0} \prod_i \left(\theta_i^{\epsilon_i}
x_i^{n_i} g_{n_i}^{(\epsilon_i)}(y,\phi)  \right) \, ,
    \end{eqnarray}
where $g_n^{(0)}(y,\phi)=g_n(y,\phi)$ and $g_n^{(1)}(y,\phi)=\tilde
g_n(y,\phi)$. By reordering the variables using Lemma
\ref{ordredestheta}, we get the desired result.\end{proof}

\begin{corollary}We have
\begin{equation}\label{genp} g_n=\sum_{\La\vdash(n|0)}z_\La(\beta)^{-1} p_\La
\aand \tilde{g}_n=\sum_{\La\vdash(n|1)}z_\La(\beta)^{-1} p_\La\,
.\end{equation}
\end{corollary}
\begin{proof}On the one hand,
\begin{equation}G(t,0;\beta)=\sum_{n\geq 0}t^n
g_n(x)=K^\beta(x,0;y,0)\Big|_{y=(t,0,0,\ldots)}\;.\end{equation}
The previous proposition and Theorem~\ref{cauchyppbeta} imply
\begin{equation}\sum_{n\geq0}t^n g_n=\sum_{\la\in
\mathrm{Par}}t^{|\la|}z_\la(\beta)^{-1}p_\la\quad\Longrightarrow\quad
g_n=\sum_{\la\in \mathrm{Par}(n)}z_\la(\beta)^{-1}p_\la\,
.\end{equation} On the other hand, \begin{equation}\partial_\tau
G(t,\tau;\beta)=\sum_{n\geq 0}t^n
\tilde{g}_n(x,\theta)=K^\beta_+(x,\theta;y,\phi)\Big|_{\substack{y=(t,0,0,\ldots)\\
\phi=(-\tau,0,0,\ldots)}}\;. \end{equation}
 Hence
\begin{equation}\sum_{n\geq 0} t^n \tilde{g}_n=\sum_{\La,\overline{\underline{\La}}=1}t^{|\La|}z_\La(\beta)^{-1}p_\La\quad\Longrightarrow\quad
\tilde{g}_n=\sum_{\La\in \mathrm{SPar}(n|1)}z_\La(\beta)^{-1}p_\La\,
,\end{equation} as claimed.
\end{proof}

Using the previous corollary and the relation
$n\tilde{p}_{n-1}\sim dp_n$, it is easy to show that the fermionic
superfunction $\tilde{g}_{n-1}$ can be represented as the exterior
derivative of $g_n$: \begin{equation}\tilde{g}_{n-1}(x,\theta)\sim d
g_n(x)\quad\mbox{for}\quad n\geq 1\, .\end{equation} This is a
direct extension  of the  $\beta=1$ case.  Also,  applying
$\omega_{\beta^{-1}}$ on equation (\ref{genp}) and  comparing with
Corollary \ref{einbasep}, we get \begin{equation}\label{dualityge}
\hat{\omega}_{\beta^{-1}} g_n=e_n\aand  \hat{\omega}_{\beta^{-1}}
\tilde{g}_n=\tilde{e}_n\, .\end{equation}
  Lemma~\ref{cauchybases} can also be trivially
generalized.
\begin{equation} \label{eqKN}
K^\beta(x,\theta;y,\phi)=\sum_\La
\versg{u_\La(x,\theta)}\,\versd{v_\La(y,\phi)}\quad\Longleftrightarrow\quad\LL
\versg{u_\La}|\versd{v_\La}\RR_\beta=\delta_{\La,\Om}\,
.\end{equation} This equation, together with Proposition~\ref{propkernelmh},
immediately implies the following.
\begin{corollary}\label{corodualitygm}
The set $\{g_\La\}_\La$ constitutes a basis of
$\mathscr{P}^{S_\infty}(\beta)$ dual to that of the supermonomials,
that is,
\begin{equation} \LL \versg{\,g_\La} | \versd{m_\Om
}\RR_\beta=\delta_{\La,\Om}\, .\end{equation}\end{corollary}

We shall need in the next section to make explicit the distinction
between an infinite and a finite number of variables.   Therefore, we
also let
\begin{equation}\LL\,\cdot \, |\, \cdot\,\RR_{\beta,N} \,:
\,\mathscr{P}^{S_N}(\beta)\times
\mathscr{P}^{S_N}(\beta)\longrightarrow\mathbb{Q}(\beta)\end{equation}
be defined by requiring that the bases
$\{ g_\La\}_{\ell(\La)\leq N}$ and $\{ m_\La\}_{\ell(\La)\leq N}$ be dual to each other:
\begin{equation} \label{defscalprodcombN} \LL
\versg{g_\La} | \versd{m_\Om }\RR_{\beta,N}:=\delta_{\La,\Om}
\, ,\end{equation}
whenever $\ell(\La)$ and $\ell(\Om)$ are not larger than $N$.
From this definition, it is thus obvious that
\begin{equation} \label{degN}
\LL\,f^{(N)} \, |\, g^{(N)} \,\RR_{\beta,N}=\LL\,f \, |\, g\,\RR_{\beta}.
\end{equation}
if $f$ and $g$ are symmetric superpolynomials of degrees
not larger than $N$, and if $f^{(N)}$ and $g^{(N)}$ are their
respective restriction to $N$
variables.  This is because $f$ and $f^{(N)}$ (resp. $g$ and $g^{(N)}$ ) then
have the same expansion in terms of the $g$ and $m$ bases.  Note that
with this definition, we have that
\begin{equation}
K^{\beta,N} = \sum_{\ell(\La) \leq N} g_{\La}(x,\theta)\, m_{\La}(y,\phi) \, ,
\end{equation}
where $K^{\beta,N}$ is the restriction of $K^{\beta}$ to $N$ variables and
where $(x,\theta)$ stand for $(x_1,\dots,x_N;\theta_1,\dots,\theta_N)$
and $(y,\phi)$ for $(y_1,\dots,y_N;\phi_1,\dots,\phi_N)$.

We complete this section by displaying a relationship
between the $g$-basis elements and the bases of monomials and homogeneous superpolynomials.

\begin{proposition}\label{propgn}Let
$n_\La!:=n_{\La^s}(1)!\,n_{\La^s}(2)!\cdots$,  and
\begin{equation}\left(%
\begin{array}{c}
   \beta \\
   n \\
\end{array}%
\right) := \frac{\beta(\beta-1)\cdots(\beta-n+1)}{n!}
\;, \qquad (\beta)_n :=  \beta(\beta-1)\cdots(\beta-n+1)
.\end{equation}

Then
\begin{eqnarray}
\label{jnenm}g_n&=&
\sum_{\La\vdash (n|0)}\prod_i \left(%
\begin{array}{c}
   \beta+\La_i-1 \\
   \La_i \\
\end{array}
\right)m_\La\, =  \sum_{\La\vdash (n|0)}
\frac{(\beta)_{\ell(\La)}}{n_\La !} \,h_\La\, ,\\
\label{jnoenm}\tilde{g}_n&=&
\sum_{\La\vdash (n|1)}(\beta+\La_1)\prod_i \left(%
\begin{array}{c}
   \beta+\La_i-1 \\
   \La_i \\
\end{array}
\right)m_\La\, = \sum_{\La\vdash
(n|1)}\frac{(\beta)_{\ell(\La)}}{n_\La!}\,h_\La\, .
 \end{eqnarray}
\end{proposition}
\begin{proof}We start with  the generating
function (\ref{generatriceg}). The product on the right hand side
%
can also be written as
\begin{eqnarray}
\lefteqn{\prod_{i\geq 1} \sum_{k\geq
0}(-1)^k\left(%
\begin{array}{c}
   -\beta \\
   k \\
\end{array}%
\right)(tx_i+\tau\theta_i)^k =}\cr
&&\prod_{i\geq 1}\left[ \sum_{k\geq
0}\left(%
\begin{array}{c}
   \beta+k-1 \\
   k \\
\end{array}%
\right)(tx_i)^k+\tau \theta_i \sum_{k\geq
1}k\left(%
\begin{array}{c}
   \beta+k-1 \\
   k \\
\end{array}%
\right)(tx_i)^{k-1}\right]
\end{eqnarray}
After some easy manipulations, this  becomes
\begin{equation}\sum_{n\geq 0}t^n\left[\sum_{\la\vdash n} \prod_i\left(%
\begin{array}{c}
   \beta+\la_i-1 \\
   \la_i \\
\end{array}%
\right)m_\la\,+\tau \sum_{\La\vdash (n|1)}(\beta+\La_1)\prod_i\left(%
\begin{array}{c}
   \beta+\La_i-1 \\
   \La_i \\
\end{array}%
\right)\,m_\La\right]\end{equation} and  the first equality in the two formulas (\ref{jnenm}) and (\ref{jnoenm}) are
seen to hold.

To prove the remaining two formulas, we use the generating function
of the homogeneous symmetric functions and proceed as follows:
\begin{eqnarray}
\prod_i(1-tx_i-\tau\theta_i)^{-\beta} &=& \left( 1+\sum_{m\geq 1}t^m h_m+\tau \sum_{n\geq 0}
t^n\tilde{h}_n
\right)^\beta\cr
&=& \sum_{k\geq 0} \left(\begin{array}{c}
   \beta \\
   k \\
\end{array}\right)
\left( \sum_{m\geq 1}t^m h_m+\tau \sum_{n\geq 0} t^n\tilde{h}_n
\right)^k\cr
&=&\sum_{n\geq 0}\sum_{\la\vdash
n}t^n\frac{(\beta)_{\ell(\La)}}{n_\la !}\, h_\la +\tau \sum_{m\geq
0}\sum_{\la\vdash m}t^m\frac{(\beta)_{\ell(\La)+1}}{\la
!}\, h_\la\sum_{n\geq 0}t^n \tilde{h}_n\cr
&=&\sum_{n\geq 0}t^n\left[\sum_{\la\vdash
n}\frac{(\beta)_{\ell(\La)}}{n_\la !}\, h_\la +\tau \sum_{\La\vdash
(n|1)}\frac{(\beta)_{\ell(\La)}}{n_{\La^s} !}\, h_\La\right]
\end{eqnarray}from which the desired expressions can be obtained.
\end{proof}


\section{Jack polynomials in superspace}

  The main goal of this section is to show that there exists a natural supersymmetric analog to the
combinatorics of the Jack polynomials. Our main result is thus the
following theorem, to be proved in section 4.2, as  part of
Corollary \ref{routes}:

\begin{theorem}
 There exists a basis $\{ \tilde J_{\La}\}_{\La}$ of $\mathscr{P}^{S_\infty}$
such that \begin{equation}\begin{array}{lll} 1)&\tilde J_{\Lambda} =
m_{\Lambda} +\sum_{\Om < \La}
\tilde c_{\La \Om }(\beta) m_{\La}&\mbox{(triangularity)};\\
2)&\LL \tilde J_{\La}| \tilde J_{\Om} \RR_\beta \propto
\delta_{\La,\Om}&\mbox{(orthogonality)}.\end{array}\end{equation}
\end{theorem}

Quite remarkably, the
resulting construction is
completely equivalent to  that of the Jack polynomials $J_\La$ in
superspace (or Jack superpolynomials)  that were defined in
\cite{DLM3} from a  physical
eigenvalue problem. Before we establish this, it is probably appropriate to first review
some relevant aspects of
that article.

\subsection{Characterizations of the physical Jack superpolynomials}

We give the main properties of Jack superpolynomials.  The
results in the first part of this subsection can all be found in \cite{DLM3}.
The section is completed with the presentation of two technical lemmas.  All the results of this
section are independent of those of section 3.

First, we define a scalar product in $\mathscr{P}$, the algebra of
superpolynomials. Given
\begin{equation}
\Delta(x) = \prod_{1\leq j<k \leq N} \left[ \frac{x_j-x_k}{x_j x_k}
\right] \, ,
\end{equation}
$\langle\cdot|\cdot\rangle_{\beta,N}$ is defined (for $\beta$ a
positive integer) on the basis elements of $\mathscr{P}$ as
\begin{equation} \label{scalarp}
\langle \, \theta_I x^{\lambda}|\theta_{J}x^{\mu} \,\rangle_{\beta,N}
= \left\{
\begin{array}{ll}
{\rm C.T.} \left[\Delta^{\beta}(\bar x) \Delta^{\beta}(x)
\bar{x}^{\mu}
x^{\lambda}\right] & {\rm if~} I=J\, , \\
0 & {\rm otherwise}\, .
\end{array}  \right.
\end{equation}
where  ${\bar x}_i= 1/x_i$, and  where ${\rm C.T.}[E]$ stands for
the constant term of the expression $E$. (This is another form of the
scalar product
(\ref{physca}). More precisely, the latter is the analytic deformation of the former for all values of $\beta$.) This gives our first characterization of the Jack
superpolynomials.
\begin{proposition}\label{charac1}
   There exists a unique basis $\{J_\Lambda\}_\La$ of $\mathscr{P}^{S_N}$
such that \begin{equation}\begin{array}{lll}
1)&J_{\Lambda} = m_{\Lambda} +\sum_{\Om < \La} c_{\La \Om }(\beta)
m_{\La}& \mbox{(triangularity)};\\
2)&\langle J_{\La}| J_{\Om}\rangle_{\beta,N} \propto
\delta_{\La,\Om}&\mbox{(orthogonality)}.\end{array}\end{equation}
\end{proposition}

In order to present the other characterizations, we need to introduce
the Dunkl-Cherednik operators (for instance, see \cite{Baker}):
\begin{equation}\label{Dunk}{\mc{D}}_j:=x_j\partial_{x_j}+\beta\sum_{k<j}{\mc{O}}_{jk}+\beta\sum_{k>j}{\mc{O}}_{jk}-\beta(j-1)\,
,\end{equation} where
\begin{equation}{\mc{O}}_{jk}=\left\{\begin{array}{ll}
\ds\frac{x_j}{x_j-x_k}(1-{K}_{jk})\,,&k<j\, ,\\
\ds\frac{x_k}{x_j-x_k}(1-{K}_{jk})\,,&k>j\, .
\end{array}\right.\end{equation}
(recall that $K_{jk}$ is the operator that exchanges the variables
$x_j$ and $x_k$). These operators can be used to define two families
of operators that preserve the space $\mathscr{P}^{S_N}_{(n|m)}$:
\begin{equation}\label{defHI} \mathcal{H}_r:=\sum_{j=1}^N {{\mc{D}}_j^{\,
r}}\,\aand \mathcal{I}_s:=\frac{1}{(N-1)!}\sum_{\sigma\in
S_N}\mc{K}_\sigma \Bigl( \theta_1\partial_{\theta_1}{\mc{D}}_1^{\,
s} \Bigr)\mc{K}_{\sigma}^{-1} \, ,\end{equation} for $r\in
\{1,2,3,\dots,N \}$ and $s\in \{0,1,2,\dots,N-1\}$ (recall this time
that $\mathcal K_{\sigma}$ is built out of the operators
  $\mathcal K_{jk}$ that exchange $x_j \leftrightarrow x_k$
and  $\theta_j \leftrightarrow \theta_k$ simultaneously). These
operators are mutually commuting when restricted to
$\mathscr{P}^{S_N}$, that is \begin{equation}[ \mathcal{H}_r,
\mathcal{H}_s]f= [ \mathcal{H}_r, \mathcal{I}_s]f= [ \mathcal{I}_r,
\mathcal{I}_s]f=0
  \qquad \forall \, r,s \, ,
\end{equation} where $f$ represents an arbitrary polynomial in
$\mathscr{P}^{S_N}$. Since they are also symmetric with respect to
the scalar product $\langle\cdot|\cdot\rangle_{\beta}$ and have,
when considered as a whole, a non-degenerate spectrum, they provide
our second characterization of the Jack superpolynomials.
\begin{proposition} \label{charac2}
  The Jack superpolynomials $\{J_{\La}\}_{\La}$
are the unique common eigenfunctions of the $2N$ operators
$\mathcal{H}_r$ and $\mathcal{I}_s$, for $r\in \{1,2,3,\dots,N \}$
and $s\in \{0,1,2,\dots,N-1\}$.
\end{proposition}

We shall now define two operators that play a special role in our study.
\begin{equation}\label{defHIs}
\mathcal{H}:=\mathcal{H}_2+ \beta (N-1)\mathcal{H}_1-\mathrm{cst}
\aand \mathcal{I}:=\mathcal{I}_1\, ,\end{equation} where
$\mathrm{cst}=\beta N(1-3N-2N^2)/6$.  When acting on symmetric
polynomials in superspace, the explicit form of $\mc{H}$ is simply
\begin{equation}
\mc{H}=\sum_i (x_i \partial_{x_i})^2+\beta
\sum_{i<j}\frac{x_i+x_j}{x_{i}-x_j}(x_i
\partial_{x_i}-x_j\partial_{x_j})-2\beta
\sum_{i<j}\frac{x_ix_j}{(x_i-x_j)^2}(1-\kappa_{ij})\, .
\end{equation}

The operator $\mc{H}$ is the Hamiltonian of the stCMS model (see
Section 1.4.2); it can be written in terms of two fermionic
operators $\mathcal{Q}$ and $\mathcal{Q}^\dagger$ as
\begin{equation}
\mathcal{H}= \{\mathcal{Q}, \mathcal{Q}^\dagger\}\, ,\end{equation}
where
\begin{equation} \mathcal{Q}:=\sum_i\theta_ix_i\partial_{x_i} \aand
\mc{Q}^\dagger=\sum_i\partial_{\theta_i}\left(x_i\partial_{x_i}
+\beta \sum_{j\neq i}\frac{x_i+x_j}{x_i-x_j}\right)\,
,\end{equation} so that $\mc{Q}^2= ( \mc{Q}^\dagger)^2=0$.
Physically, $\mc{Q}$ is seen as creating fermions while
$\mc{Q}^\dagger$ annihilates them.
 A state (superfunction) which is annihilated by the fermionic
operators is called supersymmetric.  In the case of
superpolynomials, the only supersymmetric state is the identity.

\begin{remark} The stCMS Hamiltonian $\mc{H}$
has  an elegant differential geometric interpretation  as a
Laplace-Beltrami operator.  To understand this assertion,
consider  first the real Euclidian space $\mathbb{T}^N$, where
$\mathbb{T}=[0,2\pi)$.  Then, set $x_j=e^{\mathrm{i}t_j}$ for
$t_j\in\mathbb{T}$, and
 identify the Grassmannian  variable
$\theta_i$ with the differential form $dt_i$.  This allows us to
rewrite the physical scalar product \eqref{physca} as a Hodge-de~Rham
product involving complex differential forms, that is,
\begin{equation} \langle
A(t,\theta)|B(t,\theta)\rangle_{\beta,N}\sim
\int_{\mathbb{T}^N}\overline{A(t,dt)}\wedge \ast B(t,dt)\,
,\end{equation}where the bar denotes de complex conjugation and
where  the Hodge duality operator $\ast$ is formally defined by
\begin{equation}
A(t,dt)\wedge \ast
B(t,dt)=C_{\beta,N}\prod_{i<j}\sin^{2\beta}\left(\frac{t_i-t_j}{2}\right)\sum_k
\sum_{i_1<\ldots<i_k}
A_{i_1,\ldots,i_k}B_{i_1,\ldots,i_k}dt_1\wedge\cdots\wedge dt_N\,
,\end{equation}for some constant $C_{\beta,N}$.  Note that, in the
last equation, the forms $A$ and $B$ are developed in a way similar
to that of Eq.~\eqref{developform}.  Hence, we find that the fermionic
operators $\mc{Q}$ and $\mc{Q}^\dagger$  can be respectively interpreted
as the exterior derivative and its dual: $\mc{Q}\sim-\mathrm{i}d$ and
$\mc{Q}^\dagger\sim \mathrm{i}d^\ast$. Thus
\begin{equation}
\mc{H}=\Delta:=d\,d^\ast+d^\ast\,d\, .\end{equation} In consequence,
the Jack superpolynomials can be viewed as symmetric, homogeneous,
and orthogonal  eigenforms of a Laplace-Beltrami operator. This illustrates
the known connection between supersymmetric quantum mechanics and
 differential geometry \cite{Witten,Grandjean}. \end{remark}

  If the triangularity of the Jack superpolynomial
$J_\Lambda$ with respect to the supermonomial basis is
imposed, requiring that it be a
common eigenfunction
of  $\mathcal{H}$ and $\mathcal{I}$ is sufficient to define it. This
is our third
characterization of the Jack superpolynomials.
\begin{theorem}\cite{DLM3}  \label{TheoDefJack}
The Jack superpolynomials $\{J_\Lambda\}_\La$ form the unique basis
of $\mathscr{P}^{S_N}(\beta)$ such that \begin{equation}
\mc{H}(\beta)\,J_\La=\varepsilon_\La(\beta)\, J_\La\,,\quad
\mc{I}(\beta)\,J_\La=\epsilon_\La(\beta)\, J_\La\aand
J_\La=m_\La+\sum_{\Om<\La}c_{\La\Om}(\beta)m_\Om\,.\end{equation}
The eigenvalues are given explicitly by
\begin{eqnarray}
\varepsilon_\La(\beta)&=&\sum_{j=1}^N[
(\Lambda_j^*)^2+\beta(N+1-2j)\Lambda_j^*]\cr
\epsilon_\La(\beta)&=&\sum_{i=1}^m[\Lambda_i-\beta m(m-1)-\beta
\#_\La]\, ,\end{eqnarray}
  where $\#_\La$ denotes the
number of pairs $(i,j)$ such that $\La_i<\La_j$ for $1\leq i\leq m$
and $m+1\leq j\leq N$.
\end{theorem}

When no Grassmannian variables are involved, that is when
$\overline{\underline{\La}}=0$, our characterizations of the Jack
superpolynomials specialize to known characterizations of the Jack
polynomials that can be found for instance in \cite{Stan}. There is
however in the Jack polynomials' case a more common characterization
in which the scalar product appearing in Proposition~\ref{charac1}
is replaced by the strictly combinatorial scalar product (\ref{defscalprodcomb}).
As already announced,
this
more combinatorial characterization can be extended to the supersymmetric case.
But before turning to the analysis of the behavior of $J_\La$ with respect to the combinatorial scalar product, we present two lemmas concerning properties
of the eigenvalues $\varepsilon_\La(\beta)$ and $\epsilon_\La(\beta)$.

\begin{lemma}\label{lemmaeigenvalues1} Let $\La\in\mathrm{SPar}(n|m)$
and write $\la=\La^*$.  Let also
$\varepsilon_\La(\beta)$ and $\epsilon_\La(\beta)$ be the
eigenvalues given in Theorem~\ref{TheoDefJack}. Then
\begin{eqnarray} \varepsilon_\La(\beta)&=&2\sum_j
j(\lambda'_j-\beta
\la_j)+\beta n(N+1)-n\, ,\cr
\epsilon_\La(\beta)&=&|\La^a|-\beta|{\La'}^a|-\beta\frac{m(m-1)}{2}\,
.
  \end{eqnarray}
  \end{lemma}
  \begin{proof}
Using the well known identity (\cite{Mac}, Eq. 1.6),
  \begin{equation}\sum_j(j-1)\la_j=\sum_j\left(\begin{array}{c}
                                   \la'_j \\
                                   2 \\
                                 \end{array}\right)\, ,\end{equation}
we obtain
\begin{equation}\sum_j\la_j^2=\sum_j\la_j(\la_j-1)+\sum_j\la_j=2\sum_j\left(\begin{array}{c}
                                   {\la}_j \\
                                   2 \\
                                 \end{array}\right)+n=2\sum_j
j\la_j'-n\,.\end{equation} Hence, we have
\begin{equation}\sum_j[\la_j^2+\beta(N+1-2j)\la_j]=2 \sum_j
j(\la'_j-\beta\la_j)+\beta n(N+1)-n\, ,\end{equation} as desired. As
for the second formula, we consider
\begin{equation}\#_\La=\sum_{i=1}^m \#_{\La_i}\, ,\end{equation}
where $\#_{\La_i}$ denotes the number of parts in $\La^s$ bigger
than $\La_i$.  But from the definition of the conjugation, we easily
find that
\begin{equation}\#_{\La_i}={\La'}_{m+1-i}+1-i\, ,\end{equation}
so that
\begin{equation}\#_{\La}=\sum_{i=1}^m
({\La'}_i+1-i)=|{\La'}^a|+\frac{m(m-1)}{2} \, ,
\end{equation}
from which the second formula follows.
\end{proof}

The following lemma is of interest by
 itself but it will also be used to establish the orthogonality of the Jack superpolynomials in the case where the superpartitions can be compared (cf. the discussion following Proposition \ref{propositionHI}).

\begin{lemma}\label{lemmaeigenvalues2}Let $\La$ and $\Om$ be two
superpartitions of $(n|m)$.
Then
\begin{equation}\La >_S \Om\quad\Longrightarrow\quad
\varepsilon_\La\neq\varepsilon_\Om\aand \La >_T
\Om\quad\Longrightarrow\quad \epsilon_\La\neq\epsilon_\Om\,
.\end{equation}
\end{lemma}
\begin{proof}
First, let  $\La >_S \Om$, with $\om=\Om^*$ and $\la=\La^*$.  Then,
suppose that $\om=S_{ij}\la$ for some $i<j$.  By conjugation, this
implies that $\la'=S_{i'j'}\om'$ for some $i'<j'$.  More explicitly,
we have
\begin{equation}\om_i=\la_i-1\, \geq\, \om_j=\la_j+1 \quad {\rm and} \quad
\la'_{i'}=\om'_{i'}-1\, \geq\, \la'_{j'}=\om'_{j'}+1 \, .
\end{equation} Since
\begin{equation}
\sum_k k(\om_k-\la_k)= i(\la_i-1)+j(\la_j+1)- i\la_i-j\la_j=j-i\,
,\end{equation} it follows from the previous lemma that
\begin{equation}
\varepsilon_{\La}-\varepsilon_{\Om}= 2(j'-i')+2\beta (j-i) \, .
\end{equation}
Since by supposition, we have $i<j$ and $i'<j'$, the difference is a first
order polynomial in $\beta$ with positive coefficients.  From
Lemma~\ref{lemmaST}, when $\La >_S \Om$, we know that $\om$ can be
obtained by successive applications of such $S_{ij}$'s on $\la$.
Therefore, when $\La >_S \Om$, we have that
$\varepsilon_{\La}-\varepsilon_{\Om}$ is in general a first order
polynomial in $\beta$ with positive coefficients, and thus
$\varepsilon_{\La}\neq\varepsilon_{\Om}$.

For the second case, let $\La >_T \Om$ be such that $\Om =
T_{ij}\Lambda$ for some $i\in\{1,\ldots,m\}$ and $j\in\{m+1,\ldots,
N\}$.  By conjugation, this implies that  $\La' = T_{i'j'}\Om'$ for
some $i'\in\{1,\ldots,m\}$ and $j'\in\{m+1,\ldots, N\}$, and we thus
have
\begin{equation}\Om_i=\La_j\, <\, \Om_j=\La_i
\aand {\La'}_{i'}={\Om'}_{j'}\, <\, {\La'}_{j'}={\Om'}_{i'} \,
,\end{equation} which imply
\begin{equation}|\Om^a|<|\La^a| \aand |{\La'}^a|<|{\Om'}^a| \, .\end{equation}
Therefore, from the preceding lemma, $\epsilon_{\La}-\epsilon_{\Om}$
is a first order polynomial in $\beta$ with positive coefficients.
But again, when $\La >_T \Om$, Lemma~\ref{lemmaST} assures us that
$\Om$ can be obtained by successive applications of such $T_{ij}$'s
on $\La$. This means that $\epsilon_{\La}-\epsilon_{\Om}$ is in this
case a first order polynomial in $\beta$ with positive coefficients
and we have  $\epsilon_{\La} \neq \epsilon_{\Om}$ as claimed.
\end{proof}

\subsection{Orthogonality of the  Jack superpolynomials}

In terms of the combinatorial scalar product (\ref{defscalprodcomb}),  we can directly check the self-adjointness
of our eigenvalue-problem defining operators, $\mc{H}$ and $\mc{I}$.

\begin{proposition}\label{propositionHI} The operators $\mc{H}$ and $\mc{I}$
defined in (\ref{defHIs}) are self-adjoint (symmetric) with respect
to the combinatorial scalar product
$\LL \, \cdot \, | \, \cdot \,\RR_{\beta}$ defined in
(\ref{defscalprodcomb}).  \end{proposition}
\begin{proof}
We first rewrite $\mc{H}$ and $\mc{I}$ in terms of power sums. Since
these differential operators are both of order two, it is sufficient
to determine their action on the products of the form $p_m p_n$,
$\tilde{p}_m p_n $ and $\tilde{p}_m \tilde{p}_n $.  Direct
computations give
   \begin{eqnarray}
   %
   \mc{H}&=&\sum_{n\geq 1}[
n^2+\beta n(N-n)](p_n\, \partial_{p_n}+\tilde{p}_n\,
\partial_{\tilde{p}_n})}+\beta \sum_{n,m\geq
1}[(m+n)p_m\, p_n\, \partial_{p_{m+n}}+2m\,p_n\tilde{p}_m\,
\partial_{\tilde{p}_{n+m} ] \cr
&&\quad+\sum_{n,m\geq1}m n [p_{m+n}\,
\partial_{p_n}\, \partial_{p_m}+2
\tilde{p}_{n+m}\,\partial_{\tilde{p}_m} \partial_{p_n}\,]
\end{eqnarray}
  and
    \begin{eqnarray}
    \mc{I}=\sum_{n\geq
0}(1-\beta)(n \tilde{p}_n\,
\partial_{\tilde{p}_n})+\! \frac{\beta}{2}\sum_{m,n\geq 0}\tilde{p}_m\,
\tilde{p}_n\, \partial_{\tilde{p}_m}\,
\partial_{\tilde{p}_n}+\! \! \! \sum_{m\geq 0, n\geq 1}[ n\, \tilde{p}_{m+n}\,
\partial_{\tilde{p}_m}\,
\partial_{p_n}+\beta p_n\, \tilde{p}_m\,\partial_{\tilde{p}_{m+n}}]
\, .
\end{eqnarray}
  Note that these equations are valid when $N$ is
either infinite or finite.  In the latter case, the sums over the
terms containing $\tilde{p}_m$ and $p_n$ are respectively restricted
such that $m\leq N-1$ and $n\leq N$.

Then, letting $A^\perp$ denote the adjoint of a generic operator $A$
with respect to the combinatorial scalar product
(\ref{defscalprodcomb}), it is easy to check that
\begin{equation}\beta\,p_n^\perp=n\partial_{p_n} \aand
\beta\,\tilde{p}_n^\perp=\partial_{\tilde{p}_n}\,.\end{equation}
Hence, comparing the three previous equations, we obtain that
$\mc{H}^\perp=\mc{H}$ and $\mc{I}^\perp=\mc{I}$. For these
calculations, we recall that $(ab)^\perp=
b^\perp a^\perp$ even when
$a$ and $b$ are both fermionic.
  \end{proof}

The eigenvalue problem solved in Theorem~\ref{TheoDefJack}, together
with Proposition~\ref{propositionHI} and
Lemma~\ref{lemmaeigenvalues2} readily imply the orthogonality
of the Jack superpolynomials with respect to the combinatorial
scalar product in the special case where the two superpartitions can be compared
 with respect to the Bruhat order on compositions.

In order to extend this conclusion to all superpartitions, comparable or not, the
most natural path consists in  establishing the self-adjointness of all the
operators $\mc{H}_n$ and $\mc{I}_n$. But proceeding as for $\mc{H}$ and $\mc{I}$
above, by trying to reexpress them in terms of $p_n$, ${\tilde p}_n$ and their
derivatives, seems hopeless. An indirect line of attack is mandatory.

Let us first recall that the conserved operators (\ref{defHI}) can all be expressed in terms of
the Dunkl-Cherednik operators defined in (\ref{Dunk}). The ${\mc{D}}_i$ all commute among themselves:
\begin{equation} [{\mc{D}}_i, {\mc{D}}_j]=0\;.
\end{equation}
They are not quite invariant however, as\footnote{This corrects a misprint in eq. (25) of \cite{DLM3}.}
\begin{equation}\label{commu}
{\mc{D}}_i K_{i,i+1}-K_{i,i+1}\, {\mc{D}}_{i+1}  =
\beta\,.\end{equation} We will also need the following commutation
relations:
\begin{equation}\label{comA}[{\mc{D}}_i,x_i] = x_i+\beta
\left(\sum_{j<i}x_iK_{ij}+\sum_{j>i}x_jK_{ij}\right)\;,\end{equation}
while if $i\not= k$,
\begin{equation}\label{comB}
[{\mc{D}}_i,x_k]= -\beta x_{{\rm max} (i,k)}\,
K_{ik}\;.\end{equation}

The idea of the proof for the orthogonality is the following: in a first step, we show that the conserved
operators ${\mc{H}}_n$ and ${\mc{I}}_n$ are self-adjoint with respect to the combinatorial scalar
product and then  we demonstrate that this implies the orthogonality of the
$J_\La$'s.  The self-adjointness property is established via the kernel: showing that
$F= F^\perp$ is the same as showing that
\begin{equation}F^{(x)}K^{\beta,N}= F^{(y)} K^{\beta,N}\; ,
\end{equation}
where  $K^{\beta,N}$ is defined in Theorem \ref{cauchyppbeta}, and
where $F^{(x)}$ (resp. $F^{(y)}$) stands for the quantity $F$ in the
variable $x$ (resp. $y$). In order to prove this for  our conserved
operators  ${\mc{H}}_n$ and ${\mc{I}}_n$,  we need to establish some
results on the action of symmetric monomials in the  Dunkl-Cherednik
operators  acting on the following expression:
\begin{equation}\label{defO}
{\tilde\Omega}=
\prod_{i=1}^N\frac{1}{(1-x_iy_i)}\prod_{i,j=1}^N\frac{1}{(1-x_iy_j)^\beta}\; , \end{equation}
as well as some modification of ${\tilde\Omega}$. For that
matter, we recall a result of Sahi \cite{Sahi}:
\begin{proposition}\label{sahi} The action of the Dunkl-Cherednik operators ${\mc{D}}_j$ on
${\tilde\Omega}$ defined by (\ref{defO}) satisfies:
\begin{equation}
{\mc{D}}_j^{(x)}{\tilde\Omega}= {\mc{D}}_j^{(y)}{\tilde\Omega}\;.
\end{equation}
\end{proposition}

Before turning to the core of our argument, we establish the
following lemma.
\begin{lemma}\label{lemnouv}
Given a set $J= \{j_1,\dots, j_\ell \}$, denote by $x_J$ the product
$x_{j_1} \dots x_{j_{\ell}}$.  Suppose $x_J=K_{\sigma} x_I$ for some
$\sigma \in S_N$ such that $K_{\sigma} F K_{\sigma^{-1}}=F$.  Then
\begin{equation}
\frac{1}{ x_I} F^{(x)}  x_I \, \tilde \Omega = \frac{1}{ y_I}
F^{(y)}  y_I \, \tilde \Omega \quad \implies \quad \frac{1}{ x_J}
F^{(x)}  x_J \, \tilde \Omega = \frac{1}{ y_J} F^{(y)}  y_J \,
\tilde \Omega \, .
\end{equation}
\end{lemma}
\begin{proof}  The proof is straightforward and only uses the simple property
$K_{\sigma}^{(x)} \tilde \Omega = K_{\sigma^{-1}}^{(y)} \tilde
\Omega$.  To be more precise,
\begin{eqnarray}
\frac{1}{x_J} F^{(x)} x_J \, \tilde \Omega &=&
K_{\sigma}^{(x)} \frac{1}{x_I} F^{(x)} x_I K_{\sigma^{-1}}^{(x)} \, \tilde \Omega  \nonumber=K_{\sigma}^{(y)} K_{\sigma}^{(x)} \frac{1}{x_I} F^{(x)} x_I  \, \tilde \Omega  \nonumber =K_{\sigma}^{(y)} K_{\sigma}^{(x)} \frac{1}{y_I} F^{(y)} y_I  \, \tilde \Omega \nonumber \\
&=&
K_{\sigma}^{(y)}  \frac{1}{y_I} F^{(y)} y_I   K_{\sigma^{-1}}^{(y)} \, \tilde \Omega  = \frac{1}{y_J} F^{(y)} y_J \, \tilde \Omega \, .
\end{eqnarray}
\end{proof}

We are now ready to attack the main proposition.

\begin{proposition}\label{autoco}
The mutually commuting  operators ${\mc{H}}_n$ and ${\mc{I}}_n$
satisfy
\begin{equation}\label{autocha}{\mc{H}}_n^{(x)}K^{\beta,N}= {\mc{H}}_n^{(y)} K^{\beta,N}
\qquad {\rm and}\qquad  {\mc{I}}_n^{(x,\theta)}K^{\beta,N}= {\mc{I}}_n^{(y,\phi)}
K^{\beta,N}\;,\end{equation}
with $K^{\beta,N}$ the restriction to $N$ variables of the
kernel $K^\beta$ such as defined in Theorem \ref{cauchyppbeta}.
\end{proposition}

\begin{proof}  We first expand the kernel as follows:
\begin{eqnarray}
 K^{\beta,N} &=& K_0 \prod_{i,j}\left(1+\beta \frac{\theta_i \phi_j}{(1-x_i y_j)}\right) \\
&=&
 K_0\left\{ 1+\beta e_1 \left(
\frac{\ta_i\phi_j}{(1-x_iy_j)} \right) + \cdots + \beta^N e_N
\left(\frac{\ta_i\phi_j}{(1-x_iy_j)}\right)
 \right\}
\end{eqnarray}
where $K_0$ stands for $K^{\beta,N}(x,y,0,0)$, {\it i.e.},
\begin{equation}
 K_0:=
\prod_{i,j=1}^N\frac{1}{(1-x_iy_j)^\beta}\; ,
\end{equation}
and where $ e_{\ell} (u_{i,j}) $ is the elementary symmetric
function $e_\ell$ in the variables $u_{i,j}$  where
\begin{equation}
u_{i,j}:= \frac{\ta_i\phi_j} {(1-x_iy_j)} \qquad  i,j=1,\dots,N\; .
\end{equation}
Note that, in these variables,
the maximal possible elementary symmetric function is $e_N$ given that
$\theta_i^2=\phi_i^2=0$. In the following, we will use the compact
notation $I^-=\{1,\cdots, i-1\}$ and $I^+=\{i,\cdots, N\}$ (and
similarly for $J^\pm$), together with  $w_{I^-}= w_1\cdots w_{i-1}$
and $w_{I^+}= w_i\cdots w_N$.

The action of the operators on $K^\beta$ can thus be decomposed into
their action on each monomial in this expansion.  Now observe that
$K_0$ is invariant under the exchange of any two variables $x$ or
any two variables $y$. Therefore, if an operator $F$ is such that
${\mathcal K}_{\sigma} F {\mathcal K}_{\sigma}^{-1}=F$ for all
$\sigma \in S_N$, and such that
\begin{equation} \label{eqaprouver}
F^{(x,\theta)} v_{I^-} \, K_0 = F^{(y,\phi)} v_{I^-} \, K_0\; \qquad
{\rm with}\qquad v_i:=u_{i,i}\end{equation} for all $i=1,\dots,N+1$,
then we immediately have by symmetry that $F^{(x,\theta)}
K^{\beta} =  F^{(y,\phi)} K^{\beta}$.  We will use this
observation in the case of $\mathcal H_n$ and $\mathcal I_n$.

We first consider the case $F=\mathcal H_n$. Recall from (\ref{defHI}) that
${\mc{H}}_n =p_n({\mc{D}}_i) $
is such that ${\mathcal K}_{\sigma} {\mathcal H}_n {\mathcal
K}_{\sigma}^{-1}={\mathcal H}_n$ (see \cite{DLM3}).  Since $\mathcal
H_n$ does not depend on the fermionic variables, we thus have to
prove from the previous observation that
\begin{equation}
{\mathcal H}_n^{(x)}  \frac{1}{(1-xy)_{I^-} }\, K_0 = {\mathcal
H}_n^{(y)} \frac{1}{(1-xy)_{I^-} } \, K_0\;, \end{equation} or
equivalently
\begin{equation}
{\mathcal H}_n^{(x)}  \,(1-xy)_{I^+}\, \,{\tilde\Omega} = {\mathcal
H}_n^{(y)} \,(1-xy)_{I^+}\, {\tilde\Omega}\;,
\end{equation}
for all $i=1,\dots,N+1$ (the case $i=N+1$ corresponds to the empty
product).

The underlying symmetry allows us to further simplify the problem by
focusing on the terms
\begin{equation} y_{J^+} \, {\mathcal H}_n^{(x)} \, x_{J^+} , {\tilde\Omega} =
x_{J^+} \, {\mathcal H}_n^{(y)} \, y_{J^+}  \,{\tilde\Omega}\;,
\end{equation} for $j\geq i$, or equivalently, on
\begin{equation}\frac{1}{x_{J^+}  } {\mathcal H}_n^{(x)}\, {x_{J^+} }\,{\tilde\Omega}=
\frac{1}{y_{J^+}  } {\mathcal H}_n^{(y)}\, {y_{J^+}
}\,{\tilde\Omega}\;.
\end{equation}
This follows from Lemma~\ref{lemnouv} which assures us
 that all the different terms can be obtained
from these special ones.

Now, instead of analyzing the family ${\mathcal
H}_n=p_n({\mc{D}}_i)$, it will prove simpler to  consider the
equivalent  family $e_n({\mc{D}}_i)$.  We will first show the case
$e_N({\mc{D}}_i)$, that is,
\begin{equation}\frac{1}{x_{J^+} } {\mc{D}}^{(x)}_1\cdots {\mc{D}}^{(x)}_N\, x_{J^+}
\, {\tilde\Omega} =\frac{1}{y_{J^+} }{\mc{D}}^{(y)}_1\cdots
{\mc{D}}^{(y)}_N \, y_{J^+} \, {\tilde\Omega}\;.\end{equation} Let
us concentrate on the left hand side. We note that
\begin{equation}
 \frac{1}{x_{J^+} } {\mc{D}}^{(x)}_1\cdots {\mc{D}}^{(x)}_N \, x_{J^+} \,
{\tilde{\Omega}} = \frac{1}{x_{J^+} } {\mc{D}}^{(x)}_1{x_{J^+}
}\cdots\frac{1}{x_{J^+} } {\mc{D}}^{(x)}_N\, x_{J^+}
{\tilde{\Omega}}\;. \end{equation} It thus suffices to
study each term $({x_{J^+} })^{-1} {\mc{D}}_j{x_{J^+} }$ separately.
In each case we find that
\begin{equation} {\mc{D}}_k\, {x_{J^+} }=
 {x_{J^+} } \, {\tilde {\mc{D}}}_k\;.\end{equation} The
form of ${\tilde {\mc{D}}}$ depends upon $j$ and $k$. There are two
cases:
\begin{eqnarray}\label{castr}
 k<j: &&{\tilde {\mc{D}}}_k=
 {\mc{D}}_k - \beta\sum_{\ell=j}^N K_{\ell,k}\;,\cr
k\geq j:&& {\tilde {\mc{D}}}_k= {\mc{D}}_k+ 1+
\beta\sum_{\ell=1}^{j-1} K_{\ell,k}\;
\end{eqnarray}
which can be easily checked using (\ref{comA}) and (\ref{comB}).
We can thus write
\begin{equation}\frac{1}{x_{J^+} } {\mc{D}}^{(x)}_1\cdots {\mc{D}}^{(x)}_N
\, x_{J^+}  {\tilde\Omega}
 = {\tilde {\mc{D}}}^{(x)}_1\cdots {\tilde {\mc{D}}}^{(x)}_N\,{\tilde\Omega}\;.\end{equation}
Using proposition \ref{sahi} and $K_{ij}^{(x)} \tilde \Omega =
K_{ij}^{(y)} \tilde \Omega$, the rightmost term ${\tilde
{\mc{D}}}^{(x)}_N$ can thus be changed into ${\tilde
{\mc{D}}}^{(y)}_N$. Since it commutes with the previous terms (i.e.,
it acts on the variables $y$ while the others act on $x$), we have
\begin{eqnarray}
 {\tilde {\mc{D}}}^{(x)}_1\cdots {\tilde {\mc{D}}}^{(x)}_{N-1}{\tilde
{\mc{D}}}^{(y)}_N \,{\tilde\Omega} & =&{\tilde {\mc{D}}}^{(y)}_N
{\tilde {\mc{D}}}^{(x)}_1\cdots {\tilde {\mc{D}}}^{(x)}_{N-1}
\,{\tilde\Omega}=  {\tilde {\mc{D}}}^{(y)}_N {\tilde
{\mc{D}}}^{(y)}_{N-1}\cdots {\tilde {\mc{D}}}^{(y)}_{1}
\,{\tilde\Omega}\cr &=& \frac{1}{y_{J^+} }{\mathcal D}^{(y)}_N
y_{J^+} \cdots \frac{1}{y_{J^+} }{\mathcal D}^{(y)}_1 \, y_{J^+}
{\tilde\Omega}\cr\ &=& \frac{1}{y_{J^+} }{\mathcal D}^{(y)}_N\cdots
{\mathcal D}^{(y)}_1 \, y_{J^+} {\tilde\Omega} =
\frac{1}{y_{J^+} }{\mathcal D}^{(y)}_1\cdots {\mathcal D}^{(y)}_N \,
y_{J^+} {\tilde\Omega}\;,\end{eqnarray} which is the desired result.

At this point, we have only considered a single conserved operator,
namely $e_N({\mathcal D}_i)$. But by replacing ${\mc{D}}_i$ with
${\mc{D}}_i+t$ in $e_N({\mathcal D}_i)$, we obtain a generating
function  for all the operators $e_n({\mathcal D}_i)$.  Since to
prove $e_N({\mathcal D}_i^{(x)}+t) K^{\beta,N}= e_N({\mathcal
D}_i^{(y)}+t) K^{\beta,N}$ simply amounts to replacing ${\tilde
{\mc{D}}}_i$ by ${\tilde {{\mc{D}}}}_i+t$ in the previous argument,
we have completed the proof of ${\mc{H}}_n^{(x)}K^{\beta,N}=
{\mc{H}}_n^{(y)} K^{\beta,N}$.

For the case of ${\mc{I}}_n$, we start with the expression given in (\ref{defHI})
which readily implies that  ${\mathcal K}_\sigma \,
{\mathcal I}_n \, {\mathcal K}_\sigma^{-1}= {\mathcal I}_n$.
Therefore, from the observation surrounding formula
(\ref{eqaprouver}), and because the derivative $\ta_1 \d_{\ta_1}$
annihilates the $K_0$ term in the expansion of $K^{\beta,N}$, we only
need to show that
\begin{equation}
{\mathcal I}_n^{(x,\ta)} v_{I^-}\, K_0 = {\mathcal I}_n^{(y,\phi)}
 v_{I^-} \,  K_0\;, \end{equation}
for $i=2,\dots,N+1$. Up to an overall multiplicative factor, the
only contributing part in ${\mc{I}}_n$, when acting on $ v_{I^-} $, is
\begin{equation}
{\mc O}_n:= {\mc{D}}_1^n+K_{12}{\mc{D}}_1^nK_{12} + \cdots
K_{1,i-1}{\mc{D}}_1^nK_{1,i-1}\;.\end{equation} It thus suffices to show
that
\begin{eqnarray}
{\mc O}_n^{(x)}\, (1-xy)_{I^+}\, {\tilde\Omega}
 =  {\mc O}_n^{(y)}\, (1-xy)_{I^+} \, {\tilde\Omega}
\end{eqnarray}
Once more, we can use Lemma~\ref{lemnouv} since ${\mc O}_n$ commutes
with $K_{k,\ell}$ for $k,\ell\geq i$. Thus, we only need to check that
for $j\geq i$,
\begin{eqnarray}
\frac{1}{x_{J^+} } {\mc O}_n^{(x)}{x_{J^+} }\,{\tilde\Omega}
=\frac{1}{y_{J^+} } {\mc O}_n^{(y)}{y_{J^+} }\,{\tilde\Omega}\;.
\end{eqnarray}
Since the $K_{1\ell}$'s act trivially on the  variables $x_j$ for
$j> \ell$, the previous relation reduces to proving
\begin{equation}
\frac{1}{x_{J^+} } [{\mc{D}}_1^n]^{(x)} \, {x_{J^+}
}\,{\tilde\Omega} = \frac{1}{y_{J^+} }[{\mc{D}}_1^n]^{(y)}\,
{y_{J^+}}\,{\tilde\Omega}\;.\end{equation} The left hand side takes
the form
\begin{equation}
\frac{1}{x_{J^+} } [{\mc{D}}_1^n]^{(x)} \, {x_{J^+}
}\,{\tilde\Omega}= \left\{ \frac{1}{x_{J^+} }{\mc{D}}_1^{(x)}\,
{x_{J^+} }\right\}^n \,{\tilde\Omega}\;.\end{equation} We then only
have to evaluate $({x_{J^+}})^{-1}{\mc{D}}_1^{(x)}\, {x_{J^+}}$. The
result is given by the first case in (\ref{castr}) (since $j>1$) .
The proof is completed as follows
\begin{eqnarray}
\left\{ \frac{1}{x_{J^+}}{\mc{D}}_1^{(x)}\, {x_{J^+}}\right\}^n
\,{\tilde\Omega}=  [{\tilde {\mc{D}}}^{(x)}_1]^n \,{\tilde\Omega} =
  [{\tilde {\mc{D}}}^{(y)}_1]^n \,{\tilde\Omega} = \frac{1}{y_{J^+} }
[{\mc{D}}_1^n]^{(y)}\,  {y_{J^+}
}\,{\tilde\Omega}\;.\end{eqnarray}\end{proof}

As previously mentioned, the proposition has the following
corollary.
\begin{corollary}\label{autocos}
The operators $\mc{H}_r$ and $\mc{I}_s$
defined in (\ref{defHI}) are self-adjoint (symmetric) with respect
to the combinatorial scalar product
$\LL \, \cdot \, | \, \cdot \, \RR_{\beta,N}$
given in (\ref{defscalprodcombN}).  \end{corollary}
This immediately gives our main result.
\begin{theorem}\label{mainortho} The Jack superpolynomials $\{J_{\La} \}_{\La}$ are orthogonal
with respect to the combinatorial scalar product, that is,
\begin{equation}
\LL J_{\La} | J_{\Om} \RR_{\beta} \propto \delta_{\La,\Om}\;.
\end{equation}
\end{theorem}

\begin{proof} The fact that in $N$ variables
$\LL J_{\La} | J_{\Om} \RR_{\beta,N} \propto \delta_{\La,\Om}$
follows from the equivalence between the statement of the theorem
and Corollary~\ref{autocos}. This equivalence
  follows
 from Proposition~\ref{charac2}, which says that the Jack polynomials are the
unique common eigenfunctions of the $2N$ operators appearing in
Corollary~\ref{autocos}.  Given that the
expansion coefficients of the Jack superpolynomials in terms of
supermonomials do not depend on the number of variables $N$ \cite{DLM3},
the theorem then follows from (\ref{degN}).
\end{proof}

\begin{corollary} \label{routes} The following statements are direct consequences of the orthogonality
property of the Jack polynomials in superspace.\footnote{In fact, it
can be shown that all statements of Corollary~\ref{routes} and
Theorem~\ref{theoduality} below are not only consequences of
Proposition~\ref{mainortho} but are equivalent to it.}
\begin{enumerate}
\item  There exists a basis $\{ \tilde J_{\La}\}_{\La}$ of $\mathscr{P}^{S_\infty}$
such that \begin{equation}\begin{array}{lll} 1)&\tilde J_{\Lambda} =
m_{\Lambda} +\sum_{\Om < \La}
\tilde c_{\La \Om }(\beta) m_{\La}&\mbox{(triangularity)};\\
2)&\LL \tilde J_{\La}|\tilde J_{\Om} \RR_\beta \propto
\delta_{\La,\Om}&\mbox{(orthogonality)}.\end{array}\end{equation}
\item Let $K^\beta$ be the reproducing kernel defined in
Proposition~\ref{cauchyppbeta}.  Then,
\begin{equation}K^\beta(x,\theta;y,\phi)=
\sum_{\La\in\mathrm{SPar}}j_{\La}(\beta)^{-1}\versg{J_\Lambda(x,\theta)}\,
\versd{J_\Lambda(y,\phi)}\,.\end{equation} where
\begin{equation}\label{normde}
j_\La(\beta):=\LL
\versg{J_{\La}} | \versd{J_{\La}} \RR_{\beta}\;.
\end{equation}
\item Let $\{g_{\Lambda}\}_{\Lambda}$ be the basis, defined in (\ref{defgla}),
dual to that of the supermonomials
with respect to the combinatorial scalar product.
Then, the Jack superpolynomials expand upper  triangularly in
this basis:
\begin{equation}
J_{\La} =\sum_{\Om \geq \La} u_{\La \Om}(\beta) \, g_{\Om} \, ,
\quad {\rm with} \quad u_{\La \La}(\beta) \neq 0 \, .
\end{equation}
\end{enumerate}
\end{corollary}
\begin{proof}
{\it 1}.  It was shown in \cite{DLM3} that the operators $\mathcal H$ and
$\mathcal I$ act triangularly on the supermonomial basis.  Thus,
$\mathcal H$ and $\mathcal I$ also act triangularly on the basis $\{
\tilde J_{\Lambda} \}_\La$. Furthermore, from
Proposition~\ref{propositionHI}, they are self-adjoint with respect
to the combinatorial scalar product.  Hence, we must conclude from
the previous argument that  $\tilde J_{\La}$ is an eigenfunction of
$\mathcal H$ and $\mathcal I$, from which Theorem~\ref{TheoDefJack}
implies that $\tilde J_{\La} =J_{\La}$.

\bigskip

\n {\it 2}. The proof is similar to that of  Lemma~\ref{cauchybases} (see also Section VI.2 of
\cite{Mac}).

\bigskip

\n {\it 3}.  Suppose that $\LL J_{\La}| J_{\Om} \RR_\beta \propto
\delta_{\La,\Om}$, and let $J_{\La}= \sum_{\Om \in \mathcal S}
u_{\La \Om} \, g_{\Om}$, where $\mathcal S$ is some undefined set.
If $\La$ is not the smallest element of $\mathcal S$, then there
exists at least one element $\Gamma$ of $\mathcal S$ that does not
dominate any other of its elements. In this case, we have
\begin{equation}
\LL J_{\La}| J_{\Gamma} \RR_\beta= \sum_{\Om \in \mathcal S} u_{\La
\Om}(\beta) \sum_{\Delta \leq \Gamma} c_{\Gamma \Delta}(\beta) \,
\LL g_{\Om}| m_{\Delta} \RR_\beta \, .
\end{equation}
Since $\Gamma$ does not dominate any element of $\mathcal S$, the
unique non-zero contribution in this expression is that of
  $u_{\La \Gamma}(\beta) \, c_{\Gamma \Gamma}(\beta)
\,  \LL g_{\Gamma}, m_{\Gamma} \RR_\beta= u_{\La \Gamma}(\beta)$.
Since this term is non-zero by supposition, we have the
contradiction $0=\LL J_{\La}| J_{\Gamma} \RR_\beta=u_{\La
\Gamma}(\beta)\neq 0$.
\end{proof}

\subsection{Duality}

In this subsection, we show that the homomorphism
$\hat{\omega}_\beta$, defined in Eq.~\ref{defendobeta}, has a simple
action on Jack superpolynomials.   To avoid any confusion, we make
explicit the $\beta$ dependence of the Jack superpolynomials by
writing $J^{(1/\beta)}_\La$.\footnote{The rationale for this notation is to match the one used in \cite{Mac} when $m=0$: $J^{(1/\beta)}_{\La}(x,\ta)= J^{(1/\beta)}_{\La^s}(x)= J^{(\alpha)}_{\La^s}(x)$, where
$\alpha=1/\beta$. (Similarly, in our previous works
\cite{DLM1,DLM3}, we denoted $J^{(1/\beta)}_\La$ by
$J_\La(x,\theta;1/\beta)$ to
keep our definition similar to the usual form introduced by Stanley
\cite{Stan} as $J_\la(x;\alpha)$ when $m=0$). We stress however, that when we need to make explicit the $\beta$-dependence of $j_\La,\, \mc{H}$ and $\mc{I}$, we write $j_\La(\beta)\,,\mc{H}(\beta)$and $\mc{I}(\beta)$ respectively.}

\begin{proposition}\label{propomega}
  We have
\begin{equation}\mc{H}(\beta)\hat{\omega}_{\beta}J_\La^{(\beta)}=\varepsilon_{\La'}(\beta)\,
\hat{\omega}_{\beta}J_\La^{(\beta)}\aand
\mc{I}(\beta)\hat{\omega}_{\beta}J_\La^{(\beta)}=\epsilon_{\La'}(\beta)\,
\hat{\omega}_{\beta}J_\La^{(\beta)}\, .\end{equation}
\end{proposition}
\begin{proof}Let us rewrite the special form of the operator
$\mc{H}(\beta)$ appearing in the proof of Proposition
\ref{propositionHI} as
\begin{equation}\mc{H}(\beta)=\sum_{n\geq1}[n^2+\beta
n(N-n)]\hat{A}_n+\sum_{m,n\geq 1}(\beta
\hat{B}_{m,n}+\hat{C}_{m,n})\, ,\end{equation} with
\begin{equation}\begin{array}{lll}
  \hat{A}_n&=&p_n\,
\partial_{p_n}+\tilde{p}_n\,
\partial_{\tilde{p}_n}\, ,\\
\hat{B}_{m,n}&=&(m+n)p_m\, p_n\,
\partial_{p_{m+n}}+2m\,p_n\tilde{p}_m\,
\partial_{\tilde{p}_{n+m}}\, ,\\
\hat{C}_{m,n}&=&m n  \left( p_{m+n}\,
\partial_{p_n}\, \partial_{p_m}+2 \tilde{p}_{n+m}\,\partial_{\tilde{p}_m}
\partial_{p_n} \right)\, .
\end{array}\end{equation}From these definitions, we get
\begin{equation}\hat{\omega}_{1/\beta}\hat{A}_n=\hat{A}_n \,
\hat{\omega}_{1/\beta}\, ,
\quad
\hat{\omega}_{1/\beta}\hat{B}_{m,n}=-\frac{1}{\beta}\hat{B}_{m,n} \,
\hat{\omega}_{1/\beta} \aand
\hat{\omega}_{1/\beta}\hat{C}_{m,n}=-{\beta}\hat{C}_{m,n} \,
\hat{\omega}_{1/\beta} \, .\end{equation} These relations imply
\begin{eqnarray}
\hat{\omega}_{1/\beta}\mc{H}(\beta)\hat{\omega}_{\beta}=\sum_{n\geq1}[n^2+\beta
n(N-n)]\hat{A}_n-\! \! \! \! \sum_{m,n\geq 1}(
\hat{B}_{m,n}+\beta\hat{C}_{m,n})
=(1+\beta)N\sum_{n\geq1}n\hat{A}_n-\beta\mc{H}(1/\beta)\, . \nonumber
\end{eqnarray}Now,
  considering  $\sum_{n\geq1}n\hat{A}_nm_\La=|\La|m_\La$ and Lemma
\ref{lemmaeigenvalues1}, we obtain
\begin{equation}\hat{\omega}_{1/\beta}\mc{H}(\beta)
\hat{\omega}_{\beta}J_\La^{(\beta)}=\varepsilon_{\La'}(\beta)J_\La^{(\beta)}\,
\end{equation} as claimed. The relation involving $\mc{I}(\beta)$ is
proved in a similar way.
\end{proof}

\begin{theorem}\label{theoduality} The homomorphism $\hat{\omega}_{\beta}$ is such that
\begin{equation}\hat{\omega}_{1/\beta}\versd{J^{(1/\beta)}_\La}=j_\La(\beta)\versg{J^{(\beta)}_{\La'}}\,
\end{equation}
with $ j_\La(\beta)$ such as defined in (\ref{normde}).
\end{theorem}
\begin{proof}
Let us first prove that $\hat{\omega}_{\beta}J_\La^{(\beta)}\propto J_{\La'}^{(1/\beta)}$.
 From the third point of Corollary~\ref{routes},  we know that $J_{\La}^{(1/\beta)}
=\sum_{\Om \geq \La} u_{\La \Om}(\beta) \, g_{\Om}$. But
Eq.~(\ref{dualityge}) implies $\hat{\omega}_{1/\beta}g_\La=e_\La$.
Hence,
\begin{equation}
\hat \omega_{1/{\beta}} \Bigl(J_\La^{(1/\beta)}\Bigr) = \sum_{\Om \geq
\La} u_{\La \Om}(\beta) \, e_{\Om} = \sum_{\Om \geq \La} u_{\La
\Om}(\beta) \sum_{\Gamma \leq \Om'} N_{\Om}^{\Gamma} \,  \versg{m_{\Gamma}}=
  \sum_{\Gamma \leq \La'} v_{\La \Gamma}(\beta) \, \versg{m_{\Gamma} }\, ,
\end{equation}
where we have used (\ref{dualityge}), Theorem~\ref{theoebase} and
the fact that $\Om \geq \La \iff \Om' \leq \La'$. Further, since
$N_{\Lambda}^{\Lambda'}=1$ and $u_{\La \La}(\beta)\neq 0$, we have
$v_{\La \La'} \neq 0$. Now, from Proposition~\ref{propomega},
$\hat{\omega}_{{1/\beta}} \Bigl(J_\La^{(1/\beta)}\Bigr)$ is an
eigenfunction of $\mathcal H(1/\beta)$ and $\mathcal I(1/\beta)$
with eigenvalues $\varepsilon_{\La'}(1/\beta)$ and
$\epsilon_{\La'}(1/\beta)$ respectively.   The triangularity we just
obtained ensures from Theorem~\ref{TheoDefJack}, that
$\hat{\omega}_{{1/\beta}} \Bigl(J_\La^{(1/\beta)}\Bigr)$ is
proportional to ${J_{\La'}^{(\beta)}}$.

Again from Theorem~\ref{theoebase}, we know that
$m_\La=(-1)^{m(m-1)/2}e_{\La'}+\mbox{higher terms}$, so that
\begin{equation}
J^{(1/\beta)}_\La=(-1)^{m(m-1)/2}e_{\La'}+\mbox{higher terms}\,
.\end{equation} Moreover,  from Eq.~\eqref{dualityge}, we get
\begin{equation} \hat{\omega}_\beta J_\La^{(1/\beta)}=(-1)^{m(m-1)/2}g_{\La'}+\mbox{higher terms}\, .
\end{equation}
But  the proportionality proved above  implies
\begin{equation}
\hat{\omega}_{1/\beta}\versd{J^{(1/\beta)}_\La} =A_\La(\beta)\, \versg{J_{\La'}^{(\beta)} }=A_\La(\beta)\, \versg{m_{\La'}}Ê+\mbox{lower
terms}\, ,
\end{equation} for some constant $A_\La(\beta)$.  Finally,
considering the duality between $g_\La$ and $m_\La$, we obtain
\begin{eqnarray}
(-1)^{m(m-1)/2}j_\La(\beta)&=&\LL  J_\La^{(1/\beta)} |
J_\La^{(1/\beta)}\RR_\beta \nonumber\\
&=& \LL \hat{\omega}_{\beta}
J_\La^{(1/\beta)} | \hat{\omega}_{1/\beta} J_\La^{(1/\beta)}\RR_\beta \nonumber\\
&=& \LL (-1)^{m(m-1)/2} \versd{g_{\La'} }| A_\La(\beta)\versg{m_{\La'}}\RR_\beta \nonumber\\
&=& (-1)^{m(m-1)/2} A_\La(\beta)
\end{eqnarray} as desired.
\end{proof}

\subsection{Limiting cases}

In Section 4.2, we have proved that the physical Jack
superpolynomials are orthogonal with respect to the combinatorial
scalar product. This provides a direct link with the material of
Section 3. Other links,  less general but more explicit, are
presented in this section, from the consideration of
$J_\La$ for special values of $\beta$ or for particular
superpartitions.

\begin{proposition}  \label{propogjack} For $\La= (n)$ or $(n;0)$, we have
(using the notation of Proposition~\ref{propgn}):
 \begin{equation}J_{(n)}=\frac{n!}{(\beta+n-1)_n}g_n
\aand J_{(n;0)}=\frac{n!}{(\beta+n)_{n+1}}\tilde{g}_n\,
.\end{equation}
\end{proposition}
\begin{proof}
Since $J_{(0;1^n)}=m_{(0;1^n)}=\tilde{e}_n$, we have on the one hand
$\hat{\omega}_{\beta}(J_{(0;1^n)}) =\tilde g_n$ from
(\ref{dualityge}). On the other hand, from
Proposition~\ref{propomega}, $\hat{\omega}_{\beta}(J_{(0;1^n)})$ is
an eigenfunction of $\mathcal H(\beta)$ and $\mathcal I(\beta)$ with
eigenvalues $\varepsilon_{(n;0)}(\beta)$ and
$\epsilon_{(n;0)}(\beta)$ respectively. Since $(n;0)$ is the highest
partition with one fermion in the Bruhat order, we have from Theorem
\ref{TheoDefJack}, that there exists a unique eigenfunction of
$\mathcal H$ and $\mathcal I$ with such eigenvalues. We must thus
conclude that $\tilde g_n$ is also proportional to
$J_{(n;0)}$. Looking at Proposition \ref{propgn}
and considering that the coefficient of $m_{(n;0)}$ in $J_{(n;0)}$
needs to be equal to one, we obtain
$(\beta+n)_{n+1}J_{(n;0)}={n!}\,\tilde{g}_n$.  The relation between
$J_{(n)}$ and $g_n$ is proved in a similar way.
\end{proof}

\begin{corollary}\label{coronorme} For $\La= (n)$ or $(n;0)$, the combinatorial norm of $J_\La$ is
\begin{equation} \LL J_{(n)}|J_{(n)}\RR_\beta=\frac{n!}{(\beta+n-1)_n}\aand
\LL J_{(n;0)}|J_{(n;0)}\RR_\beta=\frac{n!}{(\beta+n)_{n+1}}\,
.\end{equation}
\end{corollary}
\begin{proof}Using the previous proposition, we get
\begin{align}
(n!)^2\,\LL
g_n|g_n\RR_\beta=&(\beta+n-1)_n^{2}\;\LL J_{(n)}|J_{(n)}\RR_\beta\,
,\cr
 (n!)^2\,\LL
\tilde{g}_n|\tilde{g}_n\RR_\beta=&(\beta+n)_{n+1}^{2}\;\LL
J_{(n;0)}|J_{(n;0)}\RR_\beta\, .
 \end{align}
From
Proposition~\ref{propgn}, we know that
\begin{equation}n!\,g_n=(\beta+n-1)_n \; m_{(n)}+\ldots \, ,\quad
n!\,\tilde{g}_n=(\beta+n)_{n+1}\; m_{(n;0)}+\ldots \, ,\end{equation}
where the dots stand for lower terms in the Bruhat ordering.  Thus,
considering Corollary \ref{corodualitygm}, we get
\begin{equation} \LL g_n|g_n\RR_\beta=\frac{(\beta+n-1)_n}{n!}\,
,\quad \LL
\tilde{g}_n|\tilde{g}_n\RR_\beta=\frac{(\beta+n)_{n+1}}{n!}\end{equation}
and the proof follows.
\end{proof}

\begin{theorem}\label{limdebeta}
For $\beta=0,\,1,$ or $\beta\rightarrow \infty$, the limiting expressions of $J_\La^{(1/\beta)}$ are
\begin{equation}J_\La^{(1/\beta)}\longrightarrow\left\{\begin{array}{ll}
m_\La\, & {\rm when ~}\beta\longrightarrow0\, ,\\
\versg{e_{\La'} } \, &{\rm when ~} \beta\longrightarrow\infty\,
.\end{array}\right.\end{equation}
and
\begin{equation}J_{(n)}^{(1)}=h_n\aand
J_{(n;0)}^{(1)}=\frac{1}{n+1}\tilde{h}_n\, .\end{equation}
\end{theorem}
\begin{proof}The case $\beta \to 0$
is a direct consequence of Theorem \ref{TheoDefJack}, given that
$\mathcal H(\beta)$ and $\mathcal I(\beta)$ act diagonally on
supermonomials in this limit. The second case is also obtained from
the eigenvalue problem. Indeed, when $\beta\rightarrow\infty$,
  $\beta^{-1}\mc{H}(\beta)$ and  $\beta^{-1}\mc{I}(\beta)$ behave as
first order differential operators.
  Then, it is easy to get
  \begin{equation}\bigg[\lim_{\beta\rightarrow\infty}\frac{\mc{H}(\beta)}{\beta}\bigg]\,
  e_{\La'}=\bigg[-2\sum_j j\la_j+n(N-1)\bigg]\,e_{\La'}\wwhere
\la=\La^*\end{equation}and
  \begin{equation}\bigg[\lim_{\beta\rightarrow\infty}\frac{\mc{I}(\beta)}{\beta}\bigg]\,
  e_{\La'}=\bigg[-|\La^a|-\frac{m(m-1)}{2}\bigg]\,e_{\La'}\, .\end{equation}
These are the eigenvalues of $J_{\La}$ in the limit where $\beta\rightarrow \infty$ (cf. Lemma
\ref{lemmaeigenvalues1}). The proportionality constant between  $e_{\La'}$ and $J_\La$ is fixed by   Theorem  \ref{theoebase} and Theorem \ref{TheoDefJack}.  We have thus
 \begin{equation} \versg{e_{\La'} }=
  \lim_{\beta\rightarrow\infty} \versd{J_\La^{(1/\beta)} }\, .\end{equation}
  Finally, we note that the property concerning $h_n$ and
$\tilde{h}_n$ is an immediate
  corollary of Proposition~\ref{propogjack}.
  \end{proof}

\subsection{Normalization}

In this subsection,  $\tilde m_{\Lambda}$ shall denote the augmented supermonomial:
\begin{equation}
\tilde m_{\Lambda} = n_{\Lambda}! \, m_{\Lambda} \, ,
\end{equation}
where $n_{\Lambda}!$ is such as defined in (\ref{augmented}).

It is easy to see that the smallest
superpartition of degree $(n|m)$ in the Bruhat ordering is
\begin{equation}\label{minide}
\Lambda_{\mathrm{min}}:=(\delta_m\,;\,
1^{\ell_{n,m}}\,) \, ,\end{equation} where
\begin{equation}\label{delnm}
\ell_{n,m}:=n-|\delta_m|\, ,\quad \delta_m:=(m-1,m-2,\ldots,0) \aand
|\delta_m|=\frac{m(m-1)}{2}\, .\end{equation} Now, let
$c_\La^{\mathrm{min}}(\beta)$ stand for the coefficient of
$\tilde m_{\Lambda_{\mathrm{min}}}$ in the monomial expansion of
$J_\La^{(1/\beta)}$. We shall establish a relation between this
coefficient and the norm of the Jack superpolynomials $J_\La$.


\begin{proposition}The norm $j_\La(\beta)$  defined  in (\ref{normde}), with $\La\vdash(n|m)$, is
\begin{equation}
j_\La(\beta)= \beta^{-m-\ell_{n,m}}
 \, \frac{c_\La^{\mathrm{min}}(\beta)}{c_{\La'}^{\mathrm{min}}(1/\beta)} \end{equation}
\end{proposition}
\begin{proof}
One readily shows that
\begin{equation}
m_{\La_\mathrm{min}}=
p_{\La_\mathrm{min}}+\mbox{higher terms}\,.\end{equation}
Since
$m_{\La_\mathrm{min}}$ is the only supermonomial containing
$p_{\La_\mathrm{min}}$, we can write
\begin{equation}\label{pmini}
J_{\La}^{(1/\beta)}=c_\La^{\mathrm{min}}(\beta) \, p_{\La_\mathrm{min}}+\mbox{higher
terms}\, .\end{equation}
Let us now apply $\hat{\omega}_{1/\beta}$ on this expression. Using  Eq.~\eqref{involup} we get
\begin{equation} \hat{\omega}_{1/\beta}{J^{(1/\beta)}_\La}=
\beta^{-m-\ell_{n,m}}(-1)^{m(m-1)/2}\, {c_\La^{\mathrm{min}}(\beta)}\,
p_{\La_\mathrm{min}}+\mbox{higher terms}\, .\end{equation}
But if we apply $\hat{\omega}_{1/\beta}$ on $J_\La^{(1/\beta)}$ by using first
Theorem~\ref{theoduality} to write it as $(-1)^{m(m-1)/2 }\,  j_\La(\beta)J_{\La'}^{(\beta) }$ and expand $J_{\La'}^{(\beta)}$ using (\ref{pmini}), we get instead
\begin{equation}\hat{\omega}_{1/\beta}{J^{(1/\beta)}_\La}=j_\La(\beta)(-1)^{m(m-1)/2}\, {c_{\La'}^{\mathrm{min}}(1/\beta)}\, p_{\La_\mathrm{min}}+\mbox{higher
terms}\;.
\end{equation}
Here we have used the fact that $\Lambda_{\mathrm{min}}$, being the smallest
superpartition of degree $(n|m)$ in the Bruhat ordering, labels the smallest supermonomial in both the decomposition of $J_\La$ and $J_{\La'}$. The result follows from the comparison of the last two equations.
\end{proof}

The coefficient $c_\La^{\mathrm{min}}(\beta)$ appears from computer
experimentation to have a very simple form. We shall now introduce
the notation needed to describe it. Recall that $D[\Lambda]$ is the
diagram  used to represent $\Lambda$. Given a cell $s$
in $D[\Lambda]$, let $a_{\La}(s)$ be the number of cells (including
the possible circle at the end of the row) to the right of $s$.  Let
also $\ell_{\La}(s)$ be the number of cells (not including the
possible circle at the bottom of the column) below $s$.  Finally,
let $\Lambda^{\circ}$, be the set of cells of $D[\La]$ that do not
appear at the same time in a row containing a circle {\it and} in a
column containing a circle.
\begin{conjecture}  The coefficient  $c_\La^{\mathrm{min}}(\beta)$ of
$\tilde m_{\Lambda_{\mathrm{min}}}$
in the monomial expansion of $J_\La^{(1/\beta)}$ is given by
\begin{equation}
c_\La^{\mathrm{min}}(\beta) = \frac{1}{\prod_{s \in
\Lambda^{\circ}} \Bigl( a_{\Lambda}(s)/\beta + \ell_{\Lambda}(s)+1
\Bigr)  }
\end{equation}
with $ {\La_\mathrm{min}}$ and $\ell_{n,m} $ defined in (\ref{minide}) and (\ref{delnm}) respectively.
\end{conjecture}
For instance, if $\Lambda=(3,1,0;4,2,1)$, we can fill $D[\La]$ with
the values $\bigl( a_{\Lambda}(s)/\beta + \ell_{\Lambda}(s)+1
\bigr)$ corresponding to the cells $s \in \La^{\circ}$.  This gives
(using $\gamma=1/\beta$):
\begin{equation} {\small{\tableau[mcY]{{\mbox{\tiny
$3\gamma+5$}}&{\mbox{\tiny $2\gamma+3$}} &{\mbox{\tiny $\gamma+2$}}&
{\mbox{\tiny $1$}} \\& &{\mbox{\tiny
$\gamma+1$}}&\bl\gcercle\\{\mbox{\tiny $\gamma+3$}}&{\mbox{\tiny
$1$}}\\ & \bl\gcercle\\{\mbox{\tiny $1$}}\\  \bl \gcercle}} }
\end{equation}
Therefore, in this case,
\begin{equation}
c_\La^{\mathrm{min}}(\beta) =
\frac{1}{(3/\beta+5)(2/\beta+3)(1/\beta+2)(1/\beta+1)(1/\beta+3)}
\end{equation}

Even though the Jack superpolynomials cannot be normalized to have
positive coefficients when expanded in terms of supermonomials, we
nevertheless conjecture they satisfy the following {\it integrality} property.
\begin{conjecture} Let
\begin{equation}
J_\Lambda^{(1/\beta)}= c_\La^{\mathrm{min}}(\beta)\sum_{\Omega \leq
\Lambda} \tilde{c}_{\Lambda \Omega}(\beta) \, \tilde m_{\Omega}.
\end{equation}
Then $\tilde{c}_{\Lambda \Omega}$ is a polynomial in $1/\beta$ with
integral coefficients.
\end{conjecture}

\section{Conclusion}

\subsection{Summary}

In this work, we have presented an extension of the theory of
symmetric functions
involving fermionic variables as well as the usual bosonic variables. Our
construction being motivated by
supersymmetric considerations, we enforce from the beginning an equal
number of variables of
each type. These variables  can thus be regarded as the coordinates
of an Euclidian superspace.
Symmetric functions in superspace, or equivalently, superfunctions,
are defined to be symmetric
with respect to the diagonal action of the symmetric group.

Basically all essential objects in the theory of symmetric functions have been
extended to superspace.  If some of them had already been
introduced in previous works of ours (such as superpartitions and
supermonomials in \cite{DLM1,DLM2,DLM3} and power-sum
superpolynomials in \cite{DLM5}, section 2.5 \footnote{We have
noticed in the meantime that the elementary and
power-sum
superpolynomials had also been introduced in the second reference of
\cite{sCMS}.}), most of these extensions are new.

The resulting theory of symmetric functions in superspace, exposed
in Section 3, is quite elegant and appears to be rather rich. We
have also pointed out an interesting connection between
superpolynomials and de Rham complexes of symmetric $p$-forms.

The core results of the elementary theory of symmetric functions are
known to have a one
parameter (our $\beta$) deformation that leads to the combinatorial
definition of the Jack
polynomials.  This deformation also has a superspace lift that turns out
to be related to
our previous construction of the Jack
superpolynomials using an
approach rooted in the solution of a supersymmetric integrable quantum many-body problem
\cite{DLM1,DLM3}. In  special cases, namely, $\beta=0$,
$\beta\rightarrow\infty$ and for
$\La=(n;0)$ or $(n)$, the physical $J_\La$ were shown to reduce to
combinatorial symmetric
superfunctions (cf. Theorem
\ref{limdebeta} and Proposition
\ref{propogjack} respectively).

  The physical
Jack superpolynomials were already known to be orthogonal with
respect to the physical scalar product (\ref{physca}) or
(\ref{scalarp}) (the first one being the analytic extension of the second one, which holds when $\beta
$ is a positive integer). Here, by relying on the integrability of the physical underlying quantum
many-body problem, we have been able to prove that the physical Jack superpolynomials are also orthogonal
with respect to the combinatorial scalar product
(\ref{defscalprodcomb}).

At once, these two products are
manifestly very different from each other.  That
the Jack superpolynomials are orthogonal with respect to both
products is certainly remarkable.

Even in the absence of fermionic variables, the orthogonality of the
Jack polynomials with respect to both scalar products is a highly
non-trivial observation. In that case, one can provide a partial
rationale for the compatibility between the two scalar products, by
noticing their equivalence in the following two circumstances
\cite{Kadell, Mac}:
\begin{equation}\label{equivbeta1}
 \langle f  | g \rangle_{\beta=1,N}=\LL f  | g \RR_{\beta=1,N}
\end{equation}
(see e.g., \cite{Mac} VI.9 remark 2) and
\begin{equation}\label{equivNinfty}\lim_{N \rightarrow\infty}\frac{\langle f | g
\rangle_{\beta,N}}{\langle 1|1\rangle_{\beta,N}}= \LL f | g
\RR_\beta \, ,\end{equation} (see e.g., \cite{Mac} VI.9 (9.9)) for
$f,\, g,$ two arbitrary symmetric polynomials.

In superspace, this compatibility between the two products is even
more remarkable since the limiting-case equivalences
\eqref{equivbeta1} and \eqref{equivNinfty} are simply lost. This is
most easily seen by realizing that, after integration over the
fermionic variables, we obtain
 \begin{equation} \langle p_\la
\tilde{p}_n | p_\mu \tilde{p}_m \rangle_{\beta,ÊN}=\langle p_\la  | p_\mu
p_{m-n}\rangle_{\beta,ÊN}\, ,\quad m>n\, ,\end{equation}
and thus the super power-sums cannot be orthogonal for any value
of $N$ and
$\beta$. This shows that the connection between the two scalar
products is rather intricate.

\subsection{Outlook: Schur superpolynomials and supercombinatorics}

A central chapter in the  theory of symmetric functions concerns the
Schur polynomials $s_\la$. These are limiting case $\beta=1$ of the Jack
polynomials $J_\la^{(\beta)}$. Their
special importance lies in their deep representation theoretic
interpretation: $s_\la$ is a Lie-algebra character,  being
expressible as a sum of semistandard tableaux of shape $\la$. This
implies, in particular, that $s_\la$  has a monomial expansion with
non-negative integer coefficients.

Schur superpolynomials could similarly be defined  from the Jack
superpolynomials evaluated at
$\beta=1$:
\begin{equation}
s_\La(x,\ta)= J^{(1)}_\La(x,\ta)
\end{equation}
This is a quite natural guess.\footnote{Off hand, one could have contemplated a  connection
more intricate than $\beta=1$ between the would-be Schur
superpolynomials and the
Jack superpolynomials, such as a limiting value of
$\beta$ that depends upon the fermionic degree, for instance,
$\beta=m+1$.
 This would provide a stable definition with respect to the product of two
superpolynomials, under which the fermionic degree is additive.
But physically, it is not natural to have a coupling constant that depends upon the fermionic degree,  that is, upon a special property of the eigenfunctions.}
But is there any special property
that points
toward this identification of $J_\La^{(1)}$ with the superspace
generalization of the Schur polynomials? Unfortunately, we have not
been able to pinpoint a genuine distinguishable
  feature of the Jack superpolynomials that would single out the
special value
$\beta=1$.
For instance, the supermonomial decomposition is not integral, as shown by the following example:
\begin{equation}
s_{(2,1;0)}=m_{(2,1;0)}+\frac12 m_{(2,0;1)}-\frac18 m_{(1,0;2)}+\frac14 m_{(1,0;1,1)}\; .
\end{equation}
Similarly, the kernel $K$ can be expanded as
\begin{equation}K(x,\theta;y,\phi)=
\sum_{\La\in\mathrm{SPar}}j_{\La}(1)^{-1}\versg{s_\Lambda(x,\theta)}\,
\versd{s_\Lambda(y,\phi)}\,.\end{equation}
However, the coefficients $j_\La(1)$ are not equal to 1 as in the usual bosonic case.
The only property of the Schur polynomials that is readily  transposed to the superspace is the following one: the homogeneous superpolynomials decompose upward triangularly in terms of $s_\La$'s. Indeed, since $J_\La$ is upper triangular in the $g$- basis, the inverse is true (that is, $g_\La$ is upper triangular in the $J$-basis) and the $\beta=1$ version of this expansion is the announced property:
\begin{equation} h_\La(x,\ta)=\sum_{\Om\geq \La} v_{\La,\Om}s_\La(x,\ta)
 \end{equation}

For completeness, we mention that the proper superspace generalization of the classical definition of the Schur polynomials as a bialternant has not been found yet. (We point out in that regard that division by anticommuting variables is prohibited). Similarly, the Jacobi-Trudi identity, which expresses the Schur polynomials in terms of the $h_n$'s, has not been generalized.  Notice  that in all the instances where we have obtained a 
determinantal expression, we had at most one row or one column  made out of fermionic quantities,
something which cannot be the case for the sought for Jacobi-Trudi super-identity. 
To determine whether these properties are specific to the $m=0$ sector or not requires further study.

Note finally that, off-hand, it appears unlikely that the Schur
superpolynomials would be related to the representation theory of
special Lie superalgebras since these theories do not involve
Grassmannian variables.  Actually, it could  well be that for the
Schur superpolynomials, the representation theoretic interpretation
is simply lost.

In a different vein, with the introduction of superdiagrams, we
expect a large number of
results linked to ``Ferrer-diagram combinatorics'' to have
nontrivial extensions to the
supercase. Pieces of supercombinatorics have already been presented
at the end of Section 2.4. The
combinatorics of the diagrams $D[\La]$  also enters in the
formulation of the  norm of the Jack
superpolynomials. Another natural  ground for supercombinatorics would be to find
a Pieri formula for
the Jack superpolynomials.

\subsection{Outlook: Macdonald superpolynomials}

In this work, we have heavily stressed the existence of a
one-parameter (i.e., $\beta$) deformation
of the scalar product as the key tool for defining Jack
superpolynomials combinatorially. However,
there  also exists a two-parameter deformation ($t$ and $q$) of the
combinatorial scalar product. Again, this has a natural lift to the
superspace, namely
\begin{equation}  \label{qtspro}\LL \versg{p_\La}|
\versd{p_\Om}\RR_{q,t}:=z_\La(q,t)\delta_{\La,\Omega}\,
,\end{equation}
where
\begin{equation}
z_\La(q,t)=z_\La \prod_{i=1}^{m}\frac{1-q^{\La_i+1}}{1-t^{\La_i+1}}
\prod_{i=m+1}^{\ell(\La)}\frac{1-q^{\La_i}}{1-t^{\La_i}}\, ,\quad
m=\overline{\underline{\La}}\, .\end{equation}
This reduces to the previous scalar product $ \LL
\cdot| \cdot\RR_\beta$ when $q=t^{1/\beta}$ and $t\rightarrow 1$.
The generalized form of the reproducing kernel reads
\begin{equation}
\prod_{i,j}\frac{\left(tx_iy_j+t\theta_i\phi_j;q\right)_\infty}{\left(x_iy_j+\theta_i\phi_j;q\right)_\infty}
=\sum_\La
z_\La(q,t)^{-1}\versg{p_\La(x,\theta)}\versd{p_\La(y,\phi)}\,
,\end{equation}
with $(a;q)_\infty$ defined in (\ref{gefct}).

Now, the scalar product (\ref{qtspro}) leads directly to a
conjectured combinatorial definition of Macdonald
superpolynomials.

\begin{conjecture}
In the  space of symmetric superfunctions with rational coefficients
in
  $q$ and $t$, there exists a basis $\{  M_\La\}_{\La}$, where
  $M_\La=M_\La(x,\theta;q,t)$,
such that \begin{equation}\begin{array}{lll} 1)& M_{\Lambda} =
m_{\Lambda} +\sum_{\Om < \La}
  C_{\La \Om }(q,t) m_{\La}&\mbox{(triangularity)};\\
2)&\LL  \versg{M_{\La}}|\versd{ M_{\Om}} \RR_{q,t} \propto
\delta_{\La,\Om}&\mbox{(orthogonality)}.\end{array}\end{equation}
\end{conjecture}

Note that in this context, the combinatorial construction cannot be
compared with the physical one since the corresponding
supersymmetric eigenvalue problem has not been formulated yet. In
other words, the  proper supersymmetric version of the
Ruijsenaars-Schneider model \cite{Ruij} is still missing.

\begin{acknow}
We thank A. Joyal for pointing out the connection between ${\tilde
e}_{n-1}, \, {\tilde
h}_{n-1}$ and $ n{\tilde p}_{n-1}$ and the corresponding  one-forms  $de_n,\,
dh_{n} $ and $ dp_{n}$. This work was  supported by NSERC and
FONDECYT (Fondo Nacional de
Desarrollo Cient\'{\i}fico y Tecnol\'ogico) grant
\#1030114. P.D. is grateful to the Fondation J.A.-Vincent for a
student fellowship.
\end{acknow}

\end{document}